\documentclass[leqno,12pt]{article}
\usepackage[margin=1in]{geometry}

\usepackage{dutchcal}

\usepackage{lmodern,microtype,mathtools,enumitem,needspace}

\usepackage{amsthm,amsfonts,amsbsy,amssymb,amsmath,amscd}

\usepackage[all]{xy}
\usepackage{xr-hyper}
\usepackage{tikz-cd}
\usepackage[colorlinks=true,linkcolor=blue,urlcolor=blue,citecolor=blue]{hyperref}

\usepackage[T1]{fontenc}

\usepackage{verbatim}
\usepackage{lscape}

\usepackage[textsize=footnotesize]{todonotes}

\newcommand{\BC}{{\mathbb {C}}}

\newcommand{\BH}{{\mathbb {H}}}

\newcommand{\BP}{{\mathbb {P}}}
\newcommand{\BQ}{{\mathbb {Q}}}
\newcommand{\BR}{{\mathbb {R}}}

\newcommand{\CC  }{{\mathcal {C}}}

\newcommand{\CH}{{\mathcal {H}}}

\newcommand{\CO}{{\mathcal {O}}}

\newcommand{\an}{{\mathrm{an}}}

\newcommand{\rk}{{\operatorname{rk}}}

\newcommand{\Pic}{\operatorname{Pic}}

\renewcommand{\Re}{{\operatorname{Re}}}

\newcommand{\Spec}{{\operatorname{Spec}}}

\newcommand{\tr}{{\operatorname{Tr}}}

\newcommand{\vol}{{\operatorname{vol}}}

\newcommand{\adeg}{\wh{\deg}}

\newcommand{\wh}{\widehat}

\newtheorem{thm}{Theorem}[subsection]
\newtheorem{maintheorem}{Theorem}

\newtheorem{cor}[thm]{Corollary}
\newtheorem{lem}[thm]{Lemma}

\newtheorem {conj}[thm]{Conjecture}
\newtheorem{defn}[thm]{Definition}

\newtheorem{theorem}[thm]{Theorem}
\newtheorem{proposition}[thm]{Proposition}
\newtheorem{lemma}[thm]{Lemma}
\newtheorem{corollary}[thm]{Corollary}

\theoremstyle{definition}
\newtheorem{definition}[thm]{Definition}
\newtheorem{example}[thm]{Example}

\theoremstyle{remark}
\newtheorem{remark}[thm]{Remark}

\numberwithin{equation}{subsection}

\begin{document}

\title{Heights and energies of polarized fibrations of varieties I\\[3pt]\large Fibrations over curves}

\author{Shou-Wu Zhang}
\date{\today}

\maketitle

\begin{abstract}
Let $K$ be a number field or the function field of a curve, and
let $X\to C$ be a projective fibration over a curve over $K$.
We prove a height inequality for relatively nef adelic line
bundles $\bar L$ on $X$ with $L$ nef. For a nef, relatively ample
line bundle $L$ over a finitely generated field, we deduce a lower
bound for the local energy at almost every complex or nonarchimedean
realization of a fixed polarized model. Over $\mathbb C$, for an
ample line bundle $L$ and a smooth positive polarization, the
Euler--Lagrange equation is the Calabi--Yau equation. On the
projective bundle $X=\BP_C(E)$ of a vector bundle $E$ of degree zero, our energy on
quotient metrics on $L=\CO(1)$ agrees with Donaldson's functional when the
determinant metric is fixed. We propose a construction of canonical
heights on parameter spaces using envelopes of local energies.
\end{abstract}
\clearpage\begingroup\small\tableofcontents\endgroup\clearpage

\section*{Introduction}
\addcontentsline{toc}{section}{Introduction}
Let $K$ be a number field or the function field $k(T)$ of a
regular integral projective curve $T$ over a field $k$. Let
$f:X\to C$ be a flat projective surjective morphism from an integral
$K$-variety of dimension $d$ to a smooth projective curve over $K$.
We study the adelic self-intersection of a line
bundle on $X$ and how it varies with the metric. We fix a nef adelic
polarization $\bar A$ on $X$ with $A$ ample and $A^d=1$. In Part~II
we extend the height and energy inequalities to fibrations over
higher-dimensional polarized bases \cite{PartII}.

An adelic line bundle $\bar L$ on $X$ is \emph{strongly relatively
nef} if $\bar L+f^*\bar N$ is nef for some nef adelic bundle
$\bar N$ on $C$ with $N$ ample. It is \emph{relatively nef} if
there is a nef adelic $\mathbb Q$-line bundle $\bar M$ such that
$\bar L+\epsilon\bar M$ is strongly relatively nef for every
rational $\epsilon>0$.

\begin{maintheorem}\label{thm:main-height}
Let $\bar L$ be an integrable relatively nef adelic line
bundle whose underlying bundle $L$ is nef. Then
\[
 \bar L^{d+1}
 \le(d+1)L^d(\bar L\cdot\bar A^d).
\]
In particular, if $L^d=0$, then $\bar L^{d+1}\le0$.
\end{maintheorem}

Over a number field, the fundamental height inequalities imply
Theorem~\ref{thm:main-height} for $L$ nef and big with semipositive
metrics (Remark~\ref{rem:height-fundamental-inequality}).
For relatively nef metrics, we first prove the inequality on models.
Over a number field the proof uses a lower bound for analytic torsion
(Theorem~\ref{thm:hi}), and over $k(T)$ it uses Riemann--Roch on $T$
(Theorem~\ref{prop:function-field-height}). In both cases we pass to
adelic bundles by approximation in Theorem~\ref{thm:ahi}.

Over a number field, taking $X=C$, $f=\mathrm{id}_C$, and $\deg L=0$ gives the
Faltings--Hriljac inequality \cite{Fa84,Hr85}.
When $L_\eta^{d-1}>0$ and the normalized class $L^\circ$ is nef,
Theorem~\ref{thm:main-height} gives the discriminant inequality of
Theorem~\ref{thm:bgi}. For the
tautological quotient bundle on $\mathbb P(E)$ with $E$ semistable,
it recovers Moriwaki's arithmetic Bogomolov--Gieseker inequality
\cite{Mo94}. For related inequalities, see Miyaoka \cite{Mi88} and
Soul\'e \cite{So94}. Kawaguchi and Moriwaki proved a relative
Bogomolov inequality for arithmetic families of curves \cite[Theorem~8.1]{KM01}.

Proposition~\ref{prop:polarized-surface-blowups} gives examples on
blown-up surfaces in which $L$ is nef and relatively ample but not
semiample. The construction gives both $L^2=0$ and $L^2>0$.

We denote the difference between the right and left sides of
Theorem~\ref{thm:main-height} by
\[
 I_A(\bar L)=(d+1)L^d(\bar L\cdot\bar A^d)-\bar L^{d+1}.
\]
To describe how $I_A$ varies at a place, we first work over a field
$F$ with a real-valued absolute value and take the analytification over
its completion. We fix a semipositive metric on $A$. For two continuous
integrable metrics $h_0,h_1$ on $L$, put $\alpha_i=c_1(L,h_i)$ and
define the local energy by
\[
 \Psi_A(h_0,h_1)=\int_{X^{\mathrm{an}}}
   \log\frac{|\cdot|_{h_1}}{|\cdot|_{h_0}}
   \left(\sum_{i=0}^d\alpha_0^i\alpha_1^{d-i}
                    -(d+1)L^d c_1(\bar A)^d\right).
\]
Over a global field, we write $\Psi_v$ for this energy with the
weight of the place $v$ included. If $\bar L_0$ and $\bar L_1$ agree
away from $v$, the change-of-metric formula gives
$\Psi_v(h_{0,v},h_{1,v})=I_A(\bar L_1)-I_A(\bar L_0)$.
Thus $\Psi_v(h_{0,v},h_{1,v})\ge-I_A(\bar L_0)$ whenever
$\bar L_1$ satisfies Theorem~\ref{thm:main-height}. When $L^d=0$,
the energy is independent of $\bar A$.

For $L$ relatively ample, let $\CH^+(L)$ consist of smooth metrics
with strictly positive vertical curvature over $\mathbb C$, and
rational model metrics induced by relatively ample bundles in the
nonarchimedean case. By Lemma~\ref{lem:relative-nef-approximation},
the infimum of $\Psi_A(h_0,h)$ over $h\in\CH^+(L)$ equals its infimum
over the continuous integrable relatively nef metrics $h\in\BH^+(L)$.

To study local energies over a finitely generated field $F/K_0$,
we fix a local field $E$ containing $K_0$. For $K_0=\BQ$, we take
$E=\BC$ or a finite extension of $\BQ_p$. For $K_0=\mathbb F_p(t)$,
we take a finite extension of $\mathbb F_p((t))$ with the standard
embedding of $K_0$. Let $f:X\to C$ be a flat projective surjective
morphism of smooth geometrically integral projective varieties over
$F$, with $\dim X=d$ and $\dim C=1$.

Each $K_0$-embedding $\sigma:F\hookrightarrow E$ gives a fibration
$X_\sigma\to C_\sigma$ over $E$. To see these fibrations in one
family, let $K$ be the algebraic closure of $K_0$ in $F$ and spread
$f$ over a smooth geometrically integral $K$-variety $S$ with
function field $F$. For each
$\iota:K\hookrightarrow E$ extending the fixed embedding of $K_0$,
the embeddings $\sigma$ extending $\iota$ correspond exactly to
the $E$-points of $S_\iota$ lying over the generic point of $S$.
The fiber at such a point is $X_\sigma\to C_\sigma$. We fix an
admissible collection of polarizations $(\bar A_\sigma)_\sigma$,
induced by one metrized model of the family and normalized by
$A_\sigma^d=1$ (Definition~\ref{def:admissible-polarizations}).

\begin{maintheorem}\label{thm:main-energy}
Let $L$ be a nef, $f$-ample $\mathbb Q$-line bundle.
Then, for almost every $K_0$-embedding $\sigma:F\hookrightarrow E$
and every integrable reference metric $h_{0,\sigma}$,
\[
 \inf_{h\in\BH^+(L_\sigma)}
 \Psi_\sigma(h_{0,\sigma},h)>-\infty.
\]
\end{maintheorem}

The exceptional set may depend on the polarized model and is
independent of the metrics on $L_\sigma$, including the reference
metric. We use smooth volume on a complex parameter space and Haar
measure in nonarchimedean coordinates. Their null sets are
independent of the choice of smooth model.

To prove Theorem~\ref{thm:main-energy}, we use algebraic dynamics
on $\BP^n$ and finite maps to $\BP^n$ to construct zero-cycles on the
parameter space whose points have uniformly bounded height. Their
normalized counting measures converge to a measure dominating smooth
volume or the Haar measure class. The height inequality bounds the
average energy on these cycles. By gluing test metrics over the
parameter space, we obtain a lower bound for the integral of the
fiberwise infimum.

For $L$ nef and $f$-ample, Conjecture~\ref{conj:H-energy-boundedness}
predicts a finite lower bound for the energy on every such fibration
of smooth geometrically integral varieties over $\mathbb C$ or a complete discretely valued field. When $L^d=0$,
Conjecture~\ref{conj:sp} identifies the metrics attaining the infimum
with the globally semipositive metrics in $\BH^+(L)$. For a smooth
complex fibration, the critical equation is
$c_1(L,h)^d=L^d c_1(\bar A)^d$. When $L^d=0$ and $h\in\CH^+(L)$
is critical, $c_1(L,h)$ is semipositive of rank $d-1$
(Proposition~\ref{prop:minimizer}).

For a family of flat projective morphisms $X\xrightarrow{f}C\to S$
of relative dimensions $d-1$ and $1$, let $g:X\to S$ be the
composite and normalize $A_s^d=1$. The combination of Deligne pairings
\[
 D_A(L)=-\langle L^{[d+1]}\rangle_g
             +(d+1)L_s^d\langle L,A^{[d]}\rangle_g
\]
defines a rational line bundle on $S$, where $[j]$ denotes $j$
repeated entries. Fixing the metric on $A$, we write $D_A(h)$ for
the metric on $D_A(L)$ induced by $h$ and
$D_A(L,h)=(D_A(L),D_A(h))$. For an adelic bundle $\bar L=(L,h)$
we also write $D_A(\bar L)$. When $S=\Spec K$, we have
$I_A(\bar L)=\widehat{\deg}\,D_A(\bar L)$.

If $w_1,w_2$ are the weights of $D_A(h_1),D_A(h_2)$ in a common
local frame, their difference
\[
 \Psi_A(h_1,h_2)=w_2-w_1
\]
is the weight of the canonical section $1$ of
$D_A(L,h_2)\otimes D_A(L,h_1)^\vee$. By base change, its value at
$s$ is the energy $\Psi_{A_s}(h_{1,s},h_{2,s})$ on $X_s$.

Over a number field, when $L$ is $f$-ample and $L_s$ is nef on every geometric fiber,
we conjecture that the local envelopes define a nef adelic metric
on $D_A(L)$. Its height at $s\in S(\bar K)$ should be the infimum of
$h_{D_A(\bar L)}(s)$ over the Hermitian model bundles
$\bar L\in\CH^+(L)$ of Section~\ref{ss:family}
(Conjecture~\ref{conj:envelope}).

\paragraph*{Declaration on AI use.}
I used ChatGPT (OpenAI) to assist with language editing,
exposition, literature searches, and mathematical discussion during
the preparation of this paper. Mathematically, AI assistance was used for:
(a) adapting the analytic torsion methods of Bismut--Vasserot
\cite{BV90} from projective bundles to obtain the one-sided estimate
for fibrations $X\to C$ with $X$ and $C$ smooth and projective,
and locating relevant further references;
(b) the cutoff and gluing arguments in
Lemmas~\ref{lem:energy-envelope-gluing} and~\ref{lem:nonarch-energy-selection};
and (c) the examples obtained by blowing up points on surfaces.
I take full responsibility for the final text and its
mathematical content.

\phantomsection\label{par:notation-conventions}
\paragraph*{Notation and conventions.}
We write $L+M$ for $L\otimes M$, $-L$ for $L^\vee$, and
$L^{\otimes n}$ for the $n$th tensor power. Products such as $L^d$ and
$(\bar L\cdot\bar A^d)$ are intersection numbers. Over a number
field we use the unnormalized arithmetic degree. Over $K=k(T)$ we use
geometric adelic intersections, with the normalization of
Theorem~\ref{prop:function-field-height}. For $f:X\to C$,
$\eta$ denotes the generic point of $C$, so $L_\eta^{d-1}$ is the
top intersection number on the generic fiber.

For a flat projective morphism $g:Y\to S$ of pure relative
dimension $r$, we write $D_g(L)=D_{Y/S}(L)=\langle L^{[r+1]}\rangle_g$.
We use the same notation for metrized bundles. Our polarization comes
from an ample rational line bundle $A_0$ with a nef adelic metric
in the global setting and a semipositive metric in the local setting.
In a family with $\dim X_s=d$, we put
$\bar A=(A_{0,s}^d)^{-1/d}\bar A_0$. Flatness makes
$A_{0,s}^d$ constant on each connected component. The $d$ copies
of $A$ in $D_A(L)$ contribute the rational factor $1/A_{0,s}^d$.
We use $\bar A_0$ for rational perturbations. The subscripts in
$\Psi_v,\Psi_s,\Psi_\sigma$ indicate the local polarization.

On a local field $F$ we use the absolute value defined by
$\vol(a\Omega)=|a|_F\vol(\Omega)$ for additive Haar measure.
Thus $|\cdot|_{\BR}$ is the usual absolute value,
$|z|_{\BC}=|z|^2$, and $|\varpi|_F=(\#\kappa_F)^{-1}$ for a
uniformizer in the nonarchimedean case. A metric satisfies
$|au|_h=|a|_F|u|_h$. In a local frame $\ell$ its weight is
$\phi=-\log|\ell|_h$. We add metrics by adding their weights and
write $\bar L(\psi)$ for the metric $e^{-\psi}|\cdot|_{\bar L}$.
At a complex place $|\ell|_h=\|\ell\|_h^2$, and at a real place
$|\ell|_h=\|\ell\|_h$ on the complexification. With
$d^c=\tfrac{i}{4\pi}(\bar\partial-\partial)$, we have
\[
 c_1(L,h)=dd^c\phi_{\BC}=2dd^c\phi_{\BR}.
\]
Mixed Monge--Amp\`ere measures have their geometric intersection degrees as
masses. In particular $\Psi_{\BR}=\Psi_{\BC}/2$ for conjugation-invariant
metrics. In the analytic torsion estimates all Hermitian and $L^2$
norms have their usual meaning.

Over a complete nonarchimedean field we fix a nontrivial absolute
value. Analytic scalar extensions and completed residue fields carry
the extension of the absolute value on the ground field. For $x$ algebraic
over $F$ this is
$|x|_F=|\operatorname{Norm}_{F(x)/F}(x)|_F^{1/[F(x):F]}$.
In adelic base change, each finite extension $F'/F$ of local fields
has its own Haar normalization. Then $|a|_{F'}=|a|_F^{[F':F]}$ on $F$,
and local energies are multiplied by $[F':F]$.

Over a complete nonarchimedean field $F$, sections $s_0,\ldots,s_N$
generating $L^{\otimes m}$ define a rational Fubini--Study metric
on $L$ by the quotient of the maximum norm and passage to the $m$th
root. In a frame $\ell$, its weight is
$m^{-1}\log\max_i|f_i|_F$, where $s_i=f_i\ell^{\otimes m}$.
For discretely valued $F$, these are the metrics from rational model
bundles with a globally generated positive tensor power
\cite[Definition~5.2 and Theorem~5.14]{BE21}.

\section{Heights of fibrations}\label{sec:heights}

Over number fields, the proof of Theorem~\ref{thm:main-height}
combines a lower bound for analytic torsion, arithmetic
Riemann--Roch, and a slope estimate. Over function fields,
ordinary Riemann--Roch gives the model inequality. In both cases,
we pass to adelic bundles by approximation.

\subsection{One-sided analytic torsion estimate}\label{ss:absolute-torsion}

In the following, we prove an asymptotic torsion estimate.
Related torsion estimates for positive tensor powers, symmetric
powers, and smooth families were proved by Bismut--Vasserot
\cite{BV89,BV90} and Puchol \cite{P}.

On a smooth projective complex variety $X$ of dimension $d$,
let $L$ and $M$ be smooth Hermitian holomorphic line bundles.
We write $\nabla_L$ for the Chern connection and $\ell$ for a
local nonvanishing holomorphic section. We put
$h_L(\ell)=h_L(\ell,\ell)$ and normalize curvature by
\[
 c_1(L,h_L)=\frac{i}{2\pi}\nabla_L^2
   =\frac{1}{2\pi i}\partial\bar\partial\log h_L(\ell).
\]
Thus $h_L(\ell)=\|\ell\|_{h_L}^2$ is the squared Hermitian norm.
For a chosen K\"ahler form on $X$, let
\[
 \Delta_{n,q}=\bar\partial\bar\partial^*
                   +\bar\partial^*\bar\partial
\]
be the Dolbeault Laplacian on $(0,q)$-forms with values in
$L^{\otimes n}\otimes M$. Its adjoints use the K\"ahler metric and the
bundle metrics. We write
$\zeta_{n,q}(s)=\tr'(\Delta_{n,q}^{-s})$ for its spectral zeta
function, continued meromorphically from $\Re s>d$, where the
prime means that zero eigenvalues are omitted. All cohomology
metrics below are the ordinary Hodge $L^2$ metrics on the total space.

Let $(Y,g)$ be a closed Riemannian manifold of real dimension
$2d$, and let $V\to Y$ be a Hermitian bundle of rank $r$ with a
unitary connection $\nabla$ and a smooth self-adjoint endomorphism
$W$. For $\kappa>0$, we put
\[
 \mathcal P=\kappa\nabla^*\nabla+W,
\]
where $\nabla^*$ is the formal $L^2$-adjoint. Let $r_0>0$ be a
lower bound for the injectivity radius. Coefficient bounds
are taken uniformly in geodesic normal coordinates on balls of
radius $r_0/2$, with unitary frames given by radial parallel transport.

\begin{lemma}\label{lem:uniform-heat-remainder}
We have the following asymptotic formula.
\[
 \tr(e^{-t\mathcal P})=\sum_{j=0}^d A_jt^{j-d}
       +B(t),
 \qquad 0<t\le1,
\]
with $|A_j|\le c\operatorname{vol}(Y,g)$ and $|B(t)|\le ct\operatorname{vol}(Y,g)$, where 
 $c$ is a positive constant  depending only on $d$, $r$, $\kappa$,
$r_0$, the ellipticity bounds for $g$, and bounds for the
coefficients of $g$, $g^{-1}$, $\nabla$, $W$ and finitely many
of their derivatives.
\end{lemma}
\begin{proof}
We choose a cutoff supported inside the coordinate radius and
equal to one near the diagonal. The local parametrix is then
\[
 Q(t,x,y)=\chi(\ell)(4\pi\kappa t)^{-d}
  e^{-\ell^2/(4\kappa t)}\sum_{j=0}^d t^j u_j(x,y),
 \qquad \ell=\operatorname{dist}(x,y).
\]
We include the volume density in $u_0$. The transport equations
express the coefficients $u_j$ and $\mathcal P_xu_d$ in terms of the bounded
jets in the hypothesis; see \cite[Chapter 2]{BGV92}. It follows
that the defect $R=(\partial_t+\mathcal P_x)Q$ is uniformly bounded.
Indeed, the last transport term is bounded by
$ct^{d-d}e^{-\ell^2/(4\kappa t)}$, while the derivatives of the
cutoff are supported away from the diagonal.
Since $W\ge-c_0$, the maximum principle for vector bundles gives
$\|e^{-s\mathcal P}\|_{L^\infty\to L^\infty}\le e^{c_0s}$.
If $K$ is the heat kernel of $\mathcal P$, we apply Duhamel's formula
with $y$ fixed to obtain
\[
 K(t,\cdot,y)-Q(t,\cdot,y)
   =-\int_0^t e^{-(t-s)\mathcal P}R(s,\cdot,y)\,ds.
\]
Thus $\sup_{x,y}|K-Q|\le ct$. Taking the trace on the diagonal
and integrating over $Y$ proves the assertion.
\end{proof}

Let $f:X\to C$ be a surjective morphism from a smooth
projective complex variety of dimension $d\ge1$ to a smooth
projective complex curve. At a singular point $x$ of a fiber,
$\ker df_x=T_xX$. We fix a K\"ahler form $\omega_C$ on
$C$ and a smooth Hermitian holomorphic line bundle $M$ on $X$.
For a K\"ahler metric on $X$, we use the torsion convention
\[
 \tau_n=\sum_{q=0}^d(-1)^{q+1}q\,\zeta'_{n,q}(0).
\]
The K\"ahler metric and bundle metrics remain fixed as $n$ varies.

\begin{lemma}\label{lem:absolute-torsion}
Let $\bar L=(L,h_L)$ be a smooth Hermitian line bundle whose
curvature $\omega=c_1(\bar L)$ is strictly positive on $\ker df$
at every point, including the singular points of the fibers. There are constants
$a,c>0$, independent of $n$, such that
\[
 \Omega_a=\omega+a f^*\omega_C
\]
is K\"ahler and that the holomorphic analytic torsion of
$L^{\otimes n}\otimes M$, with the metric $\Omega_a$ on the total space, satisfies
\[
 \tau_n\ge-cn^d\log n\qquad(n\gg0).
\]
\end{lemma}
\begin{proof}
\smallskip\noindent\emph{Positivity and the spectral gap.}
The form $\beta=f^*\omega_C$ has rank at most one and kernel
$\ker df$. The form $\omega$ is positive near the zero locus of
$\beta$; on the remaining compact set, the Schur complement gives
$b>0$ with $\omega+b\beta>0$. We choose $a>3b$.
Relative to $\Omega_a=\omega+a\beta$, the
endomorphism defined by $\beta$ has rank at most one, with its
possible nonzero eigenvalue $u$ satisfying
\[
 0\le u<\frac1{a-b},
 \qquad
 \Omega_a-(a-b)\beta=\omega+b\beta>0.
\]
Consequently $\omega$ has $d-1$ eigenvalues equal to $1$ and a
remaining eigenvalue
\[
 \lambda=1-au> -\frac{b}{a-b}>-\frac12.
\]
At a point where $\beta=0$, all $d$ eigenvalues equal $1$.
Every sum of $q\ge2$ eigenvalues is therefore bounded below by a
fixed positive constant. We identify $(0,q)$-forms with values in
$L^{\otimes n}\otimes M$ with $(d,q)$-forms with values in
$L^{\otimes n}\otimes M\otimes K_X^{\otimes(-1)}$, using the induced metric on
$K_X^{\otimes(-1)}$. Then the Kodaira--Nakano identity gives, for some $c_0>0$,
\begin{equation}\label{eq:absolute-high-degree-gap}
 \Delta_{n,q}\ge c_0 n\qquad(q\ge2,\ n\gg0).
\end{equation}
Here, the curvatures of the fixed factors $M$ and $K_X^{\otimes(-1)}$ contribute
only a bounded zeroth-order term. When $d=1$, this assertion is
empty and the degree-zero argument below suffices.

\smallskip\noindent\emph{Uniform small-time estimates.}
For the K\"ahler metric $g$ of $\Omega_a$, put $g_n=ng$.
Its Dolbeault Laplacian is $\mathcal P_n=\Delta_{n,q}/n$.
The injectivity radius of $g_n$ is bounded below, and
$\operatorname{vol}(X,g_n)=n^d\operatorname{vol}(X,g)$.
We verify the required coefficient bounds in radial unitary frames
on balls of fixed $g_n$-radius. If $F_L=\nabla_L^2$ and
$F_M=\nabla_M^2$, then
\[
 \|\nabla^j(nF_L+F_M)\|_{g_n}=O(n^{-j/2}),
\]
so the connection coefficients have uniformly bounded jets of
every fixed order. The same holds for the coefficients on
differential forms and for the potential in the Weitzenb\"ock
formula. Since $\Lambda^{0,q}T^*X\otimes L^{\otimes n}\otimes M$
has fixed rank, Lemma~\ref{lem:uniform-heat-remainder} applies to
$\mathcal P_n$ with constants independent of $n$. After dividing its
heat-trace expansion by $n^d$, we obtain
\begin{equation}\label{eq:absolute-uniform-heat}
 H_{n,q}(t):=n^{-d}\tr(e^{-t\Delta_{n,q}/n})
    =\sum_{j=-d}^{0}a_{n,q,j}t^j+O(t),\qquad 0<t\le1,
\end{equation}
with uniformly bounded coefficients and remainder. In particular,
$H_{n,q}(1)\le c$. These estimates have the same rescaling as
\cite[proof of Proposition 2.4, equations (2.13)--(2.17)]{Berman11}.

\smallskip\noindent\emph{Zeta functions and torsion.}
We remove the kernel by putting $h_{n,q}=\dim\ker\Delta_{n,q}$
and $F_{n,q}(t)=H_{n,q}(t)-n^{-d}h_{n,q}$. Since
$H_{n,q}(1)\le c$, we have
$n^{-d}h_{n,q}\le c$. Subtracting the kernel from
\eqref{eq:absolute-uniform-heat} thus gives the same expansion
with bounded coefficients, say $b_{n,q,j}$, and $O(t)$ remainder.
For $\Re s>d$, we define the normalized spectral zeta function
\[
 \widetilde\zeta_{n,q}(s)
 =\frac1{\Gamma(s)}\int_0^\infty t^{s-1}F_{n,q}(t)\,dt
 =n^{s-d}\zeta_{n,q}(s).
\]
Its continuation near zero is obtained by replacing the integral
over $(0,1)$ by
\[
 \int_0^1 t^{s-1}
  \left(F_{n,q}(t)-\sum_{j=-d}^{0}b_{n,q,j}t^j\right)dt
       +\sum_{j=-d}^{0}\frac{b_{n,q,j}}{s+j}.
\]
The $O(t)$ remainder and bounded coefficients show that
$|\widetilde\zeta_{n,q}(0)|\le c$ and
\begin{equation}\label{eq:absolute-zeta-tail}
 \widetilde\zeta'_{n,q}(0)
   =R_{n,q}+\int_1^\infty F_{n,q}(t)\,\frac{dt}{t},
 \qquad |R_{n,q}|\le c.
\end{equation}
For fixed $n$, the positive spectral gap makes the integral
convergent.

For $q=0$, the function $F_{n,0}$ is nonnegative, so
\eqref{eq:absolute-zeta-tail} gives
$\widetilde\zeta'_{n,0}(0)\ge-c$. For $q\ge2$, the gap
\eqref{eq:absolute-high-degree-gap} makes $h_{n,q}=0$ and gives
\[
 0\le F_{n,q}(t)=H_{n,q}(t)
       \le c e^{-c_0(t-1)}\qquad(t\ge1).
\]
Hence $|\widetilde\zeta'_{n,q}(0)|\le c$ in those degrees.
Differentiating the scaling identity yields
\[
 \zeta'_{n,q}(0)=n^d\left(
   \widetilde\zeta'_{n,q}(0)
    -\log n\,\widetilde\zeta_{n,q}(0)\right).
\]
Thus, we have proved
\begin{equation}\label{eq:absolute-zeta-bounds}
 \zeta'_{n,0}(0)\ge-cn^d\log n,
 \qquad
 |\zeta'_{n,q}(0)|\le cn^d\log n\quad(q\ge2).
\end{equation}

Exactness of the Dolbeault complex on each positive eigenspace gives
\[
 \sum_{q=0}^d(-1)^q\zeta_{n,q}(s)=0.
\]
Eliminating degree one gives
\begin{equation}\label{eq:absolute-torsion-cancellation}
 \tau_n=\zeta'_{n,0}(0)
       +\sum_{q=2}^d(-1)^{q+1}(q-1)\zeta'_{n,q}(0).
\end{equation}
The bounds \eqref{eq:absolute-zeta-bounds} thus prove the lemma.
\end{proof}

For squared Hermitian norms on the determinant of cohomology,
our Quillen convention is
\begin{equation}\label{eq:quillen-torsion-convention}
 h_Q=e^{\tau_n}h_{L^2}.
\end{equation}
Thus, passing from Quillen degree to $L^2$ degree adds
$\tau_n/2$ at each archimedean embedding, before summing over
embeddings.

\subsection{Arithmetic Riemann--Roch inequalities}\label{ss:arithmetic-riemann}
We prove a lower bound for the arithmetic degree of the determinant
of cohomology, followed by an upper bound for the maximal slope of
spaces of sections, both with their $L^2$ metrics.

Let $K$ be a number field with ring of integers $\mathcal O_K$.
Let $\mathcal X$ be a regular integral arithmetic scheme,
projective and flat over $\mathcal O_K$, with structure morphism
$\mathcal g:\mathcal X\to\Spec\mathcal O_K$ and smooth generic
fiber $X$ of dimension $d$. Let $f:X\to C$ be a surjective
morphism to a smooth projective curve.
For each complex emebedding $\sigma: K\to \BC$, we fix conjugation-compatible K\"ahler forms $\omega_{C,\sigma}$
on $C_\sigma$ and a smooth Hermitian line bundle
$\bar{\mathcal M}$ on $\mathcal X$.
The subscript $L^2,\Omega$ below denotes the ordinary Hodge
metric on the determinant of cohomology, computed on the total
space with the K\"ahler metrics $\Omega_\sigma$. Constants may
depend on all fixed data, including the chosen K\"ahler metrics,
but not on $n$.

\begin{thm}\label{thm:ari}
Let $\bar{\mathcal L}$ be a smooth Hermitian line bundle whose
curvature $\omega_\sigma=c_1(L_\sigma,h_\sigma)$ at each
complex embedding $\sigma$ is strictly positive on $\ker df_\sigma$
at every point, including the singular points of the fibers. For sufficiently large
fixed $a_\sigma$, chosen with $a_{\bar\sigma}=a_\sigma$, the forms
\[
 \Omega_\sigma=\omega_\sigma+a_\sigma f_\sigma^*\omega_{C,\sigma}
\]
are K\"ahler, and there is a constant $c>0$ such that
\begin{equation}\label{eq:ARR}
 \widehat\deg\det R\mathcal g_*
       (\mathcal L^{\otimes n}\otimes\mathcal M)_{L^2,\Omega}
 \ge \frac{\bar{\mathcal L}^{d+1}}{(d+1)!}n^{d+1}
                           -cn^d\log n\qquad(n\gg0).
\end{equation}
\end{thm}
\begin{proof}
We choose $a_\sigma$ as in Lemma~\ref{lem:absolute-torsion}.
The absolute arithmetic Riemann--Roch theorem \cite{GS92,Fa92},
applied to $\mathcal g$ with the metrics $\Omega_\sigma$, gives
\[
 \widehat\deg_Q\det R\mathcal g_*
       (\mathcal L^{\otimes n}\otimes\mathcal M)
   =\frac{\bar{\mathcal L}^{d+1}}{(d+1)!}n^{d+1}+O(n^d).
\]
The Todd class, fixed twist, and $R$-genus contribute $O(n^d)$.
By \eqref{eq:quillen-torsion-convention}, passing to the Hodge
$L^2$ metric adds half the sum of the analytic torsions, and
Lemma~\ref{lem:absolute-torsion} gives
\[
 \tau(L_\sigma^{\otimes n}\otimes M_\sigma)\ge-c_\sigma n^d\log n.
\]
Summing over embeddings proves \eqref{eq:ARR}.
\end{proof}

We fix smooth Hermitian line bundles
$\overline{\mathcal P}$ and $\overline{\mathcal Q}$ on
$\mathcal X$ and smooth positive volume forms on the complex
fibers, which define the $L^2$ norms. For a nef adelic
$\mathbb Q$-line bundle $\bar A_0$ on $X$ with $A_0^d>0$, we put
\[
 a_{A_0}(\bar N)=\frac{1}{A_0^d}(\bar N\cdot\bar A_0^d).
\]
The ratio $a_{A_0}(\bar N)$ is unchanged when $\bar A_0$ is multiplied
by a positive real scalar. For model bundles, we identify mixed
arithmetic intersections with the corresponding adelic intersections.
With $\bar A=(A_0^d)^{-1/d}\bar A_0$ we have
\[
 a_{A_0}(\bar N)=a_A(\bar N)=(\bar N\cdot\bar A^d),\qquad A^d=1.
\]

For a nonzero Hermitian vector bundle $\bar E$ on
$\Spec\mathcal O_K$, we use unnormalized arithmetic slopes:
\[
 \widehat\mu_{\max}(\bar E)
   =\sup_{0\ne F\subset E}
       \frac{\widehat{\deg}\bar F}{\operatorname{rk}F},
\]
where $F$ runs over saturated submodules with the induced metrics.

\begin{lem}\label{lem:linear-slope}
Assume that $A_0$ is ample. For integers $n,m\ge0$ such that
$H^0(\mathcal X,\mathcal P^{\otimes n}\otimes\mathcal Q^{\otimes m})$
is nonzero, we have
\begin{equation}\label{eq:mixed-slope-bound}
 \widehat\mu_{\max}
 H^0(\mathcal X,\mathcal P^{\otimes n}\otimes\mathcal Q^{\otimes m})_{L^2}
 \le n a_{A_0}(\overline{\mathcal P})
       +m a_{A_0}(\overline{\mathcal Q})
       +O_{\mathcal P,\mathcal Q,\bar A_0}(\log(n+m+2)).
\end{equation}
\end{lem}
\begin{proof}
The two-parameter Gromov inequality
\cite[Proposition 2.13]{Yuan08} gives
\[
 \|s\|_{\sup}\le c_0(n+m+1)^{c_1}\|s\|_{L^2}
\]
at every complex embedding. Let $\bar Q$ be a Hermitian line
subbundle of the section module with the induced $L^2$ metric.
Evaluation gives a nonzero section $\ell$ of
$\mathcal P^{\otimes n}\otimes\mathcal Q^{\otimes m}\otimes\mathcal g^*Q^{\otimes(-1)}$.
We enlarge the norms of $Q$ by $c_0(n+m+1)^{c_1}$ and denote
the resulting Hermitian line by $\bar Q'$. The evaluation section $\ell$
now has sup norm at most one. Its arithmetic divisor is therefore
effective, so its intersection with $\bar A_0^d$ is nonnegative
by \cite[Lemma 5.2.13(b)]{YZ26}. We obtain
\[
 A_0^d\,\adeg\bar Q'
 \le n(\overline{\mathcal P}\cdot\bar A_0^d)
    +m(\overline{\mathcal Q}\cdot\bar A_0^d).
\]
Changing the metric back to $\bar Q$ adds
$O(\log(n+m+2))$ to its degree. We then apply Minkowski's theorem,
which bounds the maximal slope by the largest degree of a line
subbundle plus $O(\log r)$, where $r$ is the rank
\cite[Proposition 2.7]{Chen09}. Since the rank grows at most
polynomially in $n+m$, we obtain \eqref{eq:mixed-slope-bound}.
Chen's normalized degrees are multiplied by the fixed factor
$[K:\BQ]$. Compatibility of intersections allows the fixed adelic
$\bar A_0$ on every model \cite[Propositions 4.1.1--4.1.2]{YZ26}.
\end{proof}

\Needspace{17\baselineskip}

\subsection{Height inequalities}\label{ss:height-model}

We prove the height inequality for models
over number fields and function fields, then pass to adelic bundles by
approximation in Theorem~\ref{thm:ahi}.

\subsubsection*{Model inequality}

On a singular complex variety, smooth metrics are understood
in the ambient sense. Semipositivity on a singular fiber means
semipositivity after resolution of each reduced irreducible
component. Arithmetic intersections are defined through the
associated integrable adelic bundles on the projective generic fiber
\cite[Propositions 4.1.1--4.1.2]{YZ26}.

For the next theorem, let
\[
 \mathcal X\xrightarrow{\mathcal f}\mathcal C
       \xrightarrow{\varpi}\Spec\CO_K
\]
be proper flat morphisms of integral arithmetic schemes, with
projective generic fibers $X$ and $C$ of dimensions $d$ and $1$.
We fix a nef adelic polarization $\bar A$ on $X$, with $A$
ample and $A^d=1$. 
\begin{theorem}\label{thm:hi}
Let $\overline{\mathcal L}$ be a smooth Hermitian line bundle whose
underlying bundle is nef on every fiber of $\mathcal f$
and whose metric is semipositive on every complex fiber. Suppose
that $L=\mathcal L_K$ is nef on $X$, and write $\bar L$ for
the induced adelic line bundle. Then
\begin{equation}\label{eq:quantitative-model-height}
 \bar L^{d+1}\le(d+1)L^d(\bar L\cdot\bar A^d).
\end{equation}
In particular, if $L^d=0$, then $\bar L^{d+1}\le0$.
\end{theorem}
\begin{proof}
\smallskip\noindent\emph{The regular case.}
We first prove the inequality when $\mathcal X$ is regular and
projective over $\CO_K$, $C$ is smooth, and the curvature of
$\overline{\mathcal L}$ is semipositive on $\ker df$ at every
complex point, including the singular points of the fibers.
Here $\mathcal f$ need only be proper, with $\mathcal C$ integral
and proper over $\CO_K$; the fibers of $\varpi$ have dimension
at most one. This permits the alterations below.

To apply Theorem~\ref{thm:ari}, we perturb
$\overline{\mathcal L}$ by a globally ample bundle $\mathcal A'$
with generic fiber $A'=\mathcal A'_K$ and a smooth metric of strictly positive curvature. For rational
$\epsilon=p/q\in(0,1]$, with $p,q>0$, we put
\[
 \overline{\mathcal M}=q\overline{\mathcal L}
                            +p\overline{\mathcal A'},
 \qquad
 \overline{\mathcal L}_\epsilon
       =\overline{\mathcal L}+\epsilon\overline{\mathcal A'}.
\]
Then $\mathcal M$ is $\mathcal f$-ample, $M=\mathcal M_K$ is ample
on $X$, and its curvature is strictly positive on $\ker df$
everywhere.

For $k\gg0$, relative Serre vanishing gives
$R^j\mathcal f_*\mathcal M^{\otimes k}=0$ for $j>0$.
Since $\varpi$ has fibers of dimension at most one, Leray gives
\[
 H^j(\mathcal X,\mathcal M^{\otimes k})=0\qquad(j\ge2).
\]
Since $M$ is ample on $X$, generic Serre vanishing also shows
that $H^1(\mathcal X,\mathcal M^{\otimes k})$ is finite torsion.
We put $F_k=H^0(\mathcal X,\mathcal M^{\otimes k})$ and write
\[
 r_k=\rk F_k=\frac{M^d}{d!}k^d+O_{p,q}(k^{d-1}).
\]
We equip $F_k$ with the $L^2$ metric on the total space used in
Theorem~\ref{thm:ari}. The preceding cohomology calculation gives
\[
 \adeg\det Rg_*\mathcal M^{\otimes k}{}_{L^2}
 =\adeg(F_k{}_{L^2})
              -\log\#H^1(\mathcal X,\mathcal M^{\otimes k})
 \le\adeg(F_k{}_{L^2}).
\]
The K\"ahler metric on the total space is fixed for each $p,q$; it is
independent of $k$. Thus \eqref{eq:mixed-slope-bound}, with exponents
$qk$ and $pk$, gives
\[
 \widehat\mu_{\max}(F_k{}_{L^2})
     \le k\,a_A(\overline{\mathcal M})
              +O_{p,q}(\log(k+2)),
 \qquad
 a_A(\overline{\mathcal M})=q a_A(\overline{\mathcal L})
                              +p a_A(\overline{\mathcal A'}).
\]
The absolute arithmetic Riemann inequality now yields
\begin{align*}
 \frac{\overline{\mathcal M}^{d+1}}{(d+1)!}k^{d+1}
       -O_{p,q}(k^d\log k)
 &\le\adeg\det Rg_*\mathcal M^{\otimes k}{}_{L^2}\\
 &\le r_k\widehat\mu_{\max}(F_k{}_{L^2}).
\end{align*}
For fixed $p,q$, we divide by $k^{d+1}$ and let $k\to\infty$.
We obtain
\[
 \overline{\mathcal M}^{d+1}
       \le(d+1)a_A(\overline{\mathcal M})M^d.
\]
After division by $q^{d+1}$,
\begin{equation}\label{eq:perturb-height-bound}
 \begin{aligned}
 \overline{\mathcal L}_\epsilon^{d+1}
 &\le(d+1)a_A(\overline{\mathcal L}_\epsilon)(L+\epsilon A')^d\\
 &\le(d+1)c(1+\epsilon)(L+\epsilon A')^d,
 \end{aligned}
\end{equation}
where $c$ is independent of $0<\epsilon\le1$.
The analytic estimates were used with $\epsilon$ fixed.
Letting $\epsilon\to0$ in \eqref{eq:perturb-height-bound} proves
the regular case by continuity of intersections.

\smallskip\noindent\emph{Reduction by alterations.}
After normalizing the generic base curve and resolving as
necessary, semistable alteration over a curve
\cite[Corollary 5.1]{dJ97} gives a commutative diagram
\[
 \begin{tikzcd}
 X_1 \arrow[r,"p_1"] \arrow[d,"f_1"'] & X \arrow[d,"f"]\\
 C_1 \arrow[r] & C
 \end{tikzcd}
\]
where $C_1\to C$ is finite, $p_1$ is generically finite, $X_1$ and
$C_1$ are smooth projective, and $f_1$ is strictly semistable.
The curvature of $p_1^*\overline{\mathcal L}$ is semipositive on
$\ker df_1$ over the dense open subset mapping to the smooth
locus of $f$. Since $\ker df_1$ is a vector bundle on
the smooth locus of $f_1$, semipositivity extends to this locus
by continuity.

Let $x$ be a singular point of a fiber of $f_1$. In local coordinates
centered at $x$, the map is $t=z_1\cdots z_r$ with $r\ge2$.
Fix $a_1,\ldots,a_r\in\mathbb C^\times$ and approach $x$ along
$z_i=s/a_i$ for $i\le r$, with the remaining coordinates zero.
For $s\ne0$, the tangent hyperplane to the fiber is
$\ker(\sum_{i=1}^r a_i\,dz_i)$, since $dt=t\sum_{i=1}^r dz_i/z_i$.
By continuity, the curvature at $x$ is semipositive on this
hyperplane. By density the same holds for every nonzero
$(a_1,\ldots,a_r)\in\mathbb C^r$. Since $r\ge2$, every tangent
vector $v$ at $x$ satisfies $\sum_i a_i v_i=0$ for some such
$(a_1,\ldots,a_r)$. The curvature is therefore semipositive on
$T_xX_1=\ker df_{1,x}$.

To extend the diagram to arithmetic models, we let
$\mathcal C_1$ be the normalization of $\mathcal C$ in $K(C_1)$.
By excellence, $\mathcal C_1\to\mathcal C$ is finite. Thus
$\mathcal C_1$ is an integral proper arithmetic surface with
one-dimensional fibers over $\Spec\CO_K$. We choose a
projective embedding of $X_1$ and take the closure of its graph
in $\mathcal X\times\mathcal C_1\times\BP^N_{\CO_K}$,
where the products are over $\CO_K$. We denote the resulting
proper model by $\mathcal Z_1$.
The regular alteration theorem \cite[Theorem 8.2]{dJ96} gives a
regular arithmetic alteration $\mathcal Y\to\mathcal Z_1$
after a finite extension of constants. In that theorem,
$\mathcal Y$ is open in a projective model over the extended
arithmetic base. Since $\mathcal Z_1$ is proper over $\CO_K$, so is its
alteration; the separatedness of the extended arithmetic base
then makes $\mathcal Y$ proper over that base as well. This
open immersion is consequently proper and hence closed. Its
image is dense, so $\mathcal Y$ equals the projective model. In particular, its map to the separated
base $\mathcal C_1$ is projective.
We write $p:\mathcal Y\to\mathcal X$ for the composite and
$\mathcal f_Y:\mathcal Y\to\mathcal C_1$ for the induced map.
These maps fit into the commutative diagram
\[
 \begin{tikzcd}[column sep=large, row sep=large]
 \mathcal Y \arrow[r] \arrow[rr, bend left=20, "p"]
            \arrow[dr, "\mathcal f_Y"']
   & \mathcal Z_1 \arrow[r] \arrow[d]
   & \mathcal X \arrow[d, "\mathcal f"]\\
   & \mathcal C_1 \arrow[r] & \mathcal C.
 \end{tikzcd}
\]
The generic degree $D=\deg p$ is taken over the original $K$,
including constant extensions. The analytic estimates apply
componentwise after complex base change.

For the induced morphism $r:Y\to X_1$, with $f_Y=f_1r$,
\[
 v\in\ker d(f_1r)\quad\Longrightarrow\quad
 dr(v)\in\ker df_1.
\]
Thus $p^*\overline{\mathcal L}$ has semipositive curvature on
$\ker df_Y$. Its generic bundle is nef, and it is nef on every
fiber of $\mathcal f_Y$, since vertical curves map to vertical
curves or points. The regular case therefore applies.

To descend the inequality, we also approximate the pullback of
the polarization by ample bundles. Let $\bar A_0$ be the rational
polarization whose normalization is $\bar A$. Since $p^*A_0$ is
nef and big, we choose a nef adelic $\mathbb Q$-line bundle
$\bar G$ with ample underlying bundle on $Y$ and put
$\bar A_\delta=p^*\bar A_0+\delta\bar G$ for rational $\delta>0$.
We apply the regular case with the normalized polarization
$\bar A_\delta/(A_\delta^d)^{1/d}$, whose slope coefficient is
$a_{A_\delta}$ as defined before Lemma~\ref{lem:linear-slope}.
Since $(p^*A_0)^d=D A_0^d>0$, we can let $\delta\to0$.
The projection formula for adelic intersections gives
\[
 (p^*\overline{\mathcal L})^{d+1}=D \overline{\mathcal L}^{d+1},
 \qquad (p^*L)^d=D L^d,
 \qquad
 a_{p^*A_0}(p^*\overline{\mathcal L})
       =a_{A_0}(\overline{\mathcal L})
       =a_A(\overline{\mathcal L}).
\]
Thus the resulting inequality on $\mathcal Y$ is
\[
 D \overline{\mathcal L}^{d+1}
 \le D(d+1)a_A(\overline{\mathcal L})L^d.
\]
Dividing by $D>0$ proves \eqref{eq:quantitative-model-height}.
\end{proof}

\subsubsection*{Function-field inequality}

For the next theorem, let $T$ be a regular integral projective
curve over an arbitrary field $k$, put $K=k(T)$, and consider
proper flat morphisms of integral varieties
\[
 \mathcal X\xrightarrow{\mathcal f}\mathcal C
       \xrightarrow{\pi}T,
 \qquad \dim\mathcal X=d+1,\quad\dim\mathcal C=2.
\]
We fix a nef $\mathbb Q$-line bundle $\mathcal A_0$ on
$\mathcal X$ whose generic restriction $A_0$ is ample, and put
\[
 \mathcal A=\frac{\mathcal A_0}{(A_0^d)^{1/d}},\qquad
 A=\mathcal A_K,\qquad A^d=1.
\]
We write $\bar A_0$ and $\bar A$ for the induced adelic bundles.
We use geometric adelic line bundles over $K/k$, in the sense of
\cite[\S2.5.3]{YZ26}, and put
\[
 c_k=\begin{cases}
       1,&k\text{ infinite},\\
       \log\#k,&k\text{ finite}.
     \end{cases}
\]
At a closed point $v$ of $T$, our normalization is
$\log|a|_v=-c_k[\kappa(v):k]\operatorname{ord}_v(a)$.
Thus an adelic intersection of $d+1$ model bundles is $c_k$
times their geometric intersection on $\mathcal X$ over $k$;
intersections of $d$ underlying bundles on $X$ are taken over $K$.
In particular, this convention agrees with the Haar-normalized
absolute values when $k$ is finite. In the proof below,
intersections written with calligraphic model bundles remain
geometric intersection numbers over $k$.

\begin{thm}\label{prop:function-field-height}
Let $\mathcal L$ be nef on every fiber of $\mathcal f$, with nef
generic restriction $L$ on $X=\mathcal X_K$, and let $\bar L$
be its induced adelic line bundle. Then
\begin{equation}\label{eq:function-field-model-height}
 \bar L^{d+1}\le(d+1)L^d(\bar L\cdot\bar A^d).
\end{equation}
If $L^d=0$, then $\bar L^{d+1}\le0$.
\end{thm}
\begin{proof}
By Chow's lemma and the projection formula, we may replace
$\mathcal X$ by a projective birational modification. Pullback
preserves the nefness hypotheses and $A^d=1$. The morphism to
$\mathcal C$ is projective, so relative Serre vanishing applies
without flatness. We put
$a_A(\mathcal U)=(\mathcal U\cdot\mathcal A^d)$.

We choose an ample line bundle $\mathcal M$ on $\mathcal X$ to
make $\mathcal L$ relatively ample by perturbation. For rational
$\epsilon=p/q\in(0,1]$, with $p,q>0$, we put
\[
 \mathcal N=q\mathcal L+p\mathcal M,
 \qquad \mathcal L_\epsilon=\mathcal L+\epsilon\mathcal M.
\]
Then $\mathcal N$ is $\mathcal f$-ample and $N=\mathcal N_K$
is ample. We write $g=\pi\mathcal f$ and put
$F_n=g_*\mathcal N^{\otimes n}$ for $n\gg0$.
Since $F_n$ is torsion free on the regular curve $T$, it is locally
free. By generic Serre vanishing, its rank is
\[
 r_n=\rk F_n=\frac{N^d}{d!}n^d+O_{p,q}(n^{d-1}).
\]
Relative Serre vanishing for $\mathcal f$ and the fact that $\pi$
has one-dimensional fibers give $R^ig_*\mathcal N^{\otimes n}=0$ for $i\ge2$.
Generic Serre vanishing also makes $\mathcal T_n=R^1g_*\mathcal N^{\otimes n}$ torsion.
Thus
\[
 \chi(\mathcal X,\mathcal N^{\otimes n})
   =\chi(T,F_n)-\dim_k H^0(T,\mathcal T_n)\le\chi(T,F_n).
\]

To bound $\chi(T,F_n)$, let $D_T$ be an ample line bundle on $T$
and put $b=\deg D_T>0$. A nonzero section of
$F_n\otimes D_T^{\otimes(-m)}$ defines an effective Cartier
divisor on $\mathcal X$. Intersecting this divisor with the nef
bundle $\mathcal A$ gives
\[
 \bigl(nc_1(\mathcal N)-mc_1(g^*D_T)\bigr)
            \cdot\mathcal A^d\ge0.
\]
Since $g^*D_T\cdot\mathcal A^d=bA^d=b$, this forces
$mb\le n a_A(\mathcal N)$. We therefore choose
\[
 m_n=\left\lfloor\frac{n a_A(\mathcal N)}{b}\right\rfloor+1.
\]
Then $H^0(T,F_n\otimes D_T^{\otimes(-m_n)})=0$, and curve
Riemann--Roch gives
\[
 \chi(T,F_n)\le r_n m_n b
       \le r_n\bigl(n a_A(\mathcal N)+b\bigr).
\]
This holds for either sign of $m_n$. Comparing leading
coefficients in Riemann--Roch on $\mathcal X$ now gives
\[
 \mathcal N^{d+1}\le(d+1)a_A(\mathcal N)N^d.
\]
We divide by $q^{d+1}$ to obtain
\[
 \mathcal L_\epsilon^{d+1}
 \le(d+1)a_A(\mathcal L_\epsilon)(L+\epsilon M)^d.
\]
Letting $\epsilon\to0$ gives the geometric model inequality.
Multiplying both sides by $c_k$ proves
\eqref{eq:function-field-model-height}.
\end{proof}

\subsubsection*{Adelic inequality}

Let $K$ be a number field or the function field $k(T)$ of a
regular integral projective curve over an arbitrary field $k$.
In the latter case, we use the geometric adelic theory with the
factor $c_k$ fixed above. Thus nefness means global adelic nefness;
in the geometric case, the nef model bundles used for approximation
are nef on the whole projective model.

For the definitions, we allow any projective morphism
$f:X\to Y$ of projective $K$-varieties. On reducible varieties,
we impose the conditions on each reduced irreducible component.

\begin{definition}
\label{def:relative-nef-adelic}
An adelic line bundle $\bar L$ on $X$ is
\emph{strongly relatively nef for $f$} if there are nef adelic line
bundles $\bar M$ on $X$ and $\bar N$ on $Y$, with $N$ ample, such that
\begin{equation}\label{eq:strong-relative-nef-adelic}
 \bar L=\bar M-f^*\bar N.
\end{equation}

An integrable adelic line bundle $\bar L$ is
\emph{relatively nef for $f$} if there is a fixed nef
adelic $\mathbb Q$-line bundle $\bar A$
on $X$ such that
$\bar L+\epsilon\bar A$ is strongly relatively nef for every rational
$\epsilon>0$.
\end{definition}

A strongly relatively nef bundle is integrable by its defining
decomposition. In Theorem~\ref{thm:main-height} we also assume that $L$ is nef on $X$.

\paragraph*{Base change and fibers.}\label{par:relative-nef-base-change}
Both conditions are preserved by base change. Indeed, let
$u:Y'\to Y$ be a morphism of projective varieties, and write
$X'=X\times_Y Y'$, with projections $p:X'\to X$ and $f':X'\to Y'$.
We choose a nef adelic bundle $\bar Q$ on $Y'$ with $Q$ ample.
Then
\[
 p^*\bar L=(p^*\bar M+f'^*\bar Q)
                 -f'^*(u^*\bar N+\bar Q).
\]
Both factors are nef and $u^*N+Q$ is ample. Applying the same
construction to $\bar L+\epsilon\bar A$ proves relative nefness,
with fixed perturbation $p^*\bar A$. The argument also applies
locally, with semipositivity in place of nefness, and under
extension of the ground field.

We call $\bar L$ \emph{fiberwise nef} if $L$ is nef on every
geometric fiber and its metrics restrict to locally nef metrics
on every fiber at every place. At a nonarchimedean place, local
nefness here is understood in the perturbative sense of
\cite[\S3.6.6]{YZ26}; the point case is explained in
Section~\ref{ss:local-metrics}. Base change to a fiber reduces
relative nefness to this local point case. Hence
\[
 \text{strongly relatively nef}\ \Longrightarrow\
 \text{relatively nef}\ \Longrightarrow\ \text{fiberwise nef}.
\]
The first implication uses $\bar A=0$.

Relative nefness does not imply global adelic nefness: for
$Y=\Spec K$ and $\widehat\deg\bar N>0$, the strongly relatively
nef bundle $-f^*\bar N$ has negative point heights.

Strong relative nefness is preserved by every integrable adelic twist
from $Y$. Indeed, if $\bar U=\bar U_1-\bar U_2$ with $\bar U_1,\bar U_2$
nef, then
\[
 \bar L+f^*\bar U
   =(\bar M+f^*\bar U_1)-f^*(\bar N+\bar U_2),
\]
and $N+U_2$ is ample.

For the next theorem, let $f:X\to C$ be a proper flat
morphism of integral projective varieties over a number field
$K$ or a function field $K=k(T)$ as above, with $\dim X=d$
and $\dim C=1$. We fix a nef adelic polarization $\bar A$,
with $A$ ample and $A^d=1$, and write
$a_A(\bar T)=(\bar T\cdot\bar A^d)$ for adelic
intersections, with the factor $c_k$ included in the function-field
case.

\begin{thm}\label{thm:ahi}
Let $\bar L$ be an integrable relatively nef adelic line bundle
with nef underlying bundle $L$. Then
\[
 \bar L^{d+1}\le(d+1)L^d
                          (\bar L\cdot\bar A^d).
\]
In particular, $L^d=0$ implies $\bar L^{d+1}\le0$.
\end{thm}
\begin{proof}
First suppose that $\bar L$ is strongly relatively nef, and write
$\bar L=\bar M-f^*\bar N$ as in
\eqref{eq:strong-relative-nef-adelic}.

\smallskip\noindent\emph{Number fields.}
Suppose first that $K$ is a number field. We choose a generically
finite alteration $p:Y\to X$ over a finite cover $q:C'\to C$
so that $Y,C'$ are smooth projective and $g:Y\to C'$ is flat.
We also choose a nef adelic bundle $\bar G$ on $Y$, induced by
an ample Hermitian model, with $G$ ample. For a fixed rational
$t>0$, the bundles $p^*M+tG$ and $q^*N$ are ample.

We first approximate $p^*\bar M+t\bar G$ and $q^*\bar N$ by
nef model bundles with the same underlying generic bundles. At
the complex places, ampleness allows us to regularize the
continuous semipositive metrics uniformly by smooth metrics with
an arbitrarily small negative curvature error. We choose a fixed
smooth positive metric on each ample bundle and take a convex
combination of the regularized weight with this positive weight.
The coefficient of the positive weight can tend to zero while
absorbing the curvature error. We thus obtain smooth semipositive
metrics with the same uniform limit.

For either approximating bundle, we denote the original nef model
norm by $\|\cdot\|_{j,\mathrm{mod}}$ and the smooth semipositive
norm just constructed by $\|\cdot\|_j$. We choose
$\epsilon_j\to0$ so that
$|\log(\|\cdot\|_j/\|\cdot\|_{j,\mathrm{mod}})|\le\epsilon_j$
for both bundles at every archimedean embedding. We replace
$\|\cdot\|_j$ by $e^{-\epsilon_j}\|\cdot\|_j$.
The new norm is at most $\|\cdot\|_{j,\mathrm{mod}}$ at every
point, so the arithmetic degree on each horizontal curve can
only increase. Constant rescaling preserves the curvature, and
the model bundles at finite places are unchanged. We therefore
obtain smooth nef Hermitian model bundles
$\bar{\mathcal M}_{t,j}$ and $\bar{\mathcal N}_j$ converging
to $p^*\bar M+t\bar G$ and $q^*\bar N$, respectively
\cite[\S2.1.1, \S3.5.2 and \S A.5.2]{YZ26}.

We place the approximating bundles on common projective model
towers and flatten the morphisms after modifying the base model.
Pullback preserves nefness, so we can put
\[
 \bar{\mathcal D}_{t,j}
   =\bar{\mathcal M}_{t,j}-g_j^*\bar{\mathcal N}_j.
\]
Since $g_j^*\bar{\mathcal N}_j$ restricts to a constant metric
on each complex fiber, $\bar{\mathcal D}_{t,j}$ has semipositive
vertical curvature. The bundle $\mathcal D_{t,j}$ is nef on every
fiber of $g_j$, and its generic restriction is the ample bundle
$p^*L+tG$.
Let $\bar A_0$ be the rational representative of the fixed
normalized polarization $\bar A$, and set
$\bar A_\delta=p^*\bar A_0+\delta\bar G$ for rational $\delta>0$.
Theorem~\ref{thm:hi}, applied with the normalization of
$\bar A_\delta$, gives
\[
 \bar{\mathcal D}_{t,j}^{d+1}
 \le(d+1)a_{A_\delta}(\bar{\mathcal D}_{t,j})(p^*L+tG)^d.
\]
For fixed $t,\delta$, the continuity of intersections of nef
adelic bundles allows us to let $j\to\infty$. We obtain the
inequality for $p^*\bar L+t\bar G$. We then let $t\to0$ and
$\delta\to0$ in this order. The polarization denominator has
positive limit, since $(p^*A_0)^d=(\deg p)A_0^d>0$.
By the projection formula,
\[
 \begin{gathered}
 (p^*\bar L)^{d+1}=(\deg p)\bar L^{d+1},\quad
 (p^*L)^d=(\deg p)L^d,\\
 a_{p^*A_0}(p^*\bar L)=a_{A_0}(\bar L)=a_A(\bar L).
 \end{gathered}
\]
Dividing by $\deg p$ proves the strongly relatively nef case over a number
field. Degrees are taken over the original field $K$, including
any extension of constants in the alteration.

\smallskip\noindent\emph{Function fields.}
Let $K=k(T)$ and write $\bar A=(A_0^d)^{-1/d}\bar A_0$
with $\bar A_0$ rational. By the definition of nefness
\cite[Definition 2.5.2 and \S2.5.3]{YZ26}, there are strongly
nef adelic bundles $\bar U$ on $X$, $\bar Q$ on $C$, and
$\bar R$ on $X$ such that
\[
 \bar M+\epsilon\bar U,\qquad
 \bar N+\epsilon\bar Q,\qquad
 \bar A_0+\epsilon\bar R
\]
are strongly nef for every rational $\epsilon>0$. To keep the
difference nef on the generic fiber, we add the same base term
to both factors and put
\[
 \begin{aligned}
 \bar M_\epsilon&=\bar M+\epsilon\bar U
                              +\epsilon f^*\bar Q,\\
 \bar N_\epsilon&=\bar N+\epsilon\bar Q,\qquad
 \bar A_{0,\epsilon}=\bar A_0+\epsilon\bar R.
 \end{aligned}
\]
The three bundles $\bar M_\epsilon$, $\bar N_\epsilon$, and
$\bar A_{0,\epsilon}$ are strongly nef. Their difference
$\bar M_\epsilon-f^*\bar N_\epsilon=\bar L+\epsilon\bar U$
has nef underlying bundle, while $A_{0,\epsilon}$ is ample.
We normalize the polarization by
$\bar A_\epsilon=(A_{0,\epsilon}^d)^{-1/d}\bar A_{0,\epsilon}$.
For fixed $\epsilon$, we choose nef model approximations to
$\bar M_\epsilon$, $\bar N_\epsilon$, and $\bar A_{0,\epsilon}$
on common projective model towers over $T$. After modifying the
base model and flattening, we pull back the approximating bundles;
they remain nef. Thus the differences
\[
 \mathcal D_{\epsilon,j}
   =\mathcal M_{\epsilon,j}
                    -\mathcal f_j^*\mathcal N_{\epsilon,j}
\]
are nef on every fiber of $\mathcal f_j$, with fixed nef generic
bundle $L+\epsilon U$. We apply
Theorem~\ref{prop:function-field-height} with the normalized
polarization induced by $\mathcal A_{0,\epsilon,j}$ and let
$j\to\infty$. By continuity of adelic intersections, we obtain
\[
 (\bar L+\epsilon\bar U)^{d+1}
 \le(d+1)(L+\epsilon U)^d
       \bigl((\bar L+\epsilon\bar U)\cdot\bar A_\epsilon^d\bigr).
\]
The factor $c_k$ occurs on both sides of the model inequality.
Letting $\epsilon\to0$ proves the strongly relatively nef case.

To finish the proof over either field, let $\bar L$ be relatively
nef. We choose a fixed
nef adelic $\bar A'$ as in Definition~\ref{def:relative-nef-adelic}.
For every rational $\epsilon>0$, the preceding case applies to
$\bar L+\epsilon\bar A'$, whose underlying bundle is nef. Thus
\[
 (\bar L+\epsilon\bar A')^{d+1}
 \le(d+1)\bigl(a_A(\bar L)+\epsilon a_A(\bar A')\bigr)
                  (L+\epsilon A')^d.
\]
Both sides are polynomials in $\epsilon$. Letting
$\epsilon\to0$ proves the theorem.
\end{proof}

\paragraph*{Model bundles and relative nefness.}
Let $\bar{\mathcal L}$ satisfy the hypotheses of
Theorem~\ref{thm:hi}. The construction can be made on
$\mathcal X\to\mathcal C$ itself when both models are projective,
$X$ and $C$ are smooth,
and the curvature is semipositive on $\ker df$ at every complex
point. Let $\bar{\mathcal A}$
be arithmetically ample on
$\mathcal X$. For each rational $\epsilon>0$,
$\bar{\mathcal L}+\epsilon\bar{\mathcal A}$ is relatively
ample with strictly positive vertical curvature. A sufficiently
positive ample Hermitian base twist makes
\[
 \bar{\mathcal M}_\epsilon
   =\bar{\mathcal L}+\epsilon\bar{\mathcal A}
                           +\mathcal f^*\bar{\mathcal N}_\epsilon
\]
globally ample \cite[Tag 0892]{Stacks} with positive curvature.
Shrinking both norms by the same constant
makes $\bar{\mathcal M}_\epsilon$ and
$\bar{\mathcal N}_\epsilon$ nef without changing their difference.
Thus $\bar L$ is relatively nef on the original fibration.

In the generality of Theorem~\ref{thm:hi}, the models are only
proper, and semipositivity on the reduced fiber components need
not imply semipositivity on all of $\ker df$ at a singular point
of a fiber. On the strictly semistable fibration used in the proof, the
pulled-back curvature is semipositive on the full relative tangent space. On a dominating projective
model tower, the preceding construction makes the pulled-back
bundle relatively nef. The intersection inequality then descends
by the projection formula.

For a geometric model, the argument works directly on any
projective tower $\mathcal X\to\mathcal C\to T$: an ample
$\mathcal A$ on $\mathcal X$ makes
$\mathcal L+\epsilon\mathcal A$ relatively ample, and an ample
base twist makes
$\mathcal L+\epsilon\mathcal A+\mathcal f^*\mathcal N_\epsilon$
globally ample. Both factors are nef on the whole models, so the
induced adelic bundle is relatively nef. Proper towers admit
dominating projective towers with the same generic fibration;
no semistable alteration is needed here. The common models in
the approximation proof of Theorem~\ref{thm:ahi} serve to put
independently chosen approximations on compatible models.

\begin{remark}\label{rem:height-fundamental-inequality}
Over a number field $K$, let $\bar L$ be integrable with $L$ nef and
$L^d>0$. Write $e_{\mathrm{ess}}(\bar L)$ for the essential
minimum of its absolute point heights. If the metrics are
semipositive on all of $X$, then
\[
 \frac{\bar L^{d+1}}{(d+1)[K:\BQ]L^d}
 \le e_{\mathrm{ess}}(\bar L)
 \le \frac{a_A(\bar L)}{[K:\BQ]}.
\]
The first inequality follows from \cite[Theorem~5.2]{Zh95a}
by ample approximation. For the second, we apply the upper bound of
Qu and Yin \cite[Theorem~5.1]{QY24} on models and pass to the limit,
approximating the polarization by arithmetically ample bundles. Thus global semipositivity already
gives Theorem~\ref{thm:ahi}. For general metrics, the fundamental inequality
uses the $\chi$-arithmetic volume \cite[Introduction, (1.2)]{BPS15}.
The torsion bound allows Theorem~\ref{thm:hi} to retain the
self-intersection under relative positivity, including $L^d=0$.
\end{remark}

\subsubsection*{Zero fiber volume}

When both the total and generic fiber volumes vanish, strong
relative nefness gives the following equality.

\begin{corollary}
\label{cor:zero-volume-arithmetic}
Let $f:X\to C$ be as in Theorem~\ref{thm:ahi}, and let
$\bar L$ be strongly relatively nef over $C$. If
\[
 L^d=0,\qquad L_\eta^{d-1}=0,
\]
then $\bar L^{d+1}=0$.
\end{corollary}
\begin{proof}
We write $\bar L=\bar M-\bar U$, where
$\bar U=f^*\bar N$ and $\bar M,\bar N$ are nef adelic
bundles with $N$ ample. The adelic bundle
$\bar Q=\bar M+\bar U$ is nef, and
\[
 Q^d=(L+2f^*N)^d
     =L^d+2d\deg(N)L_\eta^{d-1}=0.
\]
It is also strongly relatively nef, since
$\bar Q+f^*\bar N$ is nef. Theorem~\ref{thm:ahi} gives
$\bar Q^{d+1}\le0$, while positivity of nef adelic
intersections gives the reverse inequality
\cite[Proposition 5.6]{YZ17}. Hence
\[
 0=\bar Q^{d+1}
  =\sum_{j=0}^{d+1}\binom{d+1}{j}
       (\bar M^{d+1-j}\bar U^j).
\]
Each intersection $\bar M^{d+1-j}\bar U^j$ is nonnegative,
so each is zero. Expanding $(\bar M-\bar U)^{d+1}$ now gives
$\bar L^{d+1}=0$.
\end{proof}

The following example shows why the corollary uses global strong
relative nefness. Let $K$ be a number field, let $E/K$ be an
elliptic curve, and let $M\in\Pic^0(E)(K)$ be non-torsion, with
its canonical flat adelic metric. On $X=\BP^1_K\times E$, we put
\[
 f=\operatorname{pr}_1,\qquad
 \bar L=f^*\overline{\mathcal O}(1)
                   +\operatorname{pr}_2^*\bar M,
\]
where $\overline{\mathcal O}(1)$ has semipositive metrics.
Then $L$ is nef, $L^2=\deg L_\eta=0$, and the metrics
on $L$ are semipositive at every place. Nevertheless the
projection formula and the Faltings--Hriljac theorem give
\[
 \bar L^3=3\deg\mathcal O(1)\,\bar M^2<0;
\]
see \cite[\S3.2]{YZ17}. Thus $\bar L$ is not strongly relatively
nef: local semipositivity does not suffice for the corollary.

\subsection{Examples and applications}\label{ss:height-consequences}
We begin with the Hodge index inequality of Faltings and Hriljac,
then derive a Bogomolov--Gieseker inequality and its application to
projective bundles. Finally, we construct examples by blowing up
points on surfaces.

\subsubsection*{The Hodge index and Bogomolov--Gieseker inequalities}

For $X=C$ over a number field, we recover the arithmetic Hodge
index inequality of Faltings and Hriljac \cite{Fa84,Hr85}.
Namely, an integrable adelic line bundle $\bar N$ of degree zero
on $C$ satisfies
\[
 \bar N^2\le0.
\]

Let $f:X\to C$ be as in Theorem~\ref{thm:ahi}, with $C$ smooth.
For an integrable adelic line bundle $\bar L$ on $X$, we put
\[
 \bar J=D_f(\bar L)=D_{X/C}(\bar L),
\]
where the Deligne pairing has $d$ factors. By
\cite[Theorem 4.1.3]{YZ26}, $\bar J$ is an integrable adelic
line bundle on $C$. We denote the numerical fiber class by $F$.

For an integrable adelic line bundle with $L_\eta^{d-1}>0$,
we define the relative discriminant by
\[
 \Delta_f(\bar L)
 =\frac{d+1}{2dL_\eta^{d-1}}\,
       \bar J^2-\bar L^{d+1}.
\]
For an integrable adelic line bundle $\bar N$ on $C$, the projection
formula gives
\begin{align*}
 D_f(\bar L+f^*\bar N)
    &=\bar J+dL_\eta^{d-1}\bar N,\\
 (\bar L+f^*\bar N)^{d+1}
    &=\bar L^{d+1}
      +(d+1)(\bar J\cdot\bar N)
      +\frac{d(d+1)}2 L_\eta^{d-1}\,
           \bar N^2.
\end{align*}
The cross terms cancel, so
\[
 \Delta_f(\bar L+f^*\bar N)=\Delta_f(\bar L).
\]

\begin{thm}[Bogomolov--Gieseker inequality]\label{thm:bgi}
Let $\bar L$ be relatively nef over $C$. Suppose that
$L_\eta^{d-1}>0$ and that the normalized class
\[
 L^\circ=L-\frac{L^d}{dL_\eta^{d-1}}F
\]
is nef on $X$. Then
\[
 \Delta_f(\bar L)\ge 0\]
\end{thm}
\begin{proof}
Invariance under base twists, applied to every
$\bar L+\epsilon\bar M$ with the same fixed nef $\bar M$, shows that
\[
 \bar L^\circ
 =\bar L-\frac{1}{dL_\eta^{d-1}}f^*\bar J
\]
is relatively nef. Since $\deg J=L^d$, its underlying
class is numerically $L^\circ$, which is nef and has
$(L^\circ)^d=0$. Theorem~\ref{thm:ahi} therefore gives
$(\bar L^\circ)^{d+1}\le0$. Multilinearity and the projection
formula give
\[
 (\bar L^\circ)^{d+1}
 =\bar L^{d+1}
  -\frac{d+1}{2dL_\eta^{d-1}}\,\bar J^2,
\]
which proves the inequality.
\end{proof}

To recover Moriwaki's inequality, let $\mathcal C$ be a regular
flat projective arithmetic surface with smooth generic curve $C$,
and let $\bar{\mathcal E}$ be a smooth Hermitian vector bundle of
rank $r$ on $\mathcal C$, with semistable generic restriction $E$.
On $X=\mathbb P_C(E)$, let $\bar\xi$ be the adelic bundle induced
by the tautological quotient bundle on
$\mathbb P_{\mathcal C}(\mathcal E)$ with its quotient metric.
An ample Hermitian base twist and a common constant rescaling
express $\bar\xi$ as a difference of nef bundles, so it is
strongly relatively nef over $C$.
Since $E$ is semistable, the curve criterion shows that the
normalized class $\xi-r^{-1}f^*\det E$ is nef. We can therefore
apply Theorem~\ref{thm:bgi} with $d=r$ and $L_\eta^{d-1}=1$.

For $\bar J=D_{\mathbb P_C(E)/C}(\bar\xi)$, the Bott--Chern
corrections in the comparison with Gillet--Soul\'e classes cancel
in the combination \cite[\S3.3]{G}
\[
 2r\Delta_f(\bar\xi)
 =(r+1)\bar J^2
       -2r\bar\xi^{r+1}
 =\widehat{\deg}\bigl(
      2r\widehat c_2^{\mathrm{GS}}(\bar{\mathcal E})
      -(r-1)\widehat c_1^{\mathrm{GS}}(\bar{\mathcal E})^2\bigr).
\]
Its nonnegativity recovers Moriwaki's arithmetic
Bogomolov--Gieseker inequality \cite{Mo94}.

\subsubsection*{Surface blow-ups}

By blowing up points on a curve, then resolving a pencil, one produces
nef, relatively ample bundles which are not semiample.
Let $S$ be a smooth geometrically integral
projective surface over an infinite field $k$, let $D$ be a very
ample line bundle on $S$, and let $Z\subset S$ be a smooth
geometrically integral curve.

\begin{lemma}\label{lem:integral-pencil}
For every integer $r\ge3$, there is a Lefschetz pencil in $|rD|$
defined over $k$, with transverse finite \'etale base locus
disjoint from $Z$, whose geometric members are integral.
\end{lemma}
\begin{proof}
By \cite[Th\'eor\`eme 2.5.1 and Corollaire 3.2.1]{Katz73},
a general pencil for the $r$-fold embedding is a Lefschetz pencil
in $|rD|$, with transverse finite \'etale base locus disjoint
from $Z$. These conditions define a nonempty open subset of a
Grassmannian. Since $k$ is infinite, we can choose the pencil
over $k$.

It remains to prove that its geometric members are integral.
They are connected because $rD$ is ample. If a member were
reducible, we could write it as two parts $G_1,G_2$ meeting at
its unique ordinary node, so $G_1\cdot G_2=1$. We put $b=D^2$
and $d_i=D\cdot G_i$. The Hodge index theorem gives
\[
 d_1d_2\le b,\qquad d_1+d_2=rb,
\]
On the other hand, the positive integers $d_1,d_2$ satisfy
$d_1d_2\ge rb-1>b$, a contradiction. A smooth member is integral
by connectedness, and the Lefschetz condition excludes nonreduced
members.
\end{proof}

For the following construction, let $(Y,A_Y)$ be a smooth
geometrically integral projective surface over an infinite field
$k$, with $A_Y$ very ample, and put $a=A_Y^2$.
We choose $m\in\mathbb Q_{\ge1}$ with $n=ma\in\mathbb Z$,
a smooth geometrically integral curve $B\in|A_Y|$, and distinct
points $P_1,\ldots,P_n\in B(k)$. Let $\pi:Y'\to Y$ be their
blow-up, with exceptional curves $E_i$, and let $\widetilde B$
be the strict transform of $B$. We put
\[
 L'=m\pi^*A_Y-\sum_{i=1}^nE_i.
\]
We interpret $L'$ as a rational line bundle when necessary.

The restriction
\[
 N=L'|_{\widetilde B}
   \simeq (mA_Y)|_B-\sum_{i=1}^nP_i
\]
has degree zero.

\begin{proposition}
\label{prop:polarized-surface-blowups}
Assume that $N$ is non-torsion after clearing denominators.
Then $L'$ is geometrically nef and not semiample, and
\[
 (L')^2=m(m-1)a,\qquad
 L'\cdot\widetilde B=0,\qquad L'\cdot E_i=1.
\]
Then there is a further blow-up $p:X\to Y'$ and
a flat projective fibration $f:X\to C=\BP^1_k$ with geometrically
integral fibers such that
\[
 L=p^*L'
\]
is nef, $f$-ample, and not semiample, with $L^2=m(m-1)a$.
The bundle $L$ is big if and only if $m>1$.
\end{proposition}

\begin{proof}
The blow-up formulas give the three displayed intersections.
To prove nefness, we consider a geometric integral curve
$\Gamma\subset Y$ different from $B$. Since $B$ is smooth,
we have
\[
 \sum_i\operatorname{mult}_{P_i}\Gamma\le B\cdot\Gamma
      =A_Y\cdot\Gamma.
\]
Consequently its strict transform satisfies
\[
 L'\cdot\widetilde\Gamma
   =mA_Y\cdot\Gamma-\sum_i\operatorname{mult}_{P_i}\Gamma
   \ge(m-1)A_Y\cdot\Gamma\ge0.
\]
Together with $L'\cdot E_i=1$ and $L'\cdot\widetilde B=0$,
this proves that $L'$ is nef. Its square $(L')^2=m(m-1)a$
then shows that $L'$ is big precisely when $m>1$.

To construct the fibration, we choose a very ample line bundle
$D$ on $Y'$ and an integer $r\ge3$. Applying
Lemma~\ref{lem:integral-pencil} with $Z=\widetilde B$, we obtain
a pencil in $|rD|$ with geometrically integral members and
transverse finite \'etale base locus disjoint from
$\widetilde B$. We blow up this base locus to obtain
$p:X\to Y'$ and
$f:X\to C=\BP^1_k$. The new exceptional curves are geometrically
sections, and every fiber is the strict transform of a pencil
member. The morphism is flat because $X$ is integral and $C$
is a smooth curve. For every geometric fiber $F$,
\[
 L\cdot F=r(L'\cdot D)>0.
\]
Here positivity follows from
$L'=\widetilde B+(m-1)\pi^*A_Y$ and ampleness of $D$.
Thus $L$ is ample on every geometric fiber, and hence $f$-ample
\cite[Tags 0B5Y and 0D2N]{Stacks}. Pullback preserves nefness
and the self-intersection.

Finally, the base locus is disjoint from $\widetilde B$, so
$L$ restricts to $N$ on the inverse image of $\widetilde B$.
Since a non-torsion line bundle of degree zero on a smooth
projective curve has no globally generated positive power,
neither $L'$ nor $L$ is semiample.
\end{proof}

For any ample line bundle $A_0$ on $X$, the normalization
$A=(A_0^2)^{-1/2}A_0$ gives a polarized fibration $(X/C,A)$
with $A^2=1$, together with the nef bundle $L$.

Over a number field, Lemma~\ref{lem:polarized-arithmetic-spread}
with parameter space $\Spec K$ gives a relatively ample rational
model of $L$ with positive vertical curvature. An ample Hermitian
base twist and a common rescaling give a strongly relatively nef
adelic metric $\bar L$. For a nef adelic polarization with
$A^2=1$, Theorem~\ref{thm:ahi} gives
\[
 \bar L^3\le3m(m-1)a(\bar L\cdot\bar A^2).
\]
In particular, the right-hand side is zero when $m=1$.

For $Y=\BP^2$, $A_Y=\CO(3)$, $m=s/3$, and $n=3s$, where
$s\ge3$ is an integer, the construction gives
\[
 L'=s\pi^*\CO_{\BP^2}(1)-\sum_{i=1}^{3s}E_i,
 \qquad L^2=s(s-3).
\]
The case $s=4$ is Zariski's example with twelve points
\cite[Example~6.15, p.~373]{Fuj83}. Over a number field, take
an elliptic curve $B$ embedded by $|3O|$ and points $P_i=iT$,
$1\le i\le3s$, for a point $T$ of infinite order. Then
$N=\CO_B(s)(-\sum_iP_i)$ corresponds to
$-3s(3s+1)T/2\in\Pic^0(B)$ and is non-torsion.

For $s=3$, we allow $N$ to be torsion. Then
$L'=\CO_{Y'}(\widetilde B)$ and
\[
 N=L'|_{\widetilde B}
   \simeq\CO_B(3)(-P_1-\cdots-P_9).
\]
The pencil construction of
Proposition~\ref{prop:polarized-surface-blowups} still gives $L=p^*L'$.

\begin{proposition}
\label{prop:blowups}\label{ex:nine-points}
The following are equivalent: $L'$ is semiample, $L$ is semiample,
and $N\in\Pic^0(B)$ is torsion.
\end{proposition}
\begin{proof}
If either $L'$ or $L$ is semiample, its restriction to
$\widetilde B$ or its inverse image in $X$ is semiample.
Since $\deg N=0$, this implies that $N$ is torsion.
If $N$ is non-torsion, the restriction sequences
\[
 0\longrightarrow (L')^{\otimes(r-1)}
   \longrightarrow (L')^{\otimes r}
   \longrightarrow i_*N^{\otimes r}\longrightarrow0,
 \qquad i:\widetilde B\hookrightarrow Y',
\]
give $h^0(Y',(L')^{\otimes r})=1$ for all $r\ge1$.

Conversely, suppose that $N$ has order $q$. The restriction
sequences above give
$h^0(Y',(L')^{\otimes(q-1)})=1$. Since $K_{Y'}=-L'$,
Serre duality gives
$h^2(Y',(L')^{\otimes(q-1)})=h^0(Y',(L')^{\otimes(-q)})=0$.
Riemann--Roch on the rational surface gives
$\chi(Y',(L')^{\otimes(q-1)})=1$, and hence
$h^1(Y',(L')^{\otimes(q-1)})=0$.
The restriction map
$H^0(Y',(L')^{\otimes q})\to H^0(B,N^{\otimes q})$
is therefore surjective. A section restricting to $1$ on
$\widetilde B$, together with the canonical section vanishing on
$q\widetilde B$, generates $(L')^{\otimes q}$ everywhere.
Its pullback generates $L^{\otimes q}$, proving the remaining
implication.
\end{proof}

In the setting of Theorem~\ref{prop:function-field-height}, let
$\mathcal C\to T$ be a regular projective surface model with
smooth geometrically integral generic curve $C/K$, and let
$\mathcal E$ be a vector bundle of rank $r$ on $\mathcal C$.
Assume that $E=\mathcal E_K$ is geometrically strongly semistable
of degree zero. In positive characteristic this means that all
iterated Frobenius pullbacks of $E_{\overline K}$ are semistable;
in characteristic zero ordinary geometric semistability suffices.
We put
\[
 \mathcal X=\BP_{\mathcal C}(\mathcal E),\qquad
 \mathcal L=\CO_{\mathcal X}(1),\qquad
 X=\BP_C(E),\quad L=\CO_X(1),
\]
using the quotient convention, and write
$\mathcal f:\mathcal X\to\mathcal C$ for the projection.
Let $\bar L$ be the function-field adelic bundle induced by
$\mathcal L$. A sufficiently ample twist from $\mathcal C$
makes $\mathcal L$ ample, so $\bar L$ is strongly relatively nef over
$C$. The underlying $L$ is nef, with
\[
 L^r=\deg E=0,\qquad L_\eta^{d-1}=1.
\]

We use the normalization factor $c_k>0$ defined before
Theorem~\ref{prop:function-field-height}.

\begin{example}
\label{ex:strongly-semistable-projective-bundle}
The adelic bundle $\bar L$ satisfies
\[
 \bar L^{r+1}=c_k\mathcal L^{r+1}
 =c_k\deg_{\mathcal C}\bigl(c_1(\mathcal E)^2-c_2(\mathcal E)\bigr)
 \le0.
\]
Moreover,
\[
 \deg_{\mathcal C}\bigl(2r c_2(\mathcal E)
                    -(r-1)c_1(\mathcal E)^2\bigr)\ge0.
\]
\end{example}

Nefness follows from the curve criterion \cite[Lemma 2.5]{FL22}.
Over $\overline K$, a vertical curve has positive $L$-degree.
The normalization of a horizontal curve gives a finite map
$q:D\to C$ and a quotient $q^*E\twoheadrightarrow L|_D$.
Strong semistability persists under finite pullback, including
inseparable pullback, so this quotient has nonnegative degree.

The projective-bundle formula gives the displayed expression for
$\bar L^{r+1}$. Theorem~\ref{prop:function-field-height} makes
it nonpositive, since $L$ is nef with $L^r=0$ and $\mathcal L$
is relatively ample over $\mathcal C$.
We let $\overline{\det E}$ be the adelic bundle induced by
$\det\mathcal E$ and apply Theorem~\ref{prop:function-field-height}
to $\bar L-r^{-1}f^*\overline{\det E}$. Its underlying bundle
is numerically equivalent to $L$, so the theorem gives the
displayed Bogomolov--Gieseker inequality.

Let $T$ be a smooth proper integral curve
over an algebraically closed field $k$, and let
$s_1,\ldots,s_9:T\to\BP^2$ be generically distinct morphisms
whose generic values lie on a smooth plane cubic. We write
$\Gamma_i$ for the graph of $s_i$. Assume that each fiber contains at most one
collision, involving two sections with intersection length one.
Let $p:\mathcal X\to\BP^2\times T$ be the successive blow-up
of the $\Gamma_i$, using strict transforms as centers, and let
$g:\mathcal X\to T$ be the resulting smooth projective morphism.

The generic surface is the $Y'$ of
Proposition~\ref{prop:polarized-surface-blowups} with
$Y=\BP^2_{k(T)}$, $A_Y=\CO(3)$, $m=1$, $n=9$; here the
restriction to the cubic may be torsion.
We write $E_i$ for the total transforms of the exceptional
divisors and put
\[
 H=p^*\operatorname{pr}_1^*\CO_{\BP^2}(1),\qquad
 \mathcal L=3H-\sum_{i=1}^9E_i
           =\omega_{\mathcal X/T}^{\otimes(-1)}.
\]
Let $\bar L$ be the adelic bundle on the generic surface over
$K=k(T)$ induced by $\mathcal L$. Its underlying bundle satisfies
$L^2=0$. Collisions give infinitely near points on the surface
fibers. Put
\[
 d_i=\deg s_i^*\CO_{\BP^2}(1),\qquad
 N_{\rm coll}=\sum_{i<j}\operatorname{length}_k
                         (\Gamma_i\cap\Gamma_j).
\]

\begin{example}
\label{ex:moving-nine-points}
If $\mathcal L$ is nef on every fiber of $g$, its induced
metrics are semipositive at every place, and its adelic height is
\begin{equation}\label{eq:moving-nine-points-height}
 \bar L^3=\mathcal L^3=2N_{\rm coll}-6\sum_{i=1}^9d_i\le0.
\end{equation}
In particular, the number of collision fibers is at most
$3\sum_i d_i$.
\end{example}

\begin{proof}
We first prove $\mathcal L^3\le0$ by applying the argument of
Theorem~\ref{prop:function-field-height} directly to $g$.
Since $\mathcal L$ is nef on every fiber of $g$, an ample
perturbation $\mathcal L+\epsilon\mathcal M$ is $g$-ample.
Relative Serre vanishing makes the higher direct images of
sufficiently high powers zero. The Riemann--Roch comparison in
that proof then gives $\mathcal L^3\le0$, since $L^2=0$.

We next compute $\mathcal L^3$ one blow-up at a time. Let
$q:\mathcal Y'\to\mathcal Y$ be a blow-up with center the section
$\Gamma$ and exceptional divisor $E_\Gamma$. We put
$\mathcal A=\omega_{\mathcal Y/T}^{\otimes(-1)}$ and
$\ell=\deg(\mathcal A|_\Gamma)$.
Since $\Gamma$ is a section of a smooth surface fibration,
its normal bundle has degree $\ell$. The projective-bundle
formula on $E_\Gamma$ gives
\[
 (q^*\mathcal A)^2\cdot E_\Gamma=0,\qquad
 q^*\mathcal A\cdot E_\Gamma^2=-\ell,\qquad
 E_\Gamma^3=-\ell.
\]
Consequently
\[
 (q^*\mathcal A-E_\Gamma)^3=\mathcal A^3-2\ell.
\]
At the $j$-th step, each simple intersection with an earlier
section subtracts one from the anticanonical degree on the
lifted section. Thus
\[
 \ell=3d_j-\sum_{i<j}\operatorname{length}_k
                          (\Gamma_i\cap\Gamma_j).
\]
Starting from $(\operatorname{pr}_1^*\CO_{\BP^2}(3))^3=0$
and summing over the nine blowups proves
\eqref{eq:moving-nine-points-height}.
\end{proof}

For example, $\mathcal L$ is fiberwise nef when
$s_1,\ldots,s_9$ factor through a fixed smooth plane cubic
$E\subset\BP^2$. The strict transform of $E\times T$ then
represents $\mathcal L$. Its restriction to every fiber is a
smooth integral curve of self-intersection zero: at a collision,
the infinitely near point lies in the tangent direction of $E$.
Such a curve is nef.

For a concrete example, we take $T=E$, let $s_1,\ldots,s_8$
be distinct constant maps, and let $s_9:E\hookrightarrow\BP^2$
be the inclusion. The graph of $s_9$ meets each constant section
once and simply, so
$N_{\rm coll}=8$, $\sum_i d_i=3$, and $\bar L^3=\mathcal L^3=-2$.
For the constant points $Q_1,\ldots,Q_8$ and the generic point
$\eta$ of $T=E$, the restriction of $L$ to the strict cubic is
$\CO_E(3)(-Q_1-\cdots-Q_8-\eta)$. Its class varies as a translate
of $-\mathrm{id}_E$ in $\Pic^0(E)$, so no positive multiple
is zero. Thus $L$ is not semiample by
Proposition~\ref{prop:polarized-surface-blowups}.

\section{Energies of fibrations}\label{sec:energy}

For a polarized fibration $(X/C,\bar A)$, we define an energy
$\Psi_A$, deduce lower bounds from the height inequality, and study
its critical points. When $L^d=0$, the energy is independent of
$\bar A$ and is related to classical functionals. We propose conjectures on the 
boundedness of energy and semipositivity of minimizing and critical metrics.

\subsection{Local metrics, energy, and conjectures}
\label{ss:local-metrics}

For a polarized fibration $(X/C,\bar A)$ over a valued field $F$,
we define an energy $\Psi_A$ and conjecture its boundedness. For
relatively ample $L$, its infimum is the same on $\CH^+(L)$
and $\BH^+(L)$. We also prove that the energy between strongly
relatively nef metrics vanishes when $L^d=L_\eta^{d-1}=0$.

\subsubsection*{The local energy}

Let $f:X\to C$ be a flat projective surjective morphism from a smooth
projective variety of dimension $d$ to a smooth projective curve over
a field $F$ with a real-valued absolute value. We work over the
completion of $F$ and give $A$ a semipositive metric with $A^d=1$.
For continuous integrable metrics $h_0,h_1$ on $L$ put $\alpha_i=c_1(L,h_i)$.
Their local energy is
\begin{equation}\label{eq:H-energy}
 \Psi_A(h_0,h_1)=\int_{X^{\rm an}}\log\frac{|\cdot|_{h_1}}{|\cdot|_{h_0}}
 \left(\sum_{i=0}^d\alpha_0^i\alpha_1^{d-i}-(d+1)L^d c_1(\bar A)^d\right).
\end{equation}
The measure in parentheses has mass zero, so a constant rescaling of
either metric preserves the energy. If $L^d=0$, the energy is independent of $\bar A$.
As a difference of weights on $D_A(L)$, the energy satisfies
\[
 \Psi_A(h_0,h_2)=\Psi_A(h_0,h_1)+\Psi_A(h_1,h_2).
\]

We now assume that $L$ is $f$-ample and that $F$ is $\mathbb C$ or
a complete discretely valued field. We allow singular fibers of $f$.
We denote by $\CH^+(L)$ the following space of metrics:
\begin{itemize}
\item over $\BC$, smooth metrics whose curvature is strictly
positive on each relative tangent space
$T_{X/C,x}=\ker(df_x)$, including at singular points of the fibers;
\item over a nonarchimedean field, rational model metrics
induced by $\mathcal f$-ample $\mathbb Q$-line bundles
$\mathcal L$ extending $L$ on projective model towers
$\mathcal f:\mathcal X\to\mathcal C$ over $\mathcal O_F$.
\end{itemize}

\subsubsection*{Larger metric spaces}

To allow degenerate fiberwise curvature, we introduce continuous
relatively nef metrics. Lemma~\ref{lem:relative-nef-approximation}
will identify their energy infimum with that on $\CH^+(L)$.

Over $\BC$, we call a continuous metric semipositive if its local
weights are plurisubharmonic. Over a nonarchimedean field, we call
it semipositive if it is a uniform limit of nef model metrics on
the same bundle. The latter notion is called \emph{strong
semipositivity} in \cite[\S3.6.6]{YZ26}. We will also use the
notion of local nefness obtained by perturbing the bundle.
For the definitions, let $f:X\to C$ be any projective morphism
of projective varieties over $F$.
A continuous metrized line bundle $(L,h)$ is \emph{strongly relatively nef}
for $f$ if some ample bundle $N$ on $C$ has a continuous
semipositive metric $n$ such that
\begin{equation}\label{eq:admissible-base-twist}
 \bar M=(L,h)+f^*(N,n)
 \quad\hbox{is semipositive on }X.
\end{equation}

The pair $(N,n)$ in \eqref{eq:admissible-base-twist} may depend
on $h$. When approximating $h$, we will keep $(N,n)$ fixed.
Since $f^*(N,n)$ is constant on each fiber, $(L,h)$ is integrable
and its restriction to every fiber is semipositive, including at
singular fibers. If $L$ is relatively ample, we may add a further
ample semipositive bundle on $C$ and assume that $M$ is ample.

An integrable continuous metrized line bundle $\bar L$ is
\emph{relatively nef} if $\bar L+\epsilon\bar A'$ is
strongly relatively nef for every rational $\epsilon>0$, with one fixed
globally semipositive metrized $\mathbb Q$-line bundle $\bar A'$.
The choice of $\bar A'$ may depend on $h$. For a strongly relatively
nef metric we may take $\bar A'=0$.

Both notions are preserved under base change, by the argument in
the \hyperref[par:relative-nef-base-change]{base-change discussion}.
For a morphism to $\Spec F$, strong relative nefness is
semipositivity. Relative nefness means that
$\bar L+\epsilon\bar A'$ is semipositive for every rational
$\epsilon>0$, for one fixed semipositive $\bar A'$. Over $\BC$,
this is equivalent to semipositivity, since limits of positive
currents are positive. In the nonarchimedean case, it is local
nefness, also called semipositivity in
\cite[\S3.6.6 and \S4.6.2]{YZ26}, distinct from uniform
approximation by nef model metrics on a fixed bundle. Restriction
to a fiber therefore gives a locally nef metric and a nef
underlying bundle: for every fiber curve $\Gamma$,
$\deg_\Gamma(L+\epsilon A')\ge0$ implies $\deg_\Gamma L\ge0$.

We return to the fibration $f:X\to C$ of \eqref{eq:H-energy}.

We denote by $\BH^+(L)$ the space of continuous relatively nef
metrics on $L$; these are integrable by definition.
For $f$-ample $L$, choose an ample bundle $N$ on $C$ such that
$M=L+f^*N$ is ample. Projective embeddings of suitable powers
give Fubini--Study metrics on $M$ and $N$, induced by the
standard squared Euclidean norm over $\BC$ and the maximum
norm over a nonarchimedean field. Their difference
$\bar L=\bar M-f^*\bar N$ belongs to $\CH^+(L)$: over
$\BC$ its curvature is strictly positive on $\ker(df_x)$,
including at singular fiber points; in the nonarchimedean case,
the graph construction in Lemma~\ref{lem:relative-nef-regularization}
gives a relatively ample model representing it.

For $f$-semiample $L$, choose $m>0$ and an ample bundle $M$
on $C$ so that $Q=L^{\otimes m}\otimes f^*M$ is globally
generated. Semipositive Fubini--Study metrics on $Q$ and $M$
define a strongly relatively nef metric
\[
 (L,h)=m^{-1}(\bar Q-f^*\bar M).
\]
It is smooth over $\BC$ and a rational model metric in the
nonarchimedean case. Its restriction to every fiber is semipositive.

Over $\BC$, a smooth metric whose curvature is strictly positive
on every nonzero vector in the relative tangent space $T_{X/C,x}$
is strongly relatively nef, even when $f$ has singular fibers.
To see this, we fix a smooth positive metric on an ample bundle
on $C$, denote its curvature by $\omega_C$, and put
$\alpha=c_1(L,h)$.
The form
$\alpha+a f^*\omega_C$ is positive for sufficiently large $a$.
Otherwise, compactness would give unit tangent vectors $v_j$ at
points $x_j$, with $a_j\to\infty$ and
$(\alpha+a_j f^*\omega_C)(v_j,\bar v_j)\le0$, converging to a
unit vector $v$ with $df(v)=0$ and $\alpha(v,\bar v)\le0$.
This contradicts strict positivity on $\ker(df_x)$. We may take
$a$ to be an integer and obtain the required base twist.
Over a nonarchimedean field, we start with an
$\mathcal f$-ample bundle $\mathcal L$ and choose an ample
bundle $\mathcal N$ on $\mathcal C$. For $a\gg0$, the bundle
$\mathcal L+a\mathcal f^*\mathcal N$ is ample
\cite[Tag 0892]{Stacks}. The model metrics on $\mathcal L+a\mathcal f^*\mathcal N$ and
$a\mathcal N$ are semipositive, so the metric on $\mathcal L$
is strongly relatively nef. We have proved
$\CH^+(L)\subset\BH^+(L)$ in both cases.

\paragraph*{Approximation of relatively nef metrics.}
For the next two lemmas, let $L$ be relatively ample and $C$ smooth,
with $f$ allowed to have singular fibers. Over $\BC$ we assume
that $X$ is smooth. In the nonarchimedean case, we allow any
complete nontrivially valued field $F$, singular $X$, and a
$\BQ$-line bundle $L$. The semipositive metric on $A$ is fixed.

\begin{lem}
\label{lem:relative-nef-regularization}
Every strongly relatively nef metric $h$ on $L$ is a uniform
limit of metrics $h_j\in\CH^+(L)$, with
\[
 \Psi_A(h_0,h_j)\longrightarrow \Psi_A(h_0,h)
\]
for every fixed integrable reference metric $h_0$.
If the valuation is discrete, the approximating models may be
chosen projective and flat over projective models of $C$.
\end{lem}
\begin{proof}
There is a decomposition $(L,h)=(M,a)-f^*(N,n)$ with $M,N$
ample and $a,n$ continuous semipositive. The approximation keeps
$M,N$ fixed.

Over $\BC$, the metrics $a$ and $n$ are uniform limits of
smooth semipositive metrics by \cite[\S2.1.1]{YZ26}. By taking
convex combinations with fixed smooth positive metrics, we may
choose smooth positive approximants $a_j\to a$ and $n_j\to n$.
Over a nonarchimedean field, we use the proof of
\cite[Theorem 7.8]{BE21} to approximate $a$ uniformly by ample
$\mathbb Q$-model metrics $a_j$ on $M$, and choose nef model
metrics $n_j\to n$ on $N$. In either case,
\[
 h_j=a_j-f^*n_j\longrightarrow h
\]
uniformly, and $h_j$ is strongly relatively nef and fiberwise
semipositive. Over $\BC$, the curvature of the base factor
vanishes on $\ker(df_x)$, so $h_j$ is smooth and strictly positive
there, including at singular points of the fibers. Thus $h_j\in\CH^+(L)$.

To prove the same inclusion over a nonarchimedean field, we need
a relatively ample model representing $h_j$. Let
$(\mathcal X_j,\mathcal M_j)$ and
$(\mathcal C_j,\mathcal N_j)$ represent $a_j$ and $n_j$.
We take the closure $\mathcal Z_j$ of the graph of $f$ in
$\mathcal X_j\times_{\mathcal O_F}\mathcal C_j$, and denote
the projections by $p_j$ and $\mathcal f_j$. The bundle
\[
 \mathcal L_j=p_j^*\mathcal M_j-\mathcal f_j^*\mathcal N_j
\]
extends $L$ and induces $h_j$. It is $\mathcal f_j$-ample,
since $p_j^*\mathcal M_j$ is the restriction of a relatively
ample bundle on the product and twisting from $\mathcal C_j$
preserves relative ampleness. Hence $h_j\in\CH^+(L)$.
If $F$ is discretely valued, we may also make the model morphism
flat. Since $f$ is flat, we modify $\mathcal C_j$ and replace
$\mathcal Z_j$ by its strict transform as in the proof of
\cite[Proposition 5.13]{BE21}. The strict transform is closed in
the base change, so the pullback of $\mathcal L_j$ remains
relatively ample and induces the same metric on $L$.

Finally, $h_j=a_j-f^*n_j$ is a difference of semipositive metrics
on fixed bundles. Expanding the mixed local intersections in
\eqref{eq:H-energy}, we can apply their continuity under uniform
convergence of the semipositive factors. The linear term against
$c_1(\bar A)^d$ converges by uniform convergence.
We obtain
\[
 \Psi_A(h_0,h_j)\longrightarrow \Psi_A(h_0,h)
\]
for every fixed integrable reference $h_0$; see
\cite[\S4.2.2]{YZ26} and \cite[Theorem 8.16]{BE21}.
\end{proof}

We next pass from strongly relatively nef metrics to relatively
nef metrics while keeping the underlying bundle $L$ fixed.

\begin{lem}
\label{lem:relative-nef-approximation}
Every relatively nef metric $h$ is a uniform limit of
metrics $h_j\in\CH^+(L)$ on the same bundle, with
$\Psi_A(h_0,h_j)\to \Psi_A(h_0,h)$ for every fixed integrable reference
metric $h_0$. Consequently,
\[
 \inf_{h\in\BH^+(L)}\Psi_A(h_0,h)
 =\inf_{h\in\CH^+(L)}\Psi_A(h_0,h).
\]
If the valuation is discrete, the approximating models may be
chosen projective and flat over projective models of $C$.
\end{lem}
\begin{proof}
We choose $\bar A'$ as in the definition of relative nefness.
Since relative ampleness is open, we can choose a rational
$\delta>0$ such that $L-\delta A'$ is $f$-ample. Applying the
Fubini--Study construction to $L-\delta A'$, we obtain a strongly
relatively nef metric $\bar M$. We put
$\bar L_* =\bar M+\delta\bar A'$ and, for rational
$\epsilon>0$, define
\[
 \begin{aligned}
 \bar L_\epsilon
 &=\frac{\delta}{\delta+\epsilon}(\bar L+\epsilon\bar A')
   +\frac{\epsilon}{\delta+\epsilon}\bar M\\
 &=\frac{\delta}{\delta+\epsilon}\bar L
   +\frac{\epsilon}{\delta+\epsilon}\bar L_*.
 \end{aligned}
\]
In the first expression, both terms are strongly relatively nef.
Combining their semipositive ample base twists and using convexity,
we see that $\bar L_\epsilon$ is strongly relatively nef. The
second expression shows that $\bar L_\epsilon$ has underlying
bundle $L$ and converges uniformly to $\bar L$. Since $\bar L$
and $\bar L_*$ are fixed and integrable, multilinearity of the
local pairings expresses
the energy as a polynomial in $\epsilon/(\delta+\epsilon)$.
Its value therefore converges to the energy of $h$.

We now apply Lemma~\ref{lem:relative-nef-regularization} to each
$\bar L_\epsilon$ and choose a diagonal sequence in $\CH^+(L)$
for one fixed reference metric. The cocycle identity gives energy
convergence for every other reference metric. The asserted equality
of infima follows from this approximation and the inclusion
$\CH^+(L)\subset\BH^+(L)$.
\end{proof}

Among the fibers $X_s/C_s$ satisfying the hypotheses of
Lemma~\ref{lem:relative-nef-approximation}, the set of parameters $s$
where the energy is unbounded below is the same for $\CH^+(L_s)$
and $\BH^+(L_s)$.

If $f$ is smooth over $\BC$, we may also take the infimum over
all smooth fiberwise-semipositive metrics. Indeed, convex
interpolation with a fixed smooth metric of strictly positive
vertical curvature approximates any such metric by metrics in
$\CH^+(L)$, with convergence of energy. The smoothness assumption
on $f$ ensures that positivity along the fibers controls the whole
relative tangent space. For singular fibers, positivity on
resolutions of the components need not control $\ker(df_x)$.

\subsubsection*{A lower bound for the local energy}

\begin{conj}\label{conj:H-energy-boundedness}
Suppose that $X$ and $C$ are geometrically integral and that $L$
is nef and $f$-ample. For every integrable reference metric $h_0$ on $L$,
\[
 \inf_{h\in\CH^+(L)}\Psi_A(h_0,h)>-\infty.
\]
\end{conj}
By Lemma~\ref{lem:relative-nef-approximation}, the same bound holds
on $\BH^+(L)$ if it holds on $\CH^+(L)$. This is the local finiteness
assertion of Conjecture~\ref{conj:envelope} for $S=\Spec F$.

\paragraph*{First variation.}

For a smooth reference
$h_0$ on $L$ and a smooth path $h_t$ over $\mathbb C$, we write
$\psi_t=\log(|\cdot|_{h_t}/|\cdot|_{h_0})$.
Differentiating \eqref{eq:H-energy} and integrating by parts, we obtain
\[
 \frac{d}{dt}\Psi_A(h_0,h_t)
 =(d+1)\int_{X^{\mathrm{an}}}\dot\psi_t
       \bigl(c_1(L,h_t)^d-L^dc_1(\bar A)^d\bigr).
\]
Over a nonarchimedean field, we make the corresponding calculation
for $\mathbb Q$-model metrics, also called modelic metrics.
Given $h_0,h_1$ on a common model, we put
\[
 \varphi=\log\frac{|\cdot|_{h_1}}{|\cdot|_{h_0}},
 \qquad
 |\cdot|_{h_t}=|\cdot|_{h_0}e^{t\varphi}.
\]
For rational $t$, the metric $h_t$ is a $\mathbb Q$-model metric.
By multilinearity of the Deligne pairing and linearity of the
$\bar A$-term, $\Psi_A(h_0,h_t)$ extends to a real polynomial in $t$.
Its derivative satisfies the same formula, with
$\dot\psi_t=\varphi$ and the products interpreted as signed
Chambert--Loir measures; see \cite[Theorem 8.16]{BE21}.
If $h_0$ and $h_1$ are strongly relatively nef, we choose a common
semipositive base twist. Convexity then shows that $h_t$ is
strongly relatively nef for $0\le t\le1$.

The energy vanishes under the following two intersection conditions.
For the next proposition, let $L$ be any line bundle in the local
setting. We allow $F=\BC$ or any complete nontrivially valued
nonarchimedean field, and retain the flat projective morphism $f$
with possibly singular fibers. Strong relative nefness is understood
as in \eqref{eq:admissible-base-twist}.

\begin{proposition}\label{prop:zero-volume-energy}
Assume that
\[
 L^d=0,\qquad L_\eta^{d-1}=0.
\]
Then, for all continuous strongly relatively nef metrics $h_0,h_1$
on $L$,
\[
 c_1(L,h_0)^i c_1(L,h_1)^{d-i}=0\quad(0\le i\le d),
 \qquad \Psi_A(h_0,h_1)=0.
\]
If $F=\BC$ and $f$ is smooth, the same conclusions hold for
every two smooth fiberwise-semipositive metrics.
\end{proposition}
\begin{proof}
We add semipositive ample base twists for $h_0$ and $h_1$ to
obtain a common twist $\bar U=f^*\bar N$ such that
\[
 \bar M_i=(L,h_i)+\bar U\quad(i=0,1)
\]
are semipositive. Their common underlying bundle is
$M=L+U$. Since $U^2=0$,
\[
 M^d=0,\qquad
 M^{d-1}U=\deg(N)L_\eta^{d-1}=0,\qquad
 M^{d-k}U^k=0\quad(2\le k\le d).
\]
Every mixed measure of top degree in
$\bar M_0,\bar M_1,\bar U$ is positive and has one of these
intersection numbers as its mass. Hence all these measures vanish.
We use Bedford--Taylor measures over $\BC$ and Chambert--Loir
measures in the nonarchimedean case. Their positivity,
multilinearity, and mass formula extend to continuous semipositive
metrics by approximation
\cite[Appendix 1, ``Local intersections'']{YZ17};
see also \cite[\S1.2.8 and \S\S1.3.11--1.3.12]{CL10}.
Expanding $c_1(L,h_i)=c_1(\bar M_i)-c_1(\bar U)$, we obtain
the asserted vanishing of mixed measures. Substitution in
\eqref{eq:H-energy} gives $\Psi_A(h_0,h_1)=0$.

We next assume that $f$ is smooth over $\BC$ and prove the
assertion for smooth fiberwise-semipositive metrics. Let
$\alpha=c_1(L,h)$ and let $V$ be its vertical curvature matrix.
The matrix $V$ is positive semidefinite. The top exterior power
of the vertical curvature is nonnegative and has fiber integral
$L_\eta^{d-1}=0$, so $\det V=0$ everywhere. The full curvature matrix
\[
 \begin{pmatrix}V&b\\b^*&c\end{pmatrix}
\]
has nonpositive determinant. Indeed, if $V$ has corank one,
the determinant of the full matrix is $-|b_0|^2$ times the
product of the positive eigenvalues of $V$, where $b_0$ is the
component in $\ker V$. If $V$ has larger corank, the full
determinant is zero.
Thus $\alpha^d\le0$. Since $\int_X\alpha^d=L^d=0$, this smooth
top-degree form vanishes pointwise. The convex path between two
fiberwise semipositive metrics has the same property. Expanding
$((1-t)\alpha_0+t\alpha_1)^d=0$ for $0\le t\le1$ shows that
all mixed top-degree forms vanish, and hence so does the energy.
\end{proof}

This is the local counterpart of
Corollary~\ref{cor:zero-volume-arithmetic}. Since
$L_\eta^{d-1}=0$ excludes relative ampleness, the approximation
lemma does not extend the proposition to all relatively nef metrics.

\subsection{Geometric examples and classical functionals}
\label{ss:metric-examples}

We first give examples from semiample bundles and bundles on
blow-ups of surfaces. We then compare our energy with Donaldson's
functional and the fibration functional of Feng--Liu--Wan.

\subsubsection*{Semiample line bundles}
Suppose that $L^{\otimes m}=g^*J$ for a morphism
$g:X\to Y$ with $\dim Y=d-1$ and $J$ ample. A
Fubini--Study metric on $J$ induces a smooth semipositive metric
$h_0$ on $L$ with $c_1(L,h_0)^d=0$. All mixed top-degree
products vanish for metrics pulled back from $Y$, so their energy
is zero. The same dimension argument gives zero inverse Deligne
height for adelic metrics descending through a compatible model
of $g$.

If $f$ is smooth and $h_0$ has strictly positive vertical
curvature, \cite[Theorem 0.1]{FLW} and
\eqref{eq:FLW-normalization} give $\Psi_A(h_0,h)\ge0$ for smooth
vertically positive $h$. Interpolation with $h_0$ extends the
inequality to smooth vertically semipositive metrics.

\subsubsection*{Big and nef examples}

For the surface in Proposition~\ref{prop:polarized-surface-blowups}
over a number field, let $m>1$ and $v=L^2=m(m-1)a$.
Theorem~\ref{thm:ahi} gives
\[
 \bar L^3\le3v(\bar L\cdot\bar A^2)
\]
for every integrable relatively nef adelic metric $\bar L$ and
normalized nef adelic polarization $\bar A$. A strongly relatively
nef reference $\bar L_0$ is obtained by choosing nef adelic
metrics on ample bundles $L+f^*N$ and $N$. At every place $w$,
Corollary~\ref{cor:H-local-bound} then gives
\[
 \inf_{h\in\BH^+(L_w)}\Psi_w(h_{0,w},h)
 \ge \bar L_0^3-3v(\bar L_0\cdot\bar A^2)>-\infty.
\]
Thus the corrected energy is bounded below although $L$ is not
semiample. Theorem~\ref{thm:main-energy} gives the analogous
bound for almost every realization over a finitely generated field.
For smooth globally semipositive metrics on the classical
blow-up at twelve complex points, see
\cite[Theorem~1.1 and Example~4.3]{Koike15}.

Another example is $X=\BP_C(N\oplus J)$, with
$N\in\Pic^0(C)$ non-torsion and $\deg J>0$. The quotient
bundle $L=\CO_X(1)$ is nef and relatively ample, with
$L^2=\deg J$, but its restriction to the section defined by
$N$ prevents semiampleness. The same local bound applies when
the data are defined over a number field.

\subsubsection*{Classical functionals}

Assume here that $f$ is smooth over $\BC$ and $L^d=0$.
Feng--Liu--Wan \cite{FLW} define the corresponding fibration
functional and prove that geodesic--Einstein metrics minimize it.
For projective bundles, the comparison is with Donaldson's
functional \cite{Do85}, whose arithmetic interpretation is due to
Moriwaki \cite[\S6]{Mo94}. The following identities fix the
normalizations.

\smallskip\noindent\emph{The fibration functional.}
We first compare normalizations for $h_0,h_1\in\CH^+(L)$
with strictly positive vertical curvature. Locally, we write
$|s|_h=e^{-\phi}$, so that
$\alpha_h=(\sqrt{-1}/2\pi)\partial\bar\partial\phi$, and put
$\Omega_i=\sqrt{-1}\partial\bar\partial\phi_i=2\pi\alpha_i$.
In the notation of Feng--Liu--Wan, the base dimension is $m=1$
and the fiber dimension is $n=d-1$. Since $L^d=0$, their
constant $\lambda$ is zero. Substituting in equations
(1.13)--(1.14) of \cite{FLW}, we obtain
\begin{equation}\label{eq:FLW-normalization}
\begin{aligned}
 \mathcal L_{\rm FLW}(\phi_1,\phi_0)
 &=-\frac{1}{d(d+1)}\int_X(\phi_1-\phi_0)
       \sum_{i=0}^d\Omega_0^i\Omega_1^{d-i}\\
 &=\frac{(2\pi)^d}{d(d+1)}\Psi_A(h_0,h_1).
\end{aligned}
\end{equation}
The positive factor identifies critical points and minima;
$\alpha^d=0$ is the geodesic--Einstein equation of slope zero
\cite[Proposition 1.4 and Theorem 0.1]{FLW}.

\smallskip\noindent\emph{Donaldson's functional.}
We now let $X=\mathbb P(E)$, $r=\rk E=d$, and $\deg E=0$.
For a smooth path $H_t$, $0\le t\le1$, of Hermitian squared
metrics on $E$, let $h_t$ be the induced tautological quotient
metric. We put
\[
 v_t=H_t^{-1}\dot H_t,\qquad
 R_t=\frac{\sqrt{-1}}{2\pi}
           \bar\partial(H_t^{-1}\partial H_t).
\]
\Needspace{9\baselineskip}
We normalize the degree-zero Donaldson functional and its
determinant term using the Chern curvature $R_t$ as follows:
\begin{align*}
 \mathcal D_0(H_0,H_1)
 &=\int_0^1\!\int_C
     \operatorname{Tr}\!\left(v_t
          \left(R_t-\frac{\operatorname{Tr}R_t}{r}\operatorname{Id}\right)
                    \right)\,dt,\\
 \mathcal D_{\det}(H_0,H_1)
 &=\int_0^1\!\int_C
            \operatorname{Tr}(v_t)\operatorname{Tr}(R_t)\,dt.
\end{align*}
By integration by parts, both integrals depend only on $H_0$
and $H_1$. They give the Chern--Simons description of Donaldson's
functional \cite{Do85,Mo94}. With these normalizations, we have
\begin{equation}\label{eq:Donaldson-normalization}
 \Psi_A(h_0,h_1)=\mathcal D_0(H_0,H_1)
                 +\frac{r+1}{r}\mathcal D_{\det}(H_0,H_1).
\end{equation}
To prove \eqref{eq:Donaldson-normalization}, let $\Pi$ be the
rank-one orthogonal projector corresponding to a quotient in a
fiber of $\mathbb P(E)$. We normalize the Fubini--Study form by
$\int_{\mathbb P(E_x)}\omega_{\rm FS}^{r-1}=1$ and use the
moment identity
\[
 \int_{\mathbb P(E_x)}\operatorname{Tr}(S\Pi)
              \operatorname{Tr}(T\Pi)\,\omega_{\rm FS}^{r-1}
 =\frac{\operatorname{Tr}(ST)+\operatorname{Tr}(S)\operatorname{Tr}(T)}
              {r(r+1)}.
\]
Since $\dot\phi_t=-\operatorname{Tr}(v_t\Pi)$, the horizontal
curvature of the quotient metric and this identity give
\[
 f_*(\dot\phi_t\alpha_t^r)
 =-\frac{\operatorname{Tr}(v_tR_t)
                  +\operatorname{Tr}(v_t)\operatorname{Tr}(R_t)}{r+1}.
\]
The first variation of $\Psi_A$ is
$-(r+1)\int_X\dot\phi_t\alpha_t^r$, which proves
\eqref{eq:Donaldson-normalization} by integration in $t$.
Equivalently, writing
$u=\log(\det H_1/\det H_0)$, the determinant term is
\[
 \mathcal D_{\det}(H_0,H_1)
 =\frac12\int_C u\bigl(c_1(\det E,\det H_0)
                         +c_1(\det E,\det H_1)\bigr).
\]
If the determinant metric is fixed, we choose a path preserving
it and obtain $\mathcal D_{\det}=0$ and $\mathcal D_0=\Psi_A$.
For a semistable bundle, Donaldson's functional is bounded below;
attainment by a Hermitian--Einstein metric requires polystability.

For a degree-zero semistable bundle $E$ over a number field,
$L=\CO_{\mathbb P(E)}(1)$ is nef and relatively ample with
$L^r=0$; see Example~\ref{ex:strongly-semistable-projective-bundle}.
A strongly relatively nef adelic reference and
Corollary~\ref{cor:H-local-bound} give a lower bound on
$\BH^+(L_v)$ at every place, recovering Donaldson's lower
bound for Hermitian norms and also treating metrics not induced
by such norms. The classical proof uses the Jordan--H\"older
filtration and Hermitian--Einstein metrics on stable factors
\cite{NS65,Do83,Do85,UY86}. Under long-time smooth existence of
their flow, Wan--Wang \cite[Theorem 2.8]{Wan} relate boundedness
to approximate geodesic--Einstein structures and nonlinear
semistability. In Part~II \cite[Proposition~2.2.1]{PartII} we prove
for a semistable projective bundle over a smooth
complex polarized variety that the infimum over all continuous
relatively nef metrics on the normalized tautological bundle equals
the infimum over Hermitian quotient metrics.

\smallskip\noindent\emph{Approximate criticality.}
For nef, $f$-ample $L$ with $L^d=0$, approximate criticality
also follows directly from nefness. If $B$ is ample on $C$ with
positive curvature $\omega_B$, choose
$\Omega_\varepsilon>0$ in the ample class
$c_1(L)+\varepsilon f^*c_1(B)$ and put
$\alpha_\varepsilon=\Omega_\varepsilon-\varepsilon f^*\omega_B
=c_1(L,h_\varepsilon)$. The metrics $h_\varepsilon$ have
positive vertical curvature, and $(f^*\omega_B)^2=0$ gives
\[
 \int_X|\alpha_\varepsilon^d|
 \le 2d\varepsilon\,L^{d-1}\cdot f^*B\longrightarrow0.
\]
Finski \cite[Theorem 1.1]{Finski26WZW} computes the infimum of
this curvature integral more generally using asymptotic
Harder--Narasimhan slopes; see also his work on fibered
Yang--Mills functionals and slopes of direct images
\cite{Finski26YM,Finski26HN}.

\subsection{Arithmetic lower bounds}\label{ss:energy}

Let $f:X\to C$ satisfy the hypotheses of Theorem~\ref{thm:ahi} and suppose that
$L$ is nef on $X$. On integral projective varieties we define the energy
by local intersections. For integrable adelic metrics on $L$ we put
\begin{equation}\label{eq:H-arithmetic-deficit}
 I_A(\bar L)=(d+1)L^d(\bar L\cdot\bar A^d)-\bar L^{d+1}
             =\widehat{\deg}\,D_A(\bar L).
\end{equation}
For relatively nef $\bar L$, Theorem~\ref{thm:ahi} gives $I_A(\bar L)\ge0$.
Write $h_{i,v}$ for the metric of $\bar L_i$ at $v$ and $\Psi_v$ for the
local energy \eqref{eq:H-energy} at that place.
We use unnormalized arithmetic degrees over a number field. Over $k(T)$
the factor in $-\log|a|_v=c_k[\kappa(v):k]\operatorname{ord}_v(a)$
is $c_k=\log\#k$ for finite $k$ and $c_k=1$ otherwise. With these
weights included in $\Psi_v$, the arithmetic degree of
$D_A(\bar L_1)\otimes D_A(\bar L_0)^\vee$ is the sum of
the local weights of its canonical section $1$. Thus
\begin{equation}\label{eq:H-height-energy-variation}
 I_A(\bar L_1)-I_A(\bar L_0)
    =\sum_v \Psi_v(h_{0,v},h_{1,v}).
\end{equation}
For $L^d=0$ we obtain
\begin{equation}\label{eq:height-energy-variation}
 \bar L_0^{d+1}-\bar L_1^{d+1}=\sum_v \Psi_v(h_{0,v},h_{1,v}).
\end{equation}
At a real place the integral is taken on the complexification with
conjugation-invariant metrics. The logarithmic ratio uses the extended
real absolute value and the Chern forms retain their geometric
normalization, so the energy is half its value over $\BC$.

Fix a place $v$ and an integrable adelic metric $\bar L_0$.
Write $\bar L_0^v(h_v)$ for the metric obtained by replacing $h_{0,v}$
by $h_v$ and keeping $h_{0,w}$ at every $w\ne v$.

\begin{thm}\label{thm:lb}
If $\bar L_0^v(h_v)$ is integrable and relatively nef adelically, then
\[
 \Psi_v(h_{0,v},h_v)\geq-I_A(\bar L_0).
\]
\end{thm}
\begin{proof}
The right side of \eqref{eq:H-height-energy-variation} has a single
term and $I_A(\bar L_0^v(h_v))\geq0$.
\end{proof}

When $L^d=0$ the bound is
$\Psi_v(h_{0,v},h_v)\geq\bar L_0^{d+1}$.

For a metric $h_v$ on $L_v$, strong relative nefness over $C_v$
means that some semipositive metrized ample bundle on $C_v$ makes
$h_v$ semipositive after pullback and tensor product. This condition
allows us to replace a metric at a single place.

\begin{lemma}\label{lem:one-place-globalization}
Suppose that $L$ is $f$-ample and that the adelic bundle
$\bar L_0$ is strongly relatively nef. For every place $v$ and
every strongly relatively nef metric $h_v$ on $L_v$, the adelic bundle
$\bar L_0^v(h_v)$ is strongly relatively nef.
\end{lemma}
\begin{proof}
We write $\bar L_0=\bar M_0-f^*\bar N_0$ as in
Definition~\ref{def:relative-nef-adelic}, and choose a semipositive
ample base twist $(N_v,n_v)$ for $h_v$. For sufficiently large
$m$, the bundles $L+mf^*N_0$ and $mN_{0,v}-N_v$ are ample.
We give $mN_{0,v}-N_v$ a semipositive Fubini--Study metric
$q_v$. We then equip $N=mN_0$ with $n_v+q_v$ at $v$ and
with $m\bar N_0$ at every other place. The bundles $\bar N$ and
\[
 \bar M=\bar L_0^v(h_v)+f^*\bar N
\]
have ample underlying bundles and semipositive metrics at every
place. At $v$, semipositivity follows from the choice of
$(N_v,n_v)$ and $q_v$; elsewhere we have
$\bar M=\bar M_0+(m-1)f^*\bar N_0$ locally.

It remains to make $\bar M$ and $\bar N$ nef adelically.
Their point heights are bounded below by comparison with fixed
Fubini--Study metrics. Over a number field, adding a common
sufficiently positive constant makes both bundles nef by the
successive-minima inequality and arithmetic positivity
\cite[\S A.5.4 and Theorem A.5.2]{YZ26}.
Over $k(T)$, choose semipositive model approximants uniformly
bounded relative to fixed nef metrics on the ample bundles $M,N$.
They have nonnegative degrees on vertical curves, and the uniform
height bounds control the degrees on horizontal curves divided by
their degrees over $T$. The pullback of one sufficiently positive
divisor on $T$ therefore makes all approximants nef. Passage to
the limit gives nef adelic factors. In both cases the common twist
cancels in $\bar M-f^*\bar N=\bar L_0^v(h_v)$, proving
strong relative nefness.
\end{proof}

If $f$ is smooth at a complex place $v$ and $h$ is smooth with strictly
positive vertical curvature, a sufficiently positive base twist gives
$h$ a strongly relatively nef presentation. Indeed a smooth positive
metric on an ample base bundle controls the horizontal and mixed
curvature once its tensor power is large, since the vertical block is
uniformly positive on the compact complex manifold. For $f$-ample $L$
we obtain a strongly relatively nef reference $\bar L_0$ by choosing
an ample $N$ on $C$ for which $L+f^*N$ is ample and taking the
difference of nef adelic metrics on the two ample bundles. Lemma~\ref{lem:one-place-globalization}
and Theorem~\ref{thm:lb} thus bound $\Psi_v(h_{0,v},h)$ for all smooth
vertically positive $h$. The cocycle identity allows any integrable
reference at $v$.

We keep the normalized nef polarization $\bar A$ fixed. At
archimedean places, we use the conjugation convention in the
\hyperref[par:notation-conventions]{notation conventions}.

\begin{cor}
\label{cor:H-local-bound}
Suppose that $L$ is nef and $f$-ample, and let $\bar L_0$ be a
strongly relatively nef adelic metric on $L$. At every place $v$,
one has $\Psi_v(h_{0,v},h)\ge-I_A(\bar L_0)$ for every
strongly relatively nef metric $h$ on $L_v$.
If $C_v$ is smooth and either $v$ is nonarchimedean or $X_v$ is smooth, then
\[
 \inf_{h\in\BH^+(L_v)}\Psi_v(h_{0,v},h)
 \ge-I_A(\bar L_0)>-\infty.
\]
\end{cor}
\begin{proof}
For strongly relatively nef $h$, Lemma~\ref{lem:one-place-globalization}
shows that $\bar L_0^v(h)$ is strongly relatively nef adelically,
so Theorem~\ref{thm:lb} proves the bound. Under the stated
smoothness hypotheses, we apply
Lemma~\ref{lem:relative-nef-approximation} to any
$h\in\BH^+(L_v)$ and choose $h_j\in\CH^+(L_v)$ with
$\Psi_v(h_{0,v},h_j)\to \Psi_v(h_{0,v},h)$. The bound is independent
of $j$, and passage to the limit proves the assertion for $h$.
\end{proof}

When $K$ is a number field, smooth changes at an archimedean place
with semipositive vertical curvature on a fixed fiberwise nef arithmetic model are also
covered directly by Theorem~\ref{thm:hi}, provided the other
archimedean metrics have semipositive vertical curvature.
At a real place the metric is invariant under complex conjugation;
at a complex place it is changed at both conjugate embeddings.

\Needspace{14\baselineskip}

\subsection{Critical metrics}\label{ss:critical-values}

We return to $(X/C,\bar A)$ and
$\CH^+(L)\subset\BH^+(L)$ to study the critical equation
and the relation between minima and semipositivity.
The calculation for smooth metrics suggests that, when $L^d=0$,
a metric in $\BH^+(L)$ minimizes the energy if and only if it
is globally semipositive. For surfaces, we study semipositivity
using second variations. Finally, we give examples from polystable
projective bundles and blow-ups of the projective plane at nine points.

\subsubsection*{Prescribed Monge--Amp\`ere measures}

Let $h$ be a continuous integrable metric on $L$, and put
$\alpha=c_1(L,h)$ and $\omega_A=c_1(\bar A)$. An integrable
variation is a continuous function $\varphi$ for which
$(\CO_X,e^\varphi|\cdot|)$ is integrable. Along
$|\cdot|_{h_t}=e^{t\varphi}|\cdot|_h$, put
$\alpha_t=\alpha-tdd^c\varphi$. Multilinearity of the Deligne
pairing gives
\begin{equation}\label{eq:H-energy-variations}
 \begin{aligned}
 \frac{d}{dt}\Psi_A(h_0,h_t)
   &=(d+1)\int_X\varphi
       \bigl(\alpha_t^d-L^d\omega_A^d\bigr),\\
 \frac{d^2}{dt^2}\Psi_A(h_0,h_t)
   &=-d(d+1)\int_X\varphi\,dd^c\varphi
                   \wedge\alpha_t^{d-1}.
 \end{aligned}
\end{equation}
These identities hold over $\BC$ and over any complete
nontrivially valued nonarchimedean field, with derivatives taken
from the polynomial in $t$ defined initially for rational $t$;
see \cite[Theorem~8.16]{BE21}. Over $\BC$, integration by
parts writes the second variation as
$d(d+1)\int_Xd\varphi\wedge d^c\varphi\wedge\alpha_t^{d-1}$
when the metrics and the variation are smooth.

We call $h$ critical if the first variation vanishes for every
smooth test function over $\BC$, or every model test function in
the nonarchimedean case. These functions are uniformly dense in
the continuous functions, by \cite[Theorem~5.18]{BE21} in the
nonarchimedean case. Thus, for every integrable $h$, in particular
for every $h\in\BH^+(L)$, criticality is equivalent to
\begin{equation}\label{eq:H-Monge-Ampere}
 c_1(L,h)^d=L^d\omega_A^d.
\end{equation}
The products are signed Monge--Amp\`ere measures, interpreted as
Chambert--Loir measures in the nonarchimedean case. When $L^d=0$,
the equation is $c_1(L,h)^d=0$. Here criticality tests all such
variations in the integrable metric space; a minimum on
$\BH^+(L)$ is stationary only along variations admissible for
both signs of $t$.

Suppose that $h$ is a smooth local minimum at which every smooth
variation is admissible for both signs of the parameter. Then
$h$ satisfies \eqref{eq:H-Monge-Ampere}, and its second
variation is nonnegative in every direction. The localization
argument in the proof of Proposition~\ref{prop:surface-hessian}
gives $\alpha^{d-1}\ge0$. For a surface, this is exactly
$\alpha\ge0$. If $L^d>0$ and $\omega_A>0$, we obtain
$\alpha>0$ in every dimension: equation~\eqref{eq:H-Monge-Ampere}
makes the determinant positive, and nonnegative cofactors then
force every eigenvalue to be positive. For $L^d=0$ in higher
dimensions, this argument requires further information on $\alpha$.

For the next proposition, we assume that $f$ is smooth over
$\BC$ and that the fixed polarization $\bar A$ has a smooth
semipositive metric.

\begin{proposition}\label{prop:minimizer}
Assume that $L$ is nef on $X$. Every critical metric
$h_*\in\CH^+(L)$ is globally semipositive. Its curvature has
rank $d-1$ everywhere if $L^d=0$, and is positive if $L^d>0$
and $\omega_A>0$.
\end{proposition}
\begin{proof}
Put $\alpha_*=c_1(L,h_*)$. Strict vertical positivity makes
all sufficiently small smooth variations admissible for both signs
of the parameter. The first-variation formula identifies criticality
with \eqref{eq:H-Monge-Ampere}.

Since $L$ is nef, $L^d\ge0$, and the critical equation gives
$\alpha_*^d\ge0$. We separate the vertical and base directions
and write the curvature matrix as
\[
 \begin{pmatrix}V&b\\b^*&c\end{pmatrix},\qquad V>0.
\]
Its determinant is $\det V(c-b^*V^{-1}b)\ge0$, so its Schur
complement is nonnegative and $\alpha_*$ is semipositive.
The determinant vanishes when $L^d=0$, giving rank $d-1$,
and is positive when $L^d>0$ and $\omega_A>0$, giving
positive curvature. For $d=1$ these conclusions follow directly
from \eqref{eq:H-Monge-Ampere}.
\end{proof}
The proposition applies to local minima in $\CH^+(L)$.
Since nefness is used only to ensure $L^d\ge0$, it may be
omitted when $L^d=0$. On surfaces,
Proposition~\ref{prop:surface-hessian} instead tests semipositivity
by second variations, provided those variations are admissible.

If $L$ is ample over $\BC$ and $\omega_A$ is smooth and
positive, then \eqref{eq:H-Monge-Ampere} is the Calabi--Yau
equation. By Yau's theorem \cite{Yau78}, it has a smooth positive
solution, unique up to a constant. After fixing that constant,
we obtain a canonical metric relative to the probability volume
form $\omega_A^d$. The metrics of positive curvature form a
convex space in logarithmic weights. Since the second variation
in \eqref{eq:H-energy-variations} is positive in every nonconstant
direction, the canonical metric is the unique minimum on this
space modulo constants.

The variational approach to prescribed Monge--Amp\`ere measures
for globally semipositive metrics is developed in
\cite{BBGZ13,BEGZ10} over $\BC$, and in \cite{BFJ15,FGK25}
in the nonarchimedean setting. For continuous semipositive metrics
on an ample bundle, uniqueness up to a constant follows from
\cite[Corollary~2.2]{YZ17}. Boucksom--Jonsson
\cite[Theorem~B]{BJ25} prove independence of finite-energy
measures from the ample polarization. In contrast,
Theorem~\ref{thm:main-energy} bounds the energy on the larger
space $\BH^+(L)$ for almost every realization, including when
$L^d=0$, by the height inequality.

\begin{remark}
For a line bundle $N$ on $C$, the restriction of $L$ to each
fiber is unchanged by $L_N=L+f^*N$, while
\[
 L_N^d=L^d+d\deg(N)L_\eta^{d-1}.
\]
For a smooth metric on $N$, with curvature $\beta$, the
curvature of the tensor product metric is
$\alpha_N=\alpha+f^*\beta$, and
\[
 \alpha_N^d=\alpha^d+d f^*\beta\wedge\alpha^{d-1}.
\]
When $\omega_A$ is smooth and positive and $L_N$ is ample,
Yau's theorem gives its canonical metric relative to $\omega_A^d$.
If $L$ is nef and $f$-ample and
$N$ is ample, this applies to $L+\epsilon f^*N$ for every
rational $\epsilon>0$. The behavior of these canonical metrics as $\epsilon\to0$
requires a separate convergence argument.
\end{remark}

\subsubsection*{Zero-volume minima}

We assume $L^d=0$ for the remainder of this subsection, except
where a different assumption is stated. The correction in
\eqref{eq:H-energy} then vanishes, and the critical equation is
independent of $\bar A$.

For the next conjecture, we use the local setting of
Section~\ref{ss:local-metrics}, with $X$ and $C$ geometrically
integral and $L$ relatively ample. We fix $h_0\in\BH^+(L)$.

\begin{conj}\label{conj:sp}
\label{conj:nef-minimizers}
For every $h_*\in\BH^+(L)$, one has
\[
 \Psi_A(h_0,h_*)=\inf_{h\in\BH^+(L)}\Psi_A(h_0,h)
 \quad\Longleftrightarrow\quad
 (L,h_*)\text{ is globally semipositive}.
\]
\end{conj}

This characterization is independent of the integrable reference
metric by the cocycle identity; it does not assert existence of
a minimizer.

For smooth $f$ over $\BC$, criticality and global semipositivity
are equivalent on $\CH^+(L)$. One direction is
Proposition~\ref{prop:minimizer}; conversely, a semipositive
curvature has top power of mass $L^d=0$ and hence satisfies the
critical equation.

\begin{remark}\label{rem:integrable-critical}
The following complex and model examples show that the conclusion
of Proposition~\ref{prop:minimizer} need not hold for integrable
metrics. In both, $f$ is a smooth surface fibration, $L$ is nef
and relatively ample with $L^2=0$, and $h\in\BH^+(L)$ satisfies
\[
 c_1(L,h)^2=0,\qquad (L,h)\text{ is not globally semipositive}.
\]
Thus every first variation vanishes for smooth, respectively model,
test functions.

\smallskip\noindent\emph{A complex example.}
Over $\BC$, we take $X=\BP^1\times C$, $f=p_2$, and
$L=p_1^*\mathcal O(1)$. We choose a nonconstant smooth real
function $\psi$ on $C$ and give $L$ the metric whose homogeneous
weight in coordinates $[z_0:z_1]$ is
\[
 u=\max\{\log|z_0|_{\BC}+\psi,\ \log|z_1|_{\BC}-\psi\}.
\]
Here the normalized metric value of a homogeneous linear section
$s$ is $|s|_{\BC}e^{-u}$.
We choose a positive base curvature $\omega$ large enough that
$\omega\pm dd^c\psi\ge0$. Twisting by the corresponding base
metric makes both branches, and hence their maximum, semipositive.
Thus $h$ is strongly relatively nef and integrable. On the open
regions where the respective branches dominate, its curvature is
$f^*dd^c\psi$ and $-f^*dd^c\psi$. Since $\psi$ is
nonconstant, $h$ is not globally semipositive.

To compute the signed Monge--Amp\`ere measure, we work in the
chart $z_0=1$ and put $x=\log|z_1/z_0|_{\BC}$. We regularize
the local weight by $u_\epsilon=\psi+\chi_\epsilon(x-2\psi)$,
where $\chi_\epsilon$ is a smooth convex approximation to
$\max(0,\cdot)$ agreeing with it outside $[-\epsilon,\epsilon]$.
For $q=\chi_\epsilon'(x-2\psi)$, we compute
\[
 (dd^cu_\epsilon)^2
 =2(1-2q)\chi_\epsilon''(x-2\psi)\,
       f^*(dd^c\psi)\wedge dx\wedge d^cx.
\]
The transverse coefficient has bounded total variation and integral
$2\int_0^1(1-2q)\,dq=0$, and its support shrinks to $x=2\psi$.
The signed measures therefore converge to zero. Since the fixed
base twist makes all the regularized metrics semipositive, we
may apply continuity of mixed products and conclude that
$c_1(L,h)^2=0$.

\smallskip\noindent\emph{A model example.}
Let $R$ be a complete discrete valuation ring. We blow up a
rational point of the special fiber of $\BP^1_R$, and denote
the resulting model by $\mathcal C$ and the exceptional curve
by $D$. Thus $D^2=-1$. We put
\[
 \mathcal E=\mathcal O(D)\oplus\mathcal O(-D),\qquad
 \mathcal X=\BP_{\mathcal C}(\mathcal E),\qquad
 \mathcal L=\mathcal O_{\mathcal X}(1),
\]
using the quotient convention for the projective bundle.
The generic fiber is $\BP^1\times\BP^1$, with
$L=p_1^*\mathcal O(1)$. The model metric $h$ is strongly relatively nef,
since $\mathcal L$ is relatively ample over $\mathcal C$.
For each component $C_i$ of the reduced special fiber of
$\mathcal C$,
\[
 \deg_{\BP(\mathcal E|_{C_i})}\mathcal L^2
       =\deg_{C_i}\det\mathcal E=0.
\]
The formula for the measure of a model metric \cite[\S1.3]{CL10} therefore gives
$c_1(L,h)^2=0$. Nevertheless, the quotient
$\mathcal E|_D\twoheadrightarrow\mathcal O(D)|_D$ defines a
section of degree $-1$ for $\mathcal L$, so $h$ is not globally
semipositive. Proposition~\ref{prop:surface-model-minimum}
constructs an admissible variation decreasing its energy.
\end{remark}

For the next proposition, let $f$ be flat and projective over
$\BC$ or a complete discretely valued field, with possibly
singular fibers. We continue to assume $L^d=0$.

\begin{proposition}
\label{prop:constant-energy}
The continuous globally semipositive metrics on $L$ form a
convex set on which $\Psi_A(h_0,-)$ is constant for every fixed
integrable reference metric $h_0$.
\end{proposition}
\begin{proof}
All mixed top-degree Monge--Amp\`ere measures of two such metrics
are nonnegative and have total mass $L^d=0$, so they vanish.
At a nonarchimedean place these are the corresponding
Chambert--Loir measures; see \cite[\S1.3]{CL10}.
Hence $\Psi_A(h_1,h_2)=0$, and the cocycle identity gives
$\Psi_A(h_0,h_1)=\Psi_A(h_0,h_2)$.
The common value depends on the reference metric and is zero
when that reference is globally semipositive.
\end{proof}

Combining Propositions~\ref{prop:minimizer} and
\ref{prop:constant-energy}, we see that all smooth critical
metrics with strictly positive vertical curvature have the same
energy for a fixed reference $h_0$.

\subsubsection*{Second variation on surfaces}

For a surface, we can recover semipositivity of the curvature
from the second variation. We allow arbitrary $L^2$ in the
following expansion and in Proposition~\ref{prop:surface-hessian}.
Let $h$ be a continuous integrable metric and put
$\alpha=c_1(L,h)$. For a smooth real function $\varphi$ over
$\BC$, or a model function over a complete nontrivially valued
nonarchimedean field, we set $|\cdot|_{h_t}=e^{t\varphi}|\cdot|_h$.
Expanding \eqref{eq:H-energy} by multilinearity, we obtain
\begin{equation}\label{eq:surface-energy-expansion}
 \begin{aligned}
 \Psi_A(h,h_t)&=3t\int_X\varphi
           \bigl(\alpha^2-L^2c_1(\bar A)^2\bigr)\\
       &\quad+3t^2Q_h(\varphi)
       +t^3\int_X\varphi\,(dd^c\varphi)^2,\\
 Q_h(\varphi)&=-\int_X\varphi\,dd^c\varphi\wedge\alpha.
 \end{aligned}
\end{equation}
For model variations, the products are local intersection measures
and $t$ is rational. The second derivative at zero is
$6Q_h(\varphi)$.
For $L^2=0$, the proposition below characterizes semipositivity
by criticality and nonnegative second variations over $\BC$ and
for rational model metrics. We propose the following extension.

\begin{conj}\label{conj:surface-second-variation}
Let $f:X\to C$ be a flat projective fibration from a smooth
geometrically integral surface to a smooth geometrically integral
curve over a complete nontrivially valued nonarchimedean field. Let $L$ be
$f$-ample with $L^2=0$. A metric $h\in\BH^+(L)$ is globally
semipositive if and only if it is critical and
\[
 Q_h(\varphi)\ge0\qquad\text{for every model function }\varphi.
\]
\end{conj}

The forward implication follows from the mass formula and the
local Hodge index theorem \cite[Theorem~2.1]{YZ17}. The
conjecture concerns the converse for continuous metrics. The tests
range over all model functions, independently of their admissibility
in $\BH^+(L)$.

Let $X$ be a smooth geometrically integral projective surface.
In the next proposition, we use one of the following two classes
of metrics and test functions:
\begin{itemize}
\item over $\BC$, $h$ is a continuous integrable metric on $L$,
and the tests are all smooth real functions;
\item over a complete discretely valued field, $h$ is a rational
model metric, and the tests are all model functions on all
dominating models.
\end{itemize}

\begin{proposition}
\label{prop:surface-hessian}
In either the complex or the model case,
\[
 (L,h)\text{ is semipositive}
 \quad\Longleftrightarrow\quad
 Q_h(\varphi)\ge0\text{ for every test function }\varphi.
\]
If $L^2=0$, either condition also implies $c_1(L,h)^2=0$.
\end{proposition}
\begin{proof}
\smallskip\noindent\emph{The complex case.}
Over $\BC$, integration by parts gives
\[
 Q_h(\varphi)=\int_Xd\varphi\wedge d^c\varphi\wedge\alpha.
\]
The integral is nonnegative when $\alpha$ is positive. For the
converse, we localize the inequality using a real cutoff $\eta$
supported in a coordinate chart and a real linear function $\ell$
in these coordinates. We extend
$\varphi_n=n^{-1}\eta\cos(n\ell)$ and
$\psi_n=n^{-1}\eta\sin(n\ell)$ by zero. Adding their quadratic
forms, we find
\[
 Q_h(\varphi_n)+Q_h(\psi_n)
 =\int_X\eta^2d\ell\wedge d^c\ell\wedge\alpha
  +\frac1{n^2}\int_Xd\eta\wedge d^c\eta\wedge\alpha.
\]
We let $n\to\infty$ and conclude that $\alpha$ pairs
nonnegatively with every square cutoff times a constant positive
rank-one $(1,1)$-form. Integrability expresses $\alpha$ as a
difference of positive currents, so its coefficients are measures.
The local matrix criterion for positivity on a surface now gives
$\alpha\ge0$.

\smallskip\noindent\emph{The model case.}
Over a complete discretely valued field $F$, the forward implication
is the local Hodge index theorem \cite[Theorem 2.1]{YZ17}.
For the converse, we represent $h$ by a rational line bundle
$\mathcal L$ on a normal projective model $\mathcal X$ over
the valuation ring $R$. To test the degree of $\mathcal L$ on
an integral curve $Z$ in the special fiber, we take the normalized
blowup $\pi:\mathcal X'\to\mathcal X$ of $Z$ and denote its
exceptional Cartier divisor by $E$. We have
\begin{equation}\label{eq:surface-blowup-test}
 \pi_*(E^2)=-m_Z[Z],\qquad
 Q_h(g_E)=-c\,E^2\cdot\pi^*\mathcal L
          =c\,m_Z\deg_Z\mathcal L,
\end{equation}
where $g_E=-\log|1|_{\mathcal O(E)}$,
$c=-\log|\varpi|>0$ for a uniformizer $\varpi$, and $m_Z>0$
is the Hilbert--Samuel multiplicity of the two-dimensional local
ring $\mathcal O_{\mathcal X,\eta_Z}$.
For the first identity, before normalization the exceptional
fiber over $\eta_Z$ is the projectivized associated graded ring.
It has degree $m_Z$ with respect to $\mathcal O(1)$, while
$\mathcal O(E)$ restricts to $\mathcal O(-1)$. The projection
formula preserves the identity after normalization. In particular,
$m_Z=1$ if $\mathcal X$ is regular at $\eta_Z$. The second
identity follows from the normalization of local intersections in
\cite[Appendix 1, ``Local intersections'']{YZ17}. Consequently,
$Q_h\ge0$ implies $\deg_Z\mathcal L\ge0$ for every such
curve $Z$. This is semipositivity of the model metric.

Finally, when $L^2=0$, the Monge--Amp\`ere measure of a
semipositive $h$ is positive of mass zero, and hence vanishes.
\end{proof}

To test a minimum on $\BH^+(L)$, the variation must remain
in that space. For a relatively semiample surface model, blowing
up a curve of negative degree gives such a variation. We return
to $L^d=0$, allowing singular fibers of $f$.

Let $F$ be a complete discretely valued field and let $f:X\to C$
be as in Section~\ref{ss:local-metrics}, with $\dim X=2$.
Let $h$ be a rational model metric represented by a rational
line bundle $\mathcal L$ on a projective model
$\mathcal f:\mathcal X\to\mathcal C$ of $f$ over the valuation
ring.

\begin{proposition}
\label{prop:surface-model-minimum}
Assume that $L$ is $f$-ample and that $\mathcal L$ is
$\mathcal f$-semiample. If $h$ is a local minimum of $\Psi_A(h_0,-)$
on $\BH^+(L)$ for the uniform topology on logarithmic metrics,
then $(L,h)$ is globally semipositive.
\end{proposition}
\begin{proof}
We first replace the relatively semiample model by a relatively
ample one. The relative contraction defined by a sufficiently
divisible power of $\mathcal L$ is an isomorphism on the generic
fiber, because $L$ is $f$-ample. On the contracted model, the
same metric is represented by a relatively ample rational line
bundle. After normalization, we may assume that $\mathcal X$
is normal and $\mathcal L$ is $\mathcal f$-ample.

Suppose that $h$ is not semipositive. We choose an integral
curve $Z$ in the special fiber with $\deg_Z\mathcal L<0$.
Let $\pi$ and $E$ be the normalized blowup and exceptional
divisor used in \eqref{eq:surface-blowup-test}, and put
$\mathcal M=\pi^*\mathcal L$. For every sufficiently small positive
rational $\epsilon$, the bundle
\[
 \mathcal L_\epsilon=\mathcal M-\epsilon E
\]
is ample relative to $\mathcal C$. Indeed, $\mathcal O(-E)$
is $\pi$-ample, and after taking a sufficiently large power of
$\mathcal L$ we may use composition of relative ampleness
\cite[Tags 01OF and 0C4K]{Stacks}. Adding the pullback of a
sufficiently ample bundle on $\mathcal C$ makes
$\mathcal L_\epsilon$ ample over the valuation ring. Thus its
metric $h_\epsilon$ lies in $\BH^+(L)$ and converges uniformly
to $h$.

We compute the change of energy along this admissible variation.
By the projection formula and \eqref{eq:surface-blowup-test},
\[
 \mathcal M^2\cdot E=0,\qquad
 \mathcal M\cdot E^2=-m_Z\deg_Z\mathcal L.
\]
The change of the inverse Deligne metric is therefore
\begin{equation}\label{eq:surface-model-descent}
 \begin{split}
 \Psi_A(h,h_\epsilon)
   &=3c\epsilon \mathcal M^2\cdot E
      -3c\epsilon^2 \mathcal M\cdot E^2+c\epsilon^3E^3\\
   &=3c\,m_Z\deg_Z\mathcal L\,\epsilon^2+cE^3\epsilon^3<0
 \end{split}
\end{equation}
for sufficiently small $\epsilon>0$, where $c=-\log|\varpi|$.
By the cocycle identity, this gives
$\Psi_A(h_0,h_\epsilon)<\Psi_A(h_0,h)$ for sufficiently small
$\epsilon>0$, contrary to local minimality.
\end{proof}

Conversely, a projective flat model bundle nef on every component
of the geometric special fiber induces a globally semipositive
metric \cite[\S A.5.1]{YZ26}. For $f$-ample $L$ with
$L^d=0$, Conjecture~\ref{conj:nef-minimizers} predicts that
it minimizes the energy.

\subsubsection*{Projective bundles and blow-ups}

\paragraph*{Projective bundles and Hermitian--Einstein metrics.}
For a polystable degree-zero bundle $E$ on a smooth complex
projective curve, Narasimhan--Seshadri and Donaldson
\cite{NS65,Do83} give a unitary-flat Hermitian metric. Its
quotient metric $h_*$ on $L=\CO_{\BP(E)}(1)$ is globally
semipositive, vertically positive, and satisfies
$c_1(L,h_*)^{\rk E}=0$. By \cite[Theorem 0.1]{FLW} and
\eqref{eq:FLW-normalization},
\[
 \Psi_A(h_*,h)\ge0
\]
for every smooth vertically positive metric $h$.
Lemma~\ref{lem:relative-nef-approximation} extends this to
$\BH^+(L)$, so $h_*$ attains the global minimum. Within
Hermitian metrics with fixed unitary-flat determinant, smooth
attainment requires polystability; see
\eqref{eq:Donaldson-normalization} and \cite{Do85,UY86}.

\paragraph*{Semipositive metrics on the blow-up of $\BP^2$.}
For the surface $Y'$ of Proposition~\ref{ex:nine-points},
Chen--Filip--Sun--Tosatti--Zhang \cite[Corollary 1.2]{CFSTZ26}
prove that $L'$ admits a smooth globally semipositive metric,
including in the non-semiample cases. Since $(L')^2=0$, any such
metric satisfies $c_1(L',h)^2=0$. We pull it back to $L=p^*L'$
on the pencil model of Proposition~\ref{prop:polarized-surface-blowups}.
The pullback is globally semipositive, and hence strongly relatively
nef. All these metrics have the same energy by
Proposition~\ref{prop:constant-energy}. Their vertical curvature
on the pencil may be degenerate.

\paragraph*{Nonarchimedean metrics on the blow-up of $\BP^2$.}
Let $R=\CO_F$ for a complete discretely valued field $F$.
As in Example~\ref{ex:moving-nine-points}, choose a smooth plane
cubic $\mathcal B\subset\BP^2_R$ and nine sections
$s_1,\ldots,s_9\in\mathcal B(R)$ with distinct generic points
$P_i$. Successively blowing up their lifts on the strict transform
of $\mathcal B$ gives a smooth projective model
\[
 \pi:\mathcal Y'\longrightarrow\BP^2_R,\qquad
 J_0=\pi^*\CO_{\BP^2_R}(1),\qquad
 \mathcal L'=3J_0-\sum_{i=1}^9E_i
           =\CO_{\mathcal Y'}(\widetilde{\mathcal B}),
\]
where $E_i$ denotes the total transform of the $i$-th exceptional
divisor. Collisions of sections in reduction are allowed: each
center is Cartier on the strict transform of $\mathcal B$, which
remains isomorphic to $\mathcal B$ and carries the remaining lifts.
On each geometric fiber, $\widetilde{\mathcal B}$ is smooth
integral of self-intersection zero, hence nef. Thus $\mathcal L'$
induces a globally semipositive metric on its generic fiber $L'$
\cite[\S A.5.1]{YZ26}, with
\[
 c_1(L',h_{\mathcal L'})^2=0.
\]
This works in every residue characteristic. If
$\CO_B(3)(-P_1-\cdots-P_9)\in\Pic^0(B)$ is non-torsion,
where $B=\mathcal B_F$, then $L'$ is not semiample by
Proposition~\ref{ex:nine-points}. Pullback to the pencil model
of Proposition~\ref{prop:polarized-surface-blowups} gives
semipositive metrics on $L=p^*L'$. They belong to
$\BH^+(L)$ and have the common energy of the globally semipositive
metrics in Proposition~\ref{prop:constant-energy}.

\section{Families of fibrations}\label{sec:families}
\label{ss:family-energy}

Let $X\xrightarrow{f}C\xrightarrow{\pi}S$ be a family of polarized
fibrations, and let $L$ be a rational line bundle on $X$. We use
Deligne pairings to construct a canonical rational line bundle $D_A(L)$ on $S$.
We fix the polarization $\bar A$, normalized by $A_s^d=1$, and vary
the metric on $L$. The logarithmic change of the induced norm on
$D_A(L)$ at $s$ is the energy on $X_s$.
After studying the heights and envelopes of $D_A(L)$, we construct
projective models in Section~\ref{ss:adelic-models-approximations}
and use them to prove Theorem~\ref{thm:main-energy} in
Sections~\ref{ss:complex-points} and~\ref{ss:nonarchimedean}.

\subsection{Adelic heights and envelopes}\label{ss:family}
\label{ss:local-envelopes}

Over a number field $K$, Theorem~\ref{thm:ahi} gives
$h_{D_A(\bar L)}(s)\ge 0$ for $s\in S(\bar K)$ whenever
$\bar L$ is relatively nef and $L_s$ is nef.
When $L$ is $f$-ample and every $L_s$ is nef, we conjecture that the
local envelopes define a nef adelic metric on $D_A(L)$ whose height
at each algebraic point is the infimum of the heights obtained from
the Hermitian model metrics specified below.
For a smooth complex family, we prove semipositivity of the inverse
Deligne pairing when $\alpha=c_1(L,h)$ is strictly positive on
$T_{X/C}$ and $(\alpha|_{X_s})^d=0$ on every fiber.

Let
\[
 X\xrightarrow{f}C\xrightarrow{\pi}S,\qquad g=\pi\circ f,
\]
be proper flat morphisms of integral quasi-projective
varieties over a number field $K$, with $C$ and $S$ normal. Assume
that $f$ and $\pi$ have pure relative dimensions $d-1$ and $1$,
respectively, and that $X_s$ and $C_s$ are geometrically integral
for every geometric point $s$ of $S$. The maps are projective,
and $g$ is flat of relative dimension $d$. The fibers may be
singular; these hypotheses suffice for
\cite[Theorem 4.1.3]{YZ26}.
We fix a normalized polarization $\bar A$ obtained from a nef
adelic $\mathbb Q$-line bundle $\bar A_0$ induced by a Hermitian
model, with ample underlying bundle. Thus each parameter fiber is
a polarized fibration $(X_s/C_s,\bar A_s)$ with $A_s^d=1$.
Let $L$ be a line bundle on $X$ that is fiberwise nef for $f$.
By flatness, $L_s^d$ is independent of $s$. Nefness on the whole
of $X_s$ will be imposed in the positivity statements.

For an integrable adelic metric $\bar L$ relatively nef for $f$,
the induced metric on $D_A(L)$ gives
\begin{equation}\label{eq:H-family-pairing}
 D_A(\bar L)
 =-D_g(\bar L)+(d+1)L_s^d\langle\bar L,\bar A^{[d]}\rangle_g.
\end{equation}
The mixed pairing contains one copy of $\bar L$ and $d$ of
$\bar A$, which stays fixed when the metric on $L$ varies.
For $L_s^d=0$, the pairing is $-D_g(\bar L)$.

For a finite extension $F/K$ and $s\in S(F)$,
base change and the absolute height normalization give
\begin{equation}\label{eq:family-height-normalization}
 h_{D_A(\bar L)}(s)
 =\frac{I_{A_s}(\bar L_s)}{[F:\mathbb Q]}
 =\frac{1}{[F:\mathbb Q]}
   \bigl((d+1)L_s^d(\bar L_s\cdot\bar A_s^d)
                     -\bar L_s^{d+1}\bigr).
\end{equation}
Extending $F$ multiplies both the numerator and the denominator by the
degree of the extension, so the height is independent of $F$.

To apply Theorem~\ref{thm:ahi} on $X_s/C_s$, we check that restriction
to $X_s$ preserves relative nefness.
We use adelic line bundles and the notions of nefness and
integrability of \cite[\S2.5]{YZ26}, defined by limits of models
in the boundary topology. On the quasi-projective total family,
we define strong relative nefness by adding a pullback from the base, as in
Definition~\ref{def:relative-nef-adelic}:
\[
 \bar L=\bar M-f^*\bar N,
 \qquad \bar M\text{ and }\bar N\text{ nef adelic},
 \qquad N\text{ ample on }C.
\]
Restriction to $X_s$, over any finite extension defining $s$,
preserves strong relative nefness: both factors remain nef and
$N|_{C_s}$ is ample.

An integrable adelic metric on $L$ is relatively nef if there is
one fixed nef adelic bundle $\bar A'$ on $X$ such that
$\bar L+\epsilon\bar A'$ is strongly relatively nef for every
rational $\epsilon>0$. Restriction to $X_s$ preserves this
condition, since $\bar A'_s$ is nef.
Thus, when $L_s$ is nef, we may apply Theorem~\ref{thm:ahi} to
$\bar L_s$ on $X_s/C_s$. By \eqref{eq:family-height-normalization},
we obtain $h_{D_A(\bar L)}(s)\ge 0$.

To take the infimum of these heights, we let $\CH^+(L)$ be the set of
strongly relatively nef adelic metrics induced by Hermitian models, with smooth
strictly positive vertical curvature at the archimedean places
and model metrics at the finite places. On singular complex
spaces, we interpret smoothness and vertical positivity in the
ambient sense. On smooth parameter fibers, the local restrictions
of these metrics belong to $\BH^+(L_s)$ of
Section~\ref{ss:local-metrics}. At finite places we allow all
strongly relatively nef model metrics in this family class,
including those without a relatively ample representing model.

When $\CH^+(L)$ is nonempty, we put
\[
 \widehat h_{D_A(L)}(s)=
       \inf_{\bar L\in\CH^+(L)}h_{D_A(\bar L)}(s).
\]
If $L_s$ is nef for every $s\in S(\bar K)$, the infimum is finite
and nonnegative. Any one metric in $\CH^+(L)$ gives a finite upper
bound. If $S$ is proper, it follows that $D_A(L)$ is nef.

The class $\CH^+(L)$ is nonempty when $L$ is $f$-ample.
Indeed, we choose an ample $N$ on $C$ such that
$M=L+f^*N$ is ample \cite[Tag 0892]{Stacks}. Fubini--Study
metrics from projective immersions defined by powers of $M,N$
give nef adelic bundles $\bar M,\bar N$
\cite[\S\S2.5.5--2.5.6]{YZ26}, with
\[
 \bar L=\bar M-f^*\bar N\in\CH^+(L).
\]
On smooth parameter fibers, the same construction gives
$\CH^+(L_s)\ne\varnothing$ as in Section~\ref{ss:local-metrics}.

\paragraph*{The family envelope.}
For the local formulation, let $F=\mathbb C$ or a complete
discretely valued field, and consider a tower
$X\xrightarrow{f}C\xrightarrow{\pi}S$ over $F$ with the same
geometric hypotheses as above. Fix a semipositive metrized ample
rational polarization $\bar A_0$ and normalize it to $\bar A$
with $A_s^d=1$.

For such a local family, $\CH^+(L)$ denotes the integrable
family test class consisting, over $\mathbb C$, of smooth metrics
with strictly positive curvature on the relative tangent spaces
of $f$, interpreted in the ambient sense on singular spaces.
Over a nonarchimedean field, it consists of strongly relatively
nef rational model metrics: some ample bundle $N$ on $C$ has
a continuous semipositive metric $n$ such that
$(L,h)+f^*(N,n)$ is semipositive on $X$.
On quasi-projective analytic spaces, semipositivity is understood
locally: we use continuous plurisubharmonic weights over
$\mathbb C$ and locally uniform limits of nef model metrics
in the nonarchimedean case. When the complex family and $L$
carry compatible real structures and $\bar A$ is invariant,
we also use the class of conjugation-invariant metrics. The
local assertion below applies to this class as well.

The Fubini--Study construction makes the local class nonempty
when $L$ is $f$-ample. Its restrictions to smooth parameter
fibers lie in $\BH^+(L_s)$ of Section~\ref{ss:local-metrics}.
In the nonarchimedean case, the family class allows strongly
relatively nef models, whereas $\CH^+(L_s)$ requires a relatively
ample representing model. Over a number field, $\CH^+(L)$
continues to denote the adelic Hermitian model class.

For a local family, write $|\cdot|_{h,s}$ for the norm of $D_A(h)$ on
the fiber $D_A(L)_s$. We take the pointwise envelope of these norms by setting
\begin{equation}\label{eq:local-family-envelope}
 |u|_{{\rm env},s}
   =\sup_{h\in\CH^+(L)}|u|_{h,s},
 \qquad s\in S^{\rm an},\quad u\in D_A(L)_s.
\end{equation}
Here $S^{\rm an}$ is the complex or Berkovich analytic space;
at a Berkovich point $s$, the fiber is over its completed residue
field $\kappa(s)$. For the rational line $D_A(L)$, we clear
denominators and take the corresponding root norms. The polarization stays
fixed; when $L_s^d=0$, $D_A(L,h)$ is the inverse Deligne line
with its induced metric, independent of $\bar A$.
For a family over a number field, we apply the local construction
componentwise after completion and use conjugation-invariant
metrics at real places. The global height infimum is taken over
the adelic Hermitian model class $\CH^+(L)$.

\begin{conj}\label{conj:envelope}
Assume that $L$ is $f$-ample and that $L_s$ is nef on every
geometric parameter fiber $X_s$.
\begin{enumerate}[label=(\roman*),leftmargin=*,itemsep=0.3\baselineskip]
\item \emph{Local regularity.} For a local family as above,
the envelope \eqref{eq:local-family-envelope} is finite and
defines a continuous semipositive metric on $D_A(L)$.
\item \emph{Adelic compatibility.} For a family over a number field
as above, the local envelope norms at all places define a nef
adelic line bundle
$D_A^{\rm env}(L)=(D_A(L),|\cdot|_{\rm env})\in\widehat\Pic(S)$
in the sense of \cite{YZ26}, and
\[
 h_{D_A^{\rm env}(L)}(s)=\widehat h_{D_A(L)}(s)
 =\inf_{\bar L\in\CH^+(L)}h_{D_A(\bar L)}(s)
 \qquad(s\in S(\bar K)).
\]
\end{enumerate}
\end{conj}

To relate the construction to the energy, fix $h_0\in\CH^+(L)$
in a local family and set
\[
 E_h(s)=\Psi_A(h_0,h)(s)=\log\frac{|\cdot|_{h_0,s}}{|\cdot|_{h,s}},
 \qquad
 e(s)=\inf_{h\in\CH^+(L)}E_h(s).
\]
At classical points, $E_h(s)=\Psi_s(h_0|_{X_s},h|_{X_s})$ by
base change for Deligne pairings. Where $e$ is finite,
\[
 |\cdot|_{{\rm env},s}=\exp(-e(s))|\cdot|_{h_0,s}.
\]
By the cocycle identity, replacing $h_0$ by $h_1$ adds
$\Psi_A(h_1,h_0)(s)$ to $e(s)$ and leaves the envelope norm unchanged.
For $S=\operatorname{Spec}F$ with $X,C$ smooth, local finiteness is
the assertion of Conjecture~\ref{conj:H-energy-boundedness}, by
the equality of infima in Lemma~\ref{lem:relative-nef-approximation}.
On a point, continuity and semipositivity are automatic.

\subsubsection*{An attained smooth complex family of zero volume}

We return to a complex family $(X/C,\bar A)$ with fixed
normalized semipositive polarization. Let
$X\xrightarrow{f}C\xrightarrow{\pi}S$ be smooth projective
maps of relative dimensions $d-1$ and $1$, and put $g=\pi\circ f$.
Suppose $\alpha=c_1(L,h)$ is smooth and strictly positive on
$T_{X/C}$. We study the case of zero volume, in which the
corrected pairing is independent of the polarization.

\begin{proposition}\label{prop:attained-family}
If $(\alpha|_{X_s})^d=0$ as a differential form for every $s\in S$,
then the inverse Deligne pairing
$D_A(\bar L)=-D_g(\bar L)$ in the rational Picard group has
semipositive curvature:
\[
 c_1(D_A(\bar L))=-g_*\alpha^{d+1}\ge0.
\]
\end{proposition}
\begin{proof}
We test the curvature on a disk in $S$. Positivity on $T_{X/C}$
allows us to take $\alpha$-orthogonal lifts of the directions
along $C/S$ and the disk. The identity $(\alpha|_{X_s})^d=0$
makes the diagonal entry in the $C/S$ direction zero. After
normalizing the vertical block, the Hermitian matrix is
\[
 \begin{pmatrix}I_{d-1}&0&0\\0&0&b\\0&\bar b&c\end{pmatrix}.
\]
Its determinant is $-b\bar b\le0$. The factors coming from the
exterior power and the volume form are positive, so integration over the fiber gives
$-g_*\alpha^{d+1}\ge0$.
\end{proof}

The proposition applies when smooth interior minimizers on the
fibers form a smooth family, since their top curvature forms vanish
pointwise.

\begin{remark}
For a family of projective bundles $X=\mathbb P(E)$ with
$L=\mathcal O(1)$ and $\deg E_s=0$, we can compare
$D_A(\bar L)=-D_g(\bar L)$ with Donaldson's metrized determinant
bundle \cite{Do87} through their metric variations. For quotient
metrics induced by Hermitian norms, the logarithmic variation
of the inverse Deligne pairing is the functional together with
the determinant contribution in \eqref{eq:Donaldson-normalization}.
When the determinant metric is fixed, this variation is
$\mathcal D_0$.
\end{remark}

\subsubsection*{Evidence for the local envelope}

For the local family class of Section~\ref{ss:family},
Theorems~\ref{thm:ae-energy}, \ref{thm:intrinsic-energy}, and
\ref{thm:padic-intrinsic-energy} give almost-everywhere finiteness
under their respective hypotheses. In zero volume, a smooth
family of interior minimizers induces the envelope metric.
Under the hypotheses of Proposition~\ref{prop:attained-family},
this metric is semipositive.

\subsection{Adelic models and approximations}
\label{ss:adelic-models-approximations}
\label{ss:family-critical}

Let $X\to C\to S$ be a family of fibrations over a global field
$K$, with polarization $\bar A$ and a relatively ample bundle $L$
on $X/C$. We first extend the family and $L$ to a projective model.
At both complex and nonarchimedean places, we extend metrics from a
fiber to the local family. For a projective family, we prove that
strong relative nefness is preserved when we replace the metric at
one place $v$ by a strongly relatively nef metric on $L_{K_v}$.

\subsubsection*{Arithmetic models}

Let $K$ be a global field. Write $T=\Spec\mathcal O_K$ if
$K$ is a number field, and let $T$ be the regular integral
projective curve with function field $K$ over its finite constant
field otherwise. For a finite extension $F/K$, let $T_F$ be
$\Spec\mathcal O_F$ or the normalization of $T$ in $F$,
respectively. Let $S$ be a smooth geometrically integral
quasi-projective $K$-variety, and let
\[
 X\xrightarrow{f}C\xrightarrow{\pi}S,\qquad g=\pi\circ f,
\]
be projective surjective morphisms. Assume that $f$ is flat of
relative dimension $d-1$, that $\pi$ and $g$ are smooth of relative
dimensions $1$ and $d$, and that $C_s,X_s$ are geometrically
integral for every geometric point $s$ of $S$.
Let $L$ be an $f$-ample rational line bundle.

\begin{lem}
\label{lem:polarized-arithmetic-spread}
There is a projective tower of integral schemes
\[
 \mathcal X\xrightarrow{\mathcal f}\mathcal C
       \xrightarrow{\mathcal\pi}\mathcal S
       \longrightarrow T
\]
with all morphisms flat and $\mathcal S$ normal, whose generic
base $S^c=\mathcal S_K$ is a normal projective compactification
of $S$. If $\operatorname{char}K=0$, $S^c$ may be chosen smooth.
The tower restricts over $S$ to the original one, and $L$ extends
to an $\mathcal f$-ample rational line bundle $\mathcal L$.

Every $a\in S(F)$ extends to a morphism $T_F\to\mathcal S$.
The pullback tower gives proper flat integral models
$\mathcal X_a\to\mathcal C_a\to T_F$, with
$\mathcal L_a$ relatively ample. There is a strongly relatively nef
adelic reference $\bar L_0=(L,h_0)$ whose finite-place metrics are the relatively
ample model metrics of $\mathcal L$ and whose archimedean metrics
are smooth in the ambient sense with strictly positive vertical curvature.
\end{lem}
\begin{proof}
We begin with a projective model $\mathcal S_0$ over $T$ of
a projective compactification of $S$. Spreading out the tower
and $L$ gives a projective tower
$X_U\xrightarrow{f_U}C_U\xrightarrow{\pi_U}U$, where
$U\subset\mathcal S_0$ is open, $U_K=S$, $f_U$ is flat,
and $\pi_U$ and $\pi_U\circ f_U$ are smooth. To extend $L$,
we choose a sufficiently divisible $n$. Relative Serre vanishing
and base change imply that
\[
 E_n=(f_U)_*L_U^{\otimes n}
\]
is locally free, and evaluation gives a closed immersion
$X_U\hookrightarrow\mathbb P_{C_U}(E_n)$.

We choose a projective closure
$\pi_0:\mathcal C_0\to\mathcal S_0$ of $C_U\to U$ and
extend $E_n$ to a coherent sheaf $\mathcal E_0$ on $\mathcal C_0$.
We then take the scheme-theoretic closure $\mathcal X_0$ of
$X_U$ in $\mathbb P_{\mathcal C_0}(\mathcal E_0)$ and put
\[
 \mathcal f_0:\mathcal X_0\longrightarrow\mathcal C_0,
 \qquad
 \mathcal L_0=
 \mathcal O_{\mathbb P_{\mathcal C_0}(\mathcal E_0)}(1)
       |_{\mathcal X_0}.
\]
The tautological bundle $\mathcal L_0$ is $\mathcal f_0$-ample
and extends $L_U^{\otimes n}$.

Blowing up $\mathcal C_0$ away from $C_U$ makes the pure
transform $\mathcal E_1$ locally free; Raynaud--Gruson flattening
then makes the pure transform
$\mathcal f_1:\mathcal X_1\to\mathcal C_1$ flat.
The quotient onto $\mathcal E_1$ embeds $\mathbb P(\mathcal E_1)$
into the base change of $\mathbb P(\mathcal E_0)$.
The closure $\mathcal X_1$ lies in $\mathbb P(\mathcal E_1)$,
so the pullback $\mathcal L_1$ remains relatively ample.

We flatten $\mathcal C_1\to\mathcal S_0$ by a blowup
$\mathcal S_1\to\mathcal S_0$ away from $U$. Writing
$\pi_2:\mathcal C_2\to\mathcal S_1$ for the pure transform, we set
\[
 \mathcal X_2=\mathcal X_1\times_{\mathcal C_1}\mathcal C_2,
 \qquad
 \mathcal f_2:\mathcal X_2\longrightarrow\mathcal C_2.
\]
Let $\mathcal L_2$ be the pullback of $\mathcal L_1$.
Then $\mathcal f_2$ and $\pi_2$ are flat and $\mathcal L_2$
is relatively ample. We have used the constructions of projective
models and flattening in \cite[\S4.5, proof of Theorem 4.1.3]{YZ26}.

In characteristic zero, we resolve $(\mathcal S_1)_K$ along
its boundary, extend this projective modification to
$\mathcal S_1$, and normalize. In positive characteristic, we
only normalize $\mathcal S_1$. Denote the resulting base by
$\mathcal S$ and put
\[
 \mathcal C=\mathcal C_2\times_{\mathcal S_1}\mathcal S,
 \qquad
 \mathcal X=\mathcal X_2\times_{\mathcal S_1}\mathcal S,
 \qquad
 \mathcal L=\frac1n p^*\mathcal L_2,
\]
where $p:\mathcal X\to\mathcal X_2$ is the projection.
Base change preserves flatness and relative ampleness. The
integral scheme $\mathcal S$ dominates the regular curve $T$
and is therefore flat over it. Thus $S^c=\mathcal S_K$ and
$\mathcal L$ have the required properties.
Flatness over the integral base and geometric integrality of
$X_\xi,C_\xi$ imply integrality of $\mathcal X,\mathcal C$:
their affine coordinate rings embed in their generic localizations,
which are domains.

For $a\in S(F)$, properness over $T$ extends $a$ over each
local ring of the regular one-dimensional scheme $T_F$; uniqueness
glues these extensions to $T_F\to\mathcal S$.
The pullback tower is proper and flat, with geometrically integral
generic fibers $X_a$ and $C_a$. The same argument proves that
$\mathcal X_a$ and $\mathcal C_a$ are integral.

It remains to choose the reference metric. We choose an ample
bundle $\mathcal N$ on $\mathcal C$ sufficiently positive that
\[
 \mathcal M=n\mathcal L+\mathcal f^*\mathcal N
\]
is an ample line bundle on $\mathcal X$, increasing the divisible
integer $n$ if necessary. At finite places, $\mathcal M$ and $\mathcal N$ carry their
model metrics. Over a number field, choose smooth Fubini--Study
metrics at the archimedean places and rescale them to make both
Hermitian bundles nef. Over a function field, ampleness on the
projective models already gives nef adelic bundles. Define $h_0$ by
\[
 (\mathcal L,h_0)
   =\frac1n\bigl(\overline{\mathcal M}
                    -\mathcal f^*\overline{\mathcal N}\bigr).
\]
Since $\overline{\mathcal N}$ is nef with ample generic restriction,
this decomposition makes $\bar L_0$ strongly relatively nef. At finite
places it is exactly the relatively ample model metric of
$\mathcal L$. At archimedean places it is smooth in the ambient
sense, with vertical curvature $1/n$ times the positive curvature
of $\overline{\mathcal M}$, including at singular points of fibers.
\end{proof}

\subsubsection*{Admissible collections of polarizations}

Let $K_0=\BQ$ or $\mathbb F_p(t)$, let $F/K_0$ be finitely
generated, and let $K$ be the algebraic closure of $K_0$ in $F$.
Let $f:X\to C$ be a flat projective surjective morphism of smooth
geometrically integral projective $F$-varieties, with $\dim X=d$
and $\dim C=1$. We take $E=\BC$ or a finite extension of
$\BQ_p$ when $K_0=\BQ$, and a finite extension of
$\mathbb F_p((t))$ when $K_0=\mathbb F_p(t)$, with the standard
embedding of $K_0$.

Let $S$ be a smooth geometrically integral quasi-projective $K$-model
of $F$, and let
\[
 \widetilde X\xrightarrow{\widetilde f}\widetilde C
       \xrightarrow{\widetilde\pi}S
\]
be a projective spread of $X\to C$. After shrinking $S$, we may assume
that $\widetilde f$ is flat and that $\widetilde X$ and $\widetilde C$
are smooth over $S$ with geometrically integral fibers. We fix a
projective flat model of this tower over
$\Spec\CO_K$, or over the smooth projective arithmetic base curve
of $K$ in positive characteristic. Its generic parameter base
$S^c$ is a normal projective compactification of $S$, smooth in
characteristic zero. On this model we choose a nef rational line bundle
$\bar A_0$ with ample underlying bundle. Over a number field its
archimedean metrics are smooth, semipositive, and compatible with
complex conjugation.

An embedding $\sigma\in\operatorname{Hom}_{K_0}(F,E)$ restricts to
$\iota\in\operatorname{Hom}_{K_0}(K,E)$ and gives a $K$-generic point
$s_\sigma\in S_\iota(E)$. Conversely, these generic points,
componentwise in $\iota$, are exactly the $K_0$-embeddings of $F$ into
$E$. Restriction and base change of $\bar A_0$ give
$\bar A_{0,\sigma}$ on $X_\sigma$, with the local absolute-value
normalizations fixed in the introduction.

\begin{defn}\label{def:admissible-polarizations}
A collection
$(\bar A_\sigma)_{\sigma\in\operatorname{Hom}_{K_0}(F,E)}$ is
\emph{admissible} if a single polarized model as above induces
every $\bar A_\sigma$ by the formula
\[
 \bar A_\sigma=(A_{0,\sigma}^d)^{-1/d}\bar A_{0,\sigma}.
\]
Thus $A_\sigma^d=1$ for every $\sigma$.
\end{defn}

Shrinking $S$ and pulling the fixed data back to a common model
do not change the induced collection. An auxiliary finite extension
of $K$, with a chosen embedding into $E$ extending $\iota$, uses
the pullback of these same data at the corresponding realization.

We can realize the prescribed polarization and the reference metric of
Lemma~\ref{lem:polarized-arithmetic-spread} on one tower. Write
$\mathcal X_0\to\mathcal C_0$ for the polarized model. We choose an
ample model bundle $G$ on $\mathcal X_0$ and $n$ large enough that
$nL-G$ is relatively ample on the generic fibration. After shrinking
the spread, we apply the construction of the lemma to $nL-G$, starting
with $\mathcal C_0$. This gives a tower $\mathcal X_1\to\mathcal C_1$
over a modification of the parameter base, with $\mathcal C_1\to
\mathcal C_0$ and a relatively ample rational extension $\mathcal M$ of $nL-G$.
We take the closure of the diagonal of the common open spread in
$\mathcal X_1\times_{\mathcal C_0}\mathcal X_0$, with projections
$p_1,p_0$. The bundle $(p_1^*\mathcal M+p_0^*G)/n$ extends $L$ and is
relatively ample over $\mathcal C_1$ by relative ampleness of the
external product. Flattening as in the lemma preserves relative
ampleness on the strict transform. The pullback of $\bar A_0$ induces the prescribed local polarizations. We take the reference from this extension of $L$.

\subsubsection*{Extension of metrics from a fiber}

\emph{Complex places.}
Suppose $K$ is a number field. For the tower at the start of this
subsection, fix an embedding $K\hookrightarrow\BC$, a normalized semipositive
polarization $\bar A$, and a smooth metric $h_0$ on $L$ with
strictly positive curvature on $T_{X/C}$.

\begin{lem}
\label{lem:metric-extension}
For every $s\in S(\BC)$ and $h_s\in\CH^+(L_s)$, there is
a smooth metric $h$ on $L$ with strictly positive vertical curvature
such that $h|_{X_s}=h_s$ and $h=h_0$ outside $g^{-1}(K_s)$ for
some compact $K_s\subset S(\BC)$.

For every $h_s\in\BH^+(L_s)$, there are such metrics $h_j$
whose restrictions converge uniformly to $h_s$ and satisfy
\[
 \Psi_s(h_{\rm ref},h_j|_{X_s})\longrightarrow \Psi_s(h_{\rm ref},h_s)
\]
for every integrable reference metric $h_{\rm ref}$ on $L_s$.
\end{lem}
\begin{proof}
Since $g$ is smooth and proper, $X_s$ is a compact smooth
submanifold of $X(\BC)$. We extend
$\varphi_s=\log(h_s/h_0|_{X_s})$ to a smooth real function
$\varphi$ on $X(\BC)$ and put $\widetilde h=h_0e^\varphi$.
To check vertical positivity near $X_s$, we use the closed set
\[
 \mathcal K_f=\{(x,[v])\in\mathbb P(T_X^\vee):df_x(v)=0\}.
\]
Properness of $\mathbb P(T_X^\vee)\to X$ makes the image of
the locus in $\mathcal K_f$ where $c_1(L,\widetilde h)$ is
nonpositive closed. It misses $X_s$: there
$\ker df_x\subset\ker dg_x$, and the curvature is that of $h_s$.
Thus $\widetilde h$ is vertically positive on a neighborhood
$W$ of $X_s$, even where the rank of $df$ changes.
Properness of $g$ gives $U\ni s$ with $g^{-1}(U)\subset W$.

We choose $0\le\chi\le1$ in $C_c^\infty(U)$ equal to one
near $s$ and put
\[
 h=h_0\exp\bigl((\chi\circ g)\varphi\bigr).
\]
On $\ker df$, the derivatives of $\chi\circ g$ vanish.
Thus the curvature of $h$ is the convex combination of those
of $h_0$ and $\widetilde h$, with coefficients $1-\chi\circ g$
and $\chi\circ g$. Both are positive where the second coefficient is nonzero,
so $h$ has the required positivity, restriction, and support.

For $h_s\in\BH^+(L_s)$, we use
Lemma~\ref{lem:relative-nef-approximation} to choose
$h_{s,j}\in\CH^+(L_s)$ converging uniformly and in energy
to $h_s$, and extend each $h_{s,j}$ as above. The cocycle identity
gives convergence for every $h_{\rm ref}$.
\end{proof}

At a real place we take $h_0$ invariant under conjugation. If
$s=\bar s$ and $h_s$ is invariant, averaging the logarithmic change
of metric with its conjugate gives an invariant extension. For
$s\ne\bar s$, disjoint conjugate neighborhoods let us prescribe
$h_s$ at $s$ and its conjugate at $\bar s$.

\smallskip\noindent\emph{Nonarchimedean places.}
Let $E$ be a nonarchimedean local field and
$X\xrightarrow f C\xrightarrow\pi S$ a tower of projective
$E$-varieties, with $f$ flat and $L$ an $f$-ample rational
line bundle. Fix $s\in S(E)$ such that $C_s$ is a smooth
geometrically connected curve and $X_s$ is geometrically integral,
and equip $(X_s/C_s)$ with a fixed normalized polarization
$\bar A_s$, obtained from an ample rational bundle with a
semipositive metric, so that $A_s^d=1$.

\begin{lem}
\label{lem:nonarch-metric-extension}
Every metric $h_s\in\CH^+(L_s)$ extends to a metric $h$ on $L$
induced by a relatively ample rational model of $f$:
\[
 h|_{X_s}=h_s.
\]
The representing morphism of models may be taken flat.

For every $h_s\in\BH^+(L_s)$, there are such metrics $h_j$
whose restrictions converge uniformly to $h_s$ and satisfy
\[
 \Psi_s(h_{\rm ref},h_j|_{X_s})\longrightarrow \Psi_s(h_{\rm ref},h_s)
\]
for every integrable reference metric $h_{\rm ref}$ on $L_s$.
\end{lem}
\begin{proof}
We first extend $h_s\in\CH^+(L_s)$ by writing it as the
difference of two ample model metrics. We choose a relatively ample
rational model inducing $h_s$ and a regular projective model
$\CC_s\to\Spec\CO_E$ of $C_s$ dominating its base;
see \cite[Tag 0C2R, Lemma 55.8.3]{Stacks}. By base change and
strict transform, we obtain a relatively ample rational bundle
$\mathcal L_s$ on $\mathcal f_s:\mathcal X_s\to\CC_s$
still inducing $h_s$.

We choose an ample bundle $N$ on $C$ and construct an ample rational
extension $\mathcal N_s$ of $N_s$ to $\CC_s$. To construct it,
we start with any extension and write the special fiber of
$\CC_s\to\Spec\CO_E$ as
$\sum_i b_iD_i$, with irreducible components $D_i$ and
multiplicities $b_i$. The matrix $(D_i\cdot D_j)$ has kernel
spanned by $(b_i)$, and
$\sum_i b_i\deg(\mathcal N_s|_{D_i})=\deg N_s>0$.
By adding a vertical rational divisor, we can prescribe any
positive rational vector of component degrees with the same
weighted sum.
The resulting bundle is ample. Here we use
\cite[Tag 0C5Y, Lemmas 55.9.5 and 55.9.7]{Stacks} and the
criterion for ampleness on curves in terms of degrees on
components \cite[Tag 0AYQ, Lemma 33.44.15]{Stacks}.

Choose $m$ such that $M=L+mf^*N$ and
$\mathcal L_s+m\mathcal f_s^*\mathcal N_s$ are ample.
With $N'=mN$ and $n_s$ induced by $m\mathcal N_s$, put
\[
 a_s=h_s+f_s^*n_s.
\]
Thus $h_s=a_s-f_s^*n_s$ is a difference of ample model metrics.

We extend $a_s$ and $n_s$ separately by extending their
Fubini--Study sections. For sufficiently divisible $r$, a tuple
$(\sigma_i)$ of sections of $M_s^{\otimes r}$ without common
zeros defines $ra_s$; see \cite[Theorem 5.14]{BE21}.
Increasing $r$, we may assume that $M^{\otimes r}$ is very
ample and, by Serre vanishing for the ideal of $X_s\subset X$,
that
\[
 H^0(X,M^{\otimes r})\longrightarrow H^0(X_s,M_s^{\otimes r})
\]
is surjective. We lift the $\sigma_i$ to $X$. For a continuous
reference norm on $M_s^{\otimes r}$, compactness gives
\[
 \rho=\min_{x\in X_s^{\an}}\max_i|\sigma_i(x)|>0.
\]
To make the lifted sections define an immersion, we append
uniformizer multiples of a basis of $H^0(X,M^{\otimes r})$,
chosen so small that their norms on $X_s$ are less than $\rho$.
The added sections leave the maximum on $X_s$ unchanged and
make the map to projective space a closed immersion. Taking
its closure over $\CO_E$, with $\mathcal O(1)/r$, gives an
ample rational model metric $a$ on $M$ with $a|_{X_s}=a_s$.
The same construction extends $n_s$ to an ample model metric
$n$ on $N'$. Then $h=a-f^*n$ restricts to $h_s$.
To represent $h$ by a relatively ample model of $f$, let
$(\mathcal X_M,\mathcal M)$ and
$(\mathcal C_{N'},\mathcal N')$ be the projective $\CO_E$-models
defining $a,n$. We take the scheme-theoretic closure
$\mathcal Z$ of $\Gamma_f\subset X\times_E C$ in
$\mathcal X_M\times_{\CO_E}\mathcal C_{N'}$ and write
\[
 p=\operatorname{pr}_1|_{\mathcal Z}:
       \mathcal Z\longrightarrow\mathcal X_M,\qquad
 \mathcal f=\operatorname{pr}_2|_{\mathcal Z}:
       \mathcal Z\longrightarrow\mathcal C_{N'}.
\]
The generic fibers are $p_E=\operatorname{id}_X$ and
$\mathcal f_E=f$, under $\mathcal Z_E=\Gamma_f\simeq X$.
The bundle
\[
 \mathcal L=p^*\mathcal M-\mathcal f^*\mathcal N'
\]
induces $h$. It is $\mathcal f$-ample because
$p^*\mathcal M$ is relatively ample over the second factor,
and twisting from the base preserves relative ampleness.
We flatten $\mathcal f$ by modifying its base and taking the
strict transform. Relative ampleness and the metric are preserved;
see \cite[Proposition 5.13]{BE21}.

For $h_s\in\BH^+(L_s)$, we choose
$h_{s,j}\in\CH^+(L_s)$ converging uniformly and in energy
to $h_s$ by Lemma~\ref{lem:relative-nef-approximation}, and
extend each $h_{s,j}$. The cocycle identity gives convergence
for every $h_{\rm ref}$.
\end{proof}

In both lemmas the infimum of the energy over the restrictions
equals its infimum over $\BH^+(L_s)$. Indeed, every restriction
belongs to $\BH^+(L_s)$, and the energies of the restrictions
approximate the energy of every metric in that class. By the
change-of-metric formula for Deligne pairings, the induced norms on
$D_A(L)_s$ likewise approximate the norm of every $h_s\in\BH^+(L_s)$
and have the same pointwise supremum
\cite[Theorem 4.2.3 and \S4.2.2]{YZ26}.

\subsubsection*{Adelic replacement at one place}

Let $f:X^c\to C^c$ be a projective morphism of projective
varieties over a global field $K$, let $L$ be $f$-ample, and fix
a strongly relatively nef adelic reference $\bar L_0$ on $L$.
For a place $v$ and a metric $h_v$ on $L_{K_v}$, write
$\bar L_0^v(h_v)$ for the metric obtained by replacing the
$v$-component of $\bar L_0$ by $h_v$. At a complex place we
prescribe $h_v$ at one embedding and its conjugate at the other.
At a real place we use conjugation-invariant metrics on the
complexification.

\begin{cor}\label{cor:adelic-replacement}
\label{cor:nonarch-adelic-replacement}
At every place $v$ of $K$, if $h_v$ is strongly relatively nef
for $f_{K_v}$, then $\bar L_0^v(h_v)$ is strongly relatively nef
for $f$.
\end{cor}
\begin{proof}
Write $\bar L_0=\bar M_0-f^*\bar N_0$ with both factors nef and
$N_0$ ample. Strong relative nefness of $h_v$ gives an ample
bundle $N_v$ on $C^c_{K_v}$ with a semipositive metric $n_v$
such that $h_v+f^*n_v$ is semipositive. For large $m$, both
$L+mf^*N_0$ and $mN_{0,K_v}-N_v$ are ample. Adding a semipositive
metric on the latter to $n_v$ gives a metric $n$ on $mN_{0,K_v}$
for which $n$ and $h_v+f^*n$ are semipositive. At a real place
we may average the weights of $n$ with their conjugates, since
$h_v$ is invariant.

We give $mN_0$ the metric $n$ at $v$ and the metric of $m\bar N_0$
elsewhere, obtaining $\bar N$. Both $\bar N$ and
$\bar M=\bar L_0^v(h_v)+f^*\bar N$ have ample underlying bundles
and semipositive metrics at every place. As in
Lemma~\ref{lem:one-place-globalization}, a common twist from
$\Spec K$ makes both factors nef and cancels in $\bar M-f^*\bar N$.
This argument applies in every dimension of $C^c$.
\end{proof}

\subsection{Extension to complex points}\label{ss:complex-points}

We prove Theorem~\ref{thm:main-energy} for complex embeddings.
To pass from algebraic points to almost every complex point of $S$,
we use Latt\`es maps and finite projections to construct a measure
that dominates smooth volume. We obtain it as a limit of averages
of Galois orbit measures of points with uniformly bounded height.

Let $X\xrightarrow{f}C\xrightarrow{\pi}S$ be a tower over a
number field $K$ satisfying the geometric hypotheses at the start
of Section~\ref{ss:adelic-models-approximations}, with $g=\pi\circ f$
and an $f$-ample rational line bundle $L$. We equip the
family with a fixed normalized polarization $\bar A$, obtained from a nef
Hermitian rational line bundle $\bar A_0$ on an arithmetic
compactification as in Lemma~\ref{lem:polarized-arithmetic-spread},
with ample underlying bundle and smooth semipositive archimedean
metrics. We fix these polarized model data, with $A_s^d=1$.

Let $\xi$ be the generic point of $S$.
We fix an embedding $K\hookrightarrow\mathbb C$ and a smooth metric
$h_0$ on $L$ whose curvature is strictly positive on every nonzero
vector in the relative tangent space $T_{X/C,x}$. For $s\in S(\mathbb C)$, put
\[
 m_\infty(s)=\inf_h \Psi_s(h_{0,s},h),\qquad
 h_{0,s}=h_0|_{X_s},
\]
where $h$ runs over the smooth class $\CH^+(L_s)$.

\begin{thm}\label{thm:ae-energy}
If $L_\xi$ is nef, then
\[
 m_\infty\in L^1_{\mathrm{loc}}(S(\mathbb C)),
\]
and the same infimum is obtained over $\BH^+(L_s)$.
\end{thm}

The proof also shows that $m_\infty\le0$ is upper semicontinuous.
If $f$ is smooth, its value is the infimum over all smooth
fiberwise-semipositive metrics. Outside one set of measure zero,
the energy is bounded below on each of these classes; for the
fixed polarized model, this set is independent of the metric
over which the infimum is taken.

For the model reference of Lemma~\ref{lem:polarized-arithmetic-spread},
the proof gives the stronger assertion
\[
 m_\infty\in L^1\bigl(S(\mathbb C),dV^{\mathrm c}|_S\bigr)
\]
for every smooth positive volume form $dV^{\mathrm c}$ on the
chosen smooth projective compactification $S^{\mathrm c}$.
If $S$ is projective, this holds for every smooth reference with
strictly positive vertical curvature.

The proof extends metrics from parameter fibers, applies the
height inequality to points with bounded height, and glues
metrics to pass from an average bound to the infimum.
Perturbation by $\epsilon A_0$ allows us to use generic nefness.

Let $S/K$ be smooth, geometrically integral, and projective,
with ample height $h_B$. Fix $K\hookrightarrow\BC$; all
embeddings in the following orbit measures extend it.

\begin{lem}
\label{lem:bounded-height-volume}
For every smooth positive volume form $dV$
on $S(\BC)$, there are constants $H,c>0$ and a probability measure
$\mu$ such that $\mu\ge c\,dV$ and $\mu$ belongs to the weak
closed convex hull of the normalized orbit measures
\[
 \nu_a=\frac1{[K(a):K]}
      \sum_{\tau:K(a)\hookrightarrow\BC}\delta_{\tau(a)},
 \qquad h_B(a)\le H.
\]
Moreover, $\mu$ gives zero mass to every proper closed
$K$-subvariety of $S$. For any dense $K$-open subset $U\subset S$,
the same $\mu$ belongs to the weak closed convex hull using only
points $a\in U(\bar K)$ with the same bound $h_B(a)\le H$.
\end{lem}
\begin{proof}
\smallskip\noindent\emph{A measure on projective space.}
The assertion is immediate when $\dim S=0$. For $n=\dim S>0$,
we first construct the required measure on $\BP^n$ using a
Latt\`es map. We fix an elliptic curve $E/K$ with origin $0$ and let
\[
 \rho:E\longrightarrow E/\{\pm1\}\simeq\BP^1
\]
be the quotient map. Multiplication by $2$ descends to an
endomorphism $\varphi:\BP^1\to\BP^1$ of degree $4$, satisfying
$\varphi\circ\rho=\rho\circ[2]$. The product
$\varphi^{\times n}$ commutes with permutations of the factors,
so it descends through
\[
 \sigma:(\BP^1)^n\longrightarrow
       \operatorname{Sym}^n(\BP^1)\simeq\BP^n
\]
to an endomorphism $\Phi:\BP^n\to\BP^n$ of algebraic degree
$4$. We put $A=E^n$ and $v=\sigma\circ\rho^n$ and obtain
\[
 \begin{tikzcd}
 A \arrow[r,"{[2]}"] \arrow[d,"v"'] & A \arrow[d,"v"]\\
 \BP^n \arrow[r,"\Phi"'] & \BP^n,
 \end{tikzcd}
 \qquad
 \deg v=\delta:=2^n n!,\qquad \deg\Phi=4^n.
\]
Here $\deg\Phi$ is its degree as a finite morphism
\cite[Propositions 1.2 and 4.1]{GHK23}.

Let $m_E$ and $m_A=m_E^{\otimes n}$ be the Haar probability
measures. Since the equilibrium measure of $\varphi$ is
$\rho_*m_E$ \cite[\S7]{Milnor06}, the formula for the measure
on the symmetric product gives
\[
 \lambda:=\mu_\Phi
   =\sigma_*\bigl((\rho_*m_E)^{\otimes n}\bigr)=v_*m_A;
\]
see \cite[proof of Lemma 1.3]{GHK23}. We check that $\lambda$
dominates a smooth positive volume form. Let $\Omega$ be the
normalized Fubini--Study volume on $\BP^n$. By compactness,
$v^*\Omega\le c_0m_A$ for some $c_0>0$. Applying the
change-of-variables formula $v_*v^*\Omega=\delta\Omega$, we obtain
\[
 \lambda\ge\frac{\delta}{c_0}\Omega.
\]
Also, $\lambda$ gives zero mass to every proper algebraic
subset of $\BP^n$, since its inverse image under $v$ is a
proper algebraic subset of $A$ and has Haar measure zero.

To approximate $\lambda$ by algebraic cycles, we use the
inverse images of the fixed point $b=v(0)\in\BP^n(K)$:
\[
 D_N=(\Phi^N)^*[b]\qquad(N\ge1),
\]
with scheme-theoretic multiplicities. The finite maps $v$ and
$\Phi$ are flat, since their sources and targets are smooth of the
same dimension \cite[Tag 00R4]{Stacks}. Since $v^{-1}(b)$ is
supported at $0$, one has $v^*[b]=\delta[0]$. Hence
\[
 v^*D_N=[2^N]^*v^*[b]=\delta[A[2^N]],
 \qquad D_N=v_*[A[2^N]],
\]
where the second equality follows by pushing forward and using
$v_*v^*=\delta$. In particular, $\deg D_N=4^{nN}$.
The normalized counting measures on the torsion grids
$A[2^N](\BC)$ converge to $m_A$, so the probability measures
associated with $4^{-nN}D_N$ converge to $\lambda$.

For a Weil height $h$ of $\CO_{\BP^n}(1)$, the
Call--Silverman height \cite{CS93} is
\[
 \widehat h_\Phi(y)=\lim_{m\to\infty}4^{-m}h(\Phi^m(y)).
\]
It satisfies $\widehat h_\Phi\circ\Phi=4\widehat h_\Phi$
and $|h-\widehat h_\Phi|\le c_\Phi$ for some $c_\Phi\ge0$
(see also \cite[\S4.1]{Zh06} and \cite[\S2]{Zh95b} for the
adelic construction).
Every support point of $D_N$ is the image of a torsion point
and is preperiodic, so
\[
 h(y)\le c_\Phi
\]
for every such point and every $N$.

\smallskip\noindent\emph{Pullback to the parameter space.}
We pull the measures on $\BP^n$ back to $S$ by finite linear
projections. For each $x\in S(\BC)$, we can choose a finite
$K$-projection $u:S\to\BP^n$ that is \'etale at $x$.
Indeed, centers avoiding $S$ and the embedded tangent space
at $x$ form a nonempty Zariski-open subset of the Grassmannian,
and the $K$-rational centers are Zariski dense over $\BC$.
By compactness, we can choose finitely many such maps
$u_i:S\to\BP^n$ whose \'etale loci cover $S(\BC)$.
We write $r_i=\deg u_i$. As above, smoothness of the source
and target implies that $u_i$ is flat. The effective zero-cycles
\[
 Z_{i,N}=u_i^*D_N
\]
are defined over $K$ and have degree $r_i4^{nN}$.
For every point $s$ in their support, height comparison under the
finite morphism gives
\[
 h_B(s)\le a_i h(u_i(s))+b_i
       \le a_i c_\Phi+b_i.
\]
Taking the maximum over the finitely many $i$ gives one height
bound $H>0$, independent of $N$.

For a continuous function $\psi$ on $S(\BC)$, the fiber trace
\[
 (\operatorname{Tr}_{u_i}\psi)(y)
 =\sum_{x\in u_i^{-1}(y)}
          \operatorname{length}_x(u_i^{-1}(y))\,\psi(x)
\]
is continuous: properness confines nearby fibers to disjoint
neighborhoods of the original fiber points, with constant total
multiplicities. Thus the measures of $Z_{i,N}/(r_i4^{nN})$ converge to
\[
 \mu_i=\frac1{r_i}u_i^*\lambda
      =\frac1{r_i}u_i^*(v_*m_A).
\]
Pullback of measures is defined by this trace. The normalized
measure of $Z_{i,N}$ is a convex combination of $\nu_a$, weighted
by cycle multiplicities and residue-field degrees.

The inequality for $\lambda$ implies
\[
 \mu_i\ge\frac{\delta}{c_0r_i}u_i^*\Omega.
\]
For the smooth form $\Omega$, pullback by the trace operator
agrees with differential pullback. Since the \'etale loci of
the $u_i$ cover $S(\BC)$, the smooth nonnegative top form
$\sum_i r_i^{-1}u_i^*\Omega$ is everywhere positive. Its ratio
to $dV$ has a positive minimum. We take $\mu$ to be the finite
average of the $\mu_i$. Then $\mu\ge c\,dV$ for some $c>0$,
and $\mu$ is approximated by convex combinations of orbit
measures of points with height at most $H$.

For a proper closed $K$-subvariety $Z\subsetneq S$, finiteness
of $u_i$ implies that $u_i(Z)$ is proper in $\BP^n$. Hence
$\lambda(u_i(Z))=0$, and the trace formula gives
$\mu_i(Z)=0$. Averaging over $i$, we obtain $\mu(Z)=0$.

For any prescribed dense $K$-open $U\subset S$, delete the
whole Galois orbits in $S\setminus U$. Their mass tends to zero
by Portmanteau, since $\mu(S\setminus U)=0$. Renormalization
preserves the weak limit and the bound $H$ on every support point.
\end{proof}

For the gluing lemma, retain the geometric setting, polarization,
embedding, and reference of Theorem~\ref{thm:ae-energy}. Let $\mathcal H_c$ consist of the smooth
metrics with strictly positive vertical curvature that agree
with $h_0$ outside $g^{-1}(K_h)$ for some compact
$K_h\subset S(\mathbb C)$. We put
\[
 E_h(s)=\Psi_s(h_{0,s},h|_{X_s}),\qquad
 e_\infty(s)=\inf_{h\in\mathcal H_c}E_h(s).
\]
Let $\mu$ be a Radon probability on a projective compactification
$S^c(\mathbb C)$, with $\mu(S^c\setminus S)=0$.
We extend each $E_h$ by zero to $S^c$.

\begin{lem}
\label{lem:energy-envelope-gluing}
If
\[
 \int E_h\,d\mu\ge-B\qquad(h\in\mathcal H_c)
\]
for one constant $B\ge0$, then
\[
 \int_{S(\mathbb C)}e_\infty\,d\mu\ge-B.
\]
\end{lem}
\begin{proof}
The functions $E_h$ are continuous by smoothness and properness
of $g$. They have compact support, and $E_{h_0}=0$. We first
choose a countable family
$\{h_0,h_1,h_2,\ldots\}\subset\mathcal H_c$ such that
\[
 e_\infty(s)=\inf_{j\ge0}E_{h_j}(s)
 \qquad(s\in S(\mathbb C)).
\]
For each rational $a$, second countability gives a countable
subcover of $\{e_\infty<a\}$ by sets $\{E_h<a\}$.
The corresponding metrics, together with $h_0$, give this family.
The finite minima
\[
 F_k=\min(0,E_{h_1},\ldots,E_{h_k})
\]
decrease to $e_\infty\le0$. By monotone convergence applied to
the nonnegative increasing sequence $-F_k$, it suffices to prove
$\int F_k\,d\mu\ge-B$ for every $k$.

For fixed $k$, it suffices to approximate $F_k$ in $L^1(\mu)$
by the energy of a single metric; compare \cite[Theorem 14.60]{RW98}.
Let $W_j$ be the measurable set where $j$ is the smallest index
attaining $F_k$, with index zero corresponding to $h_0$.
By inner regularity, we choose compact subsets
$K_j\subset W_j$ whose union has measure at least $1-\delta$.
The $K_j$ are pairwise disjoint. For $1\le j\le k$, we choose
smooth compactly supported functions $0\le\chi_j\le1$, equal
to one on $K_j$, whose supports are disjoint and avoid $K_i$
for every $i\ne j$, including $K_0$. With
$\varphi_j=\log(h_j/h_0)$, we put
\[
 h_\chi=h_0\exp\!\left(\sum_{j=1}^k
                         (\chi_j\circ g)\varphi_j\right).
\]
The derivatives of $\chi_j\circ g$ vanish on the relative
tangent spaces of $f$. The curvature there is a convex
combination of strictly positive curvatures, so
$h_\chi\in\mathcal H_c$. Its energy equals $F_k$ on
$\bigcup_{j=0}^k K_j$.

To control the energy on the complement, we use the functions
\[
 (s,u)\longmapsto
 \Psi_s\bigl(h_{0,s},(h_0\exp(u\varphi_j))|_{X_s}\bigr),
 \qquad 1\le j\le k,\quad u\in[0,1],
\]
which vanish outside a common compact set
$K\subset S(\mathbb C)$. By continuity on $K\times[0,1]$,
their absolute values have a common bound $C_k$, depending only
on the fixed metrics and independent of the cutoffs and $\delta$.
At most one cutoff is nonzero on each fiber. Hence both
$|E_{h_\chi}|$ and $|F_k|$ are bounded by $C_k$, and
\[
 \|E_{h_\chi}-F_k\|_{L^1(\mu)}\le2C_k\delta.
\]
Letting $\delta\to0$ gives the required bound
$\int F_k\,d\mu\ge-B$, completing the proof.
\end{proof}

Under the hypotheses of Theorem~\ref{thm:ae-energy}, take the
model, reference $h_0$, and compactification $S^c$ of
Lemma~\ref{lem:polarized-arithmetic-spread}. Let $\mu$ be the
measure of Lemma~\ref{lem:bounded-height-volume}, with bound
$H_*$ for an ample height $h_B\ge0$ on $S^c$.

\begin{lem}
\label{lem:generic-nef-average}
With $\mathcal H_c$ and $E_h$ as in
Lemma~\ref{lem:energy-envelope-gluing}, there is a constant
$B_0\ge0$, depending on the fixed data above and independent
of $h\in\mathcal H_c$, such that
\[
 \int_{S(\BC)}E_h\,d\mu\ge-B_0
 \qquad(h\in\mathcal H_c).
\]
\end{lem}
\begin{proof}
\smallskip\noindent\emph{Changing the metric at one place.}
We work on the projective flat arithmetic tower of
Lemma~\ref{lem:polarized-arithmetic-spread}. Let $\bar L_0$
be its relatively ample Hermitian model bundle extending $L$.
We keep the nef Hermitian rational bundle $\bar A_0$ representing
the polarization $\bar A$ of Theorem~\ref{thm:ae-energy}.

To vary an arbitrary complex metric at a real embedding, we
put $K'=K(i)$ when the chosen embedding is real, and $K'=K$
otherwise. We fix an extension of the embedding to $K'$ and
denote its complex place by $v$. After pulling the model and
polarization back to $K'$, we change the metric at $v$, with
the conjugate change at its conjugate embedding, and use the
Haar normalization. To retain the original measure $\mu$, we
decompose the $K$-orbit of each $a\in S(\bar K)$ into
$K'$-orbits represented by $a_i$. We write
\[
 e=[K(a):K],\qquad e_i=[K'(a_i):K'],\qquad
 \nu_a=\sum_i\frac{e_i}{e}\nu'_{a_i},\qquad \sum_i e_i=e.
\]
The measures $\nu'_{a_i}$ use the chosen embedding of $K'$.
Since absolute heights over $K'$ have denominator $[K':\BQ]$,
we will include the factor $[K':K]$ when comparing them with
the original $K$-orbit averages.
For $h\in\mathcal H_c$, the function
\[
 \varphi_h=\log\frac{|\cdot|_h}{|\cdot|_{h_0}}
\]
vanishes outside $g^{-1}(K_h)$ for a compact
$K_h\subset S(\BC)$. Its support is compact by properness of
$g$, so we extend it smoothly by zero to the compactified total
space. The extended metric has strictly positive vertical curvature, since
it agrees with $h_0$ near the boundary. An ample base twist with
sufficiently positive curvature makes it semipositive on the
compactification. Corollary~\ref{cor:adelic-replacement} therefore
gives the strongly relatively nef adelic metric $\bar L_0^v(h)$.
The finite-place models are those of $\bar L_0$.

We now perturb $L$ to apply the height inequality at algebraic
points. For rational $0<\epsilon\le1$, we put
\[
 \bar L_\epsilon(h)=\bar L_0^v(h)+\epsilon\bar A_0,
 \qquad L_\epsilon=L+\epsilon A_0.
\]
The bundle $L_{\epsilon,\xi}$ is ample, since $L_\xi$ is nef
and $A_{0,\xi}$ is ample. We can therefore choose a dense
$K$-open subset $U_\epsilon\subset S$, depending on $\epsilon$,
over which $L_\epsilon$ is ample relative to $g$. By flatness,
\[
 v_\epsilon=(L_\epsilon|_{X_s})^d
 =\sum_{j=0}^d\binom dj\epsilon^j L_\xi^{d-j}A_{0,\xi}^j,
 \qquad v_\epsilon\longrightarrow L_\xi^d.
\]
The numbers $v_\epsilon$ are independent of $s$ and bounded
for $0<\epsilon\le1$. We perturb by the rational bundle $A_0$
and keep $\bar A$ fixed in the mixed pairings, which are rational
by our normalization. We consider
\[
 \bar R_{\epsilon,h}=D_A(\bar L_\epsilon(h)).
\]
At every $a\in U_\epsilon(F)$ over a finite extension $F/K'$,
the restriction $\bar L_{\epsilon,h,a}$ is relatively nef by base
change and its underlying bundle is ample on $X_a$.
Theorem~\ref{thm:ahi} and \eqref{eq:H-family-pairing} give
\[
 h_{\bar R_{\epsilon,h}}(a)
 =\frac1{[F:\BQ]}\left(
   -\bar L_{\epsilon,h,a}^{d+1}
   +(d+1)v_\epsilon
      (\bar L_{\epsilon,h,a}\cdot\bar A_a^d)
                      \right)\ge0.
\]

\smallskip\noindent\emph{The average energy bound.}
We compare these nonnegative heights with the reference heights.
By multilinearity, $\bar R_{\epsilon,h_0}$ is a polynomial in
$\epsilon$ whose coefficients are finitely many fixed Hermitian
model pairings on the projective base. Ample height comparison
therefore gives a constant $\gamma\ge0$ such that
\[
 [K':K]\,h_{\bar R_{\epsilon,h_0}}(a)\le \gamma(h_B(a)+1)
 \qquad(0<\epsilon\le1).
\]
We write $h_s=h|_{X_s}$ and let $h_{A_0}$ be the metric of
$\bar A_0$. We denote by $E_{\epsilon,h}^{A}(s)$ the corrected
local energy on $L_\epsilon$, with the normalized polarization
$\bar A_s$ fixed:
\[
 E_{\epsilon,h}^{A}(s)
 =\Psi_s(h_{0,s}\otimes h_{A_0,s}^{\otimes\epsilon},
             h_s\otimes h_{A_0,s}^{\otimes\epsilon}).
\]
The change-of-metric formula for $\bar R_{\epsilon,h}$ gives
\[
 \begin{aligned}
 \frac1{[K:\BQ]}\int E_{\epsilon,h}^{A}\,d\nu_a
 &=[K':K]\sum_i\frac{e_i}{e}
   \bigl(h_{\bar R_{\epsilon,h}}(a_i)
          -h_{\bar R_{\epsilon,h_0}}(a_i)\bigr)\\
 &\ge-B_*,
 \end{aligned}
\]
whenever $a\in U_\epsilon(\bar K)$ and $h_B(a)\le H_*$, where
$B_*=\gamma(H_*+1)$ is independent of $h\in\mathcal H_c$, of
rational $0<\epsilon\le1$, and of $a$ satisfying these conditions.
Indeed, $U_\epsilon$ is defined over $K$, so every $a_i$ belongs
to $U_\epsilon$ and $h_B(a_i)=h_B(a)$.

We fix $h$ and $\epsilon$ and approximate $\mu$ by convex
combinations of orbit measures of points $a\in U_\epsilon$
with height at most $H_*$. The function $E_{\epsilon,h}^{A}$
is continuous, supported in $K_h$, and extends continuously
by zero to $S^c(\BC)$. We can thus pass to the weak limit:
\[
 \frac1{[K:\BQ]}\int E_{\epsilon,h}^{A}\,d\mu\ge-B_*.
\]
With $\alpha_0=c_1(L,h_0)$, $\alpha_h=c_1(L,h)$, and
$\omega_{A_0}=c_1(\bar A_0)$, this energy is
\[
 \begin{aligned}
 &\int_{X_s}\varphi_h
    \sum_{i=0}^d(\alpha_0+\epsilon\omega_{A_0})^i
                 (\alpha_h+\epsilon\omega_{A_0})^{d-i}\\
 &\qquad-(d+1)v_\epsilon\int_{X_s}\varphi_h\,c_1(\bar A_s)^d.
 \end{aligned}
\]
For fixed $h$, this is a polynomial in $\epsilon$ with continuous
coefficients supported in $K_h$. It converges uniformly to
$E_h(s)=\Psi_s(h_{0,s},h_s)$ as $\epsilon\to0$. Hence
\[
 \int E_h\,d\mu\ge-[K:\BQ]B_*.
\]
We obtain the assertion with $B_0=[K:\BQ]B_*$. For each
$\epsilon$ we used points in its own open set $U_\epsilon$
to take the limit to $\mu$, and only afterwards let $\epsilon$
tend to zero.
\end{proof}

\begin{proof}[Proof of Theorem~\ref{thm:ae-energy}]
We first use the model reference of
Lemma~\ref{lem:polarized-arithmetic-spread}. We choose $\mu$ as in
Lemma~\ref{lem:bounded-height-volume}, dominating a smooth
positive volume form on $S^c$.
By Lemma~\ref{lem:generic-nef-average}, all $h\in\mathcal H_c$
satisfy the same integral bound $-B_0$. Applying
Lemma~\ref{lem:energy-envelope-gluing}, we obtain
\[
 \int_{S(\BC)}(-e_\infty)\,d\mu\le B_0,
 \qquad e_\infty(s)=\inf_{h\in\mathcal H_c}E_h(s)\le0.
\]
Since $\mu\ge c\,dV^c$ on $S^c(\BC)$, we have
$e_\infty\in L^1(S(\BC),dV^c|_S)$. In particular,
$e_\infty$ is locally integrable for every smooth volume form
on $S(\BC)$. By Lemma~\ref{lem:metric-extension},
$e_\infty(s)$ is the infimum over $\CH^+(L_s)$ at every $s$,
including when
$f_s$ has singular fibers. By
Lemma~\ref{lem:relative-nef-approximation}, this is also
$m_\infty(s)$, the infimum over $\BH^+(L_s)$. When $f$ is
smooth, we interpolate a smooth fiberwise-semipositive metric
with $h_{0,s}$ to obtain the additional equality of infima.

As an infimum of continuous functions, $m_\infty$ is upper
semicontinuous. For another smooth reference $h'_0$ with strictly
positive vertical curvature, the cocycle identity gives
\[
 m_{h'_0}(s)=m_{h_0}(s)+\Psi_s(h'_{0,s},h_{0,s}).
\]
The function $\Psi_s(h'_{0,s},h_{0,s})$ is smooth. Thus the new
infimum is locally integrable and upper semicontinuous. When
$S$ is projective, this function is bounded, and we also obtain
global integrability for $h'_0$.
\end{proof}

\subsubsection*{An intrinsic formulation over a finitely generated field}

Let $F$ be a finitely generated extension of $\BQ$, and let
$K$ be the algebraic closure of $\BQ$ in $F$. For an embedding
$\sigma:F\hookrightarrow\BC$, write $X_\sigma/C_\sigma$ and
$L_\sigma$ for the corresponding realizations. Their polarizations
form an admissible collection $(\bar A_\sigma)$ in the sense of
Definition~\ref{def:admissible-polarizations}.

To interpret almost every embedding, we choose a smooth
geometrically integral quasi-projective $K$-model $S$ with
$K(S)=F$. For each $\iota:K\hookrightarrow\BC$, embeddings
extending $\iota$ identify with the $K$-generic points of
$S_\iota(\BC)$; see also \cite[Lemma 3.1.1 and its proof]{YZ26}.
Their complement is a countable union of proper $K$-subvarieties,
so it has measure zero for smooth volume. Any two models have
isomorphic dense open subsets, with complements of measure zero.
Thus they define the same null sets on the space of embeddings.
We use this measure class on each of the finitely many components
indexed by the embeddings of $K$.

For the following theorem, let $F/\BQ$ be finitely generated,
and let $f:X\to C$ be a
flat projective surjective morphism of smooth geometrically
integral projective varieties over $F$, with $\dim X=d$ and
$\dim C=1$. Fix an admissible collection $(\bar A_\sigma)$
and one set of polarized model data inducing it, as in
Definition~\ref{def:admissible-polarizations}.

\begin{thm}
\label{thm:intrinsic-energy}
Let $L$ be a geometrically nef, $f$-ample $\BQ$-line bundle.
For almost every embedding $\sigma:F\hookrightarrow\BC$,
one has, for every integrable reference metric $h_{0,\sigma}$
on $L_\sigma$,
\[
 \inf_{h\in\BH^+(L_\sigma)}
 \Psi_\sigma(h_{0,\sigma},h)>-\infty.
\]
The same infimum is obtained over $\CH^+(L_\sigma)$.
If $f$ is smooth, it also equals the infimum over all smooth
fiberwise-semipositive metrics. For the fixed polarized spread,
the exceptional set is independent of the metrics in these classes
and of the integrable reference metric.
\end{thm}
\begin{proof}
We use the fixed polarized spread over a smooth geometrically
integral quasi-projective $K$-model $S$ of $F$ carrying $L$.
After shrinking $S$, we may assume that $f:X\to C$ is flat
and projective, that $g:X\to S$ and $\pi:C\to S$ are smooth
and projective with geometrically integral fibers $X_s,C_s$,
and that $L$ is $f$-ample. We retain the arithmetic
compactification and polarization $\bar A$ in
Theorem~\ref{thm:ae-energy}.

Since the geometric generic restriction of $L$ is nef, we may
apply Theorem~\ref{thm:ae-energy} at each complex embedding of
$K$. It gives boundedness outside one set of volume zero, with
the same infimum over $\BH^+(L_s)$ and $\CH^+(L_s)$. By
the refinement following that theorem, if $f$ is smooth, we shrink
the spread to make $f$ smooth there as well and obtain equality
with the infimum over smooth
fiberwise-semipositive metrics. The finite extension of constants
used at real places preserves the assertion, since restriction
of embeddings is a finite covering on the generic complex
realizations. We identify these generic points with embeddings
of $F$ to obtain the theorem for the chosen reference metric.
For each embedding outside the exceptional set, the cocycle
identity changes the energy by a finite constant when we change
the integrable reference metric. Thus the exceptional set is
independent of the reference; the lower bound may depend on it.
\end{proof}

By the change-of-metric formula, Theorem~\ref{thm:intrinsic-energy}
also gives a finite envelope norm on the corrected Deligne line
for almost every embedding. Boundedness at every embedding and
regularity of these norms remain the assertions of
Conjectures~\ref{conj:H-energy-boundedness} and~\ref{conj:envelope}.

\subsection{Non-archimedean extensions}\label{ss:nonarchimedean}

We now prove Theorem~\ref{thm:main-energy} over a nonarchimedean
local field $E$. As in Section~\ref{ss:complex-points}, we use
dynamics on projective space to construct a measure on the
analytification of a projective compactification of $S$. Its
restriction to $S(E)$ dominates the Haar measure class. To glue
local metrics, we use maxima of semipositive weights and the
invariance of the corrected energy under constant rescaling.

\subsubsection*{Finite minima of energies}

Let $v$ be a nonarchimedean place of a global field $K$, and put
$E=K_v$. Fix a projective flat model tower
$\mathcal X\to\mathcal C\to\mathcal S$ as in
Lemma~\ref{lem:polarized-arithmetic-spread}, with a relatively ample
rational model bundle $\mathcal L$. Write
$X^c\xrightarrow f C^c\xrightarrow\pi S^c$ for the projective generic
tower, with $S^c$ integral, $g=\pi f$, and $W=(S^c)_E^{\an}$.
The relative dimensions are $d-1$ and $1$. On the open $S\subset S^c$
the maps $g$ and $\pi$ are smooth with geometrically integral fibers.
We fix a strongly relatively nef adelic reference $\bar L_0$ as in
the lemma and denote its metric at $v$ by $h_0$. We also fix the
normalized polarization $\bar A$ induced by a nef rational model
bundle $\bar A_0$ with ample underlying bundle, so that $A_s^d=1$.
Here $\BH^+(L_E)$ denotes the continuous integrable relatively nef
metrics on $X^c_E/C^c_E$.

For $h\in\BH^+(L_E)$, put
\[
 E_h(s)=\Psi_s(h_{0,s},h|_{X_s})\qquad(s\in W).
\]
By \eqref{eq:H-family-pairing}, $E_h=\Psi_A(h_0,h)$ is the difference
between the weights of $D_A(h)$ and $D_A(h_0)$. It is continuous on all of $W$. On a nonreduced fiber we use the multiplicities
of its fundamental cycle; at a non-rigid point $s$ we take the
fiber over the completed residue field $\kappa(s)$.
For every $s\in S(E)$, Lemmas~\ref{lem:nonarch-metric-extension}
and~\ref{lem:relative-nef-approximation} imply
\begin{equation}\label{eq:nonarch-local-infimum}
 \inf_{\substack{h\in\BH^+(L_E)\\
                   h\ \mathrm{strongly\ rel.\ nef}}}E_h(s)
   =\inf_{h_s\in\BH^+(L_s)}\Psi_s(h_{0,s},h_s).
\end{equation}
Indeed, relatively ample model metrics extend by
Lemma~\ref{lem:nonarch-metric-extension}, strong relative nefness
is preserved by restriction, and
Lemma~\ref{lem:relative-nef-approximation} gives the full fiber class.

\begin{lemma}
\label{lem:nonarch-energy-selection}
Let $\mu$ be a Radon probability measure on $W$, and let
$h_0,h_1,\ldots,h_k\in\BH^+(L_E)$ be strongly relatively nef.
There is a constant $c\ge0$, depending only on this finite
collection and the fixed family data, such that for every
$\delta>0$ there is a strongly relatively nef
$h_\delta\in\BH^+(L_E)$ satisfying
\[
 \bigl\|E_{h_\delta}-\min(0,E_{h_1},\ldots,E_{h_k})
       \bigr\|_{L^1(\mu)}\le c\delta.
\]
Consequently, if $\int E_h\,d\mu\ge-B_0$ for every strongly
relatively nef $h\in\BH^+(L_E)$, then
\[
 e=\inf_{\substack{h\in\BH^+(L_E)\\
                   h\ \mathrm{strongly\ rel.\ nef}}}E_h\in L^1(\mu),
 \qquad \int e\,d\mu\ge-B_0.
\]
In particular $e> -\infty$ outside one $\mu$-null set.
\end{lemma}
\begin{proof}
\smallskip\noindent\emph{Choosing the metric.}
We choose one semipositive metric $n$ on an ample bundle $N$
on $C_E^c$ such that every $a_j=h_j+f_E^*n$ is semipositive.
After enlarging the twist, we may assume that $M=L_E+f_E^*N$
is ample. We fix $n$ and put $U=f_E^*N$; these fixed bundles
will also be used to bound the energy.

Let $\phi_j$ be the negative logarithmic weight of $h_j$, and
put $c_0=\max_j\|\phi_j-\phi_0\|_\infty$ and $r=c_0+1$.
We partition $W$ by the smallest index attaining
$\min_jE_{h_j}$ and choose disjoint compact subsets $K_j$
of these Borel regions whose union has mass at least $1-\delta$.
To select $h_j$ over $K_j$, we use Urysohn's lemma on the
compact Hausdorff Berkovich space $W$ to choose continuous functions
\[
 u_j:W\longrightarrow[-r,r],\qquad
 u_j|_{K_j}=r,\qquad u_j|_{K_i}=-r\quad(i\ne j).
\]
By \cite[Theorem 5.18]{BE21}, we can approximate $u_j$ uniformly
within $1/2$ by a rational model function $\chi_j$, with all
models over $\CO_E$. Then $|\chi_j|\le r+1/2$, and for
$s\in K_j$ and $i\ne j$ we have
\[
 \chi_j(s)-\chi_i(s)\ge 2r-1=2c_0+1>2c_0.
\]
We put
\[
 \phi_\delta=\max_{0\le j\le k}
                  (\phi_j+g^*\chi_j).
\]
Let $h_\delta$ be the metric with weight $\phi_\delta$,
so $|\ell|_{h_\delta}=e^{-\phi_\delta}$ in a local frame.
Since $|\phi_j-\phi_i|\le2c_0$,
the $j$-th term attains the maximum strictly on every fiber
over $K_j$. Hence, for $s\in K_j$,
\[
 \phi_\delta|_{X_s}=\phi_j|_{X_s}+\chi_j(s).
\]

To prove strong relative nefness of $h_\delta$, we choose a
common projective model of $S_E^c$ over $\CO_E$ and an ample
model metric $q_Q$ on a line bundle $Q$ such that $q_Q$ and
all $\chi_j+q_Q$ are rational Fubini--Study metrics. Then
$a_j+g^*(\chi_j+q_Q)$ are continuous semipositive metrics on
the same ample bundle $M+g^*Q$. Their maximum is semipositive
by \cite[Corollary 7.5 and Theorem 7.8]{BE21}. Indeed, we can
approximate them uniformly by Fubini--Study metrics; after taking
a common power, the union of their section tuples defines the
maximum. Subtracting $f_E^*(n+\pi^*q_Q)$ from this maximum
gives $\phi_\delta$. Since $n+\pi^*q_Q$ is semipositive on
the ample bundle $N+\pi^*Q$, the metric $h_\delta$ is strongly
relatively nef and belongs to $\BH^+(L_E)$.

\smallskip\noindent\emph{Bounding the energy.}
For $s\in K_j$, constant-rescaling invariance and the cocycle
identity give
\[
 \begin{aligned}
 E_{h_\delta}(s)-E_{h_j}(s)
 &=\Psi_s(h_j|_{X_s},h_\delta|_{X_s})\\
 &=-\chi_j(s)\bigl((d+1)L_s^d-(d+1)L_s^d\bigr)=0.
 \end{aligned}
\]
By the choice of $K_j$, we have $E_{h_j}(s)=\min_iE_{h_i}(s)$
there. Hence $E_{h_\delta}=\min_iE_{h_i}$ on $\bigcup_jK_j$.

We now bound the energy on the complement of $\bigcup_jK_j$
using the original twist $n$. On each fiber over $s\in W$,
the numbers $\chi_j(s)$ are constants. Hence
$h_\delta+f_s^*n_s$ is semipositive, and
$\|\phi_\delta-\phi_0\|_\infty\le c_*:=c_0+r+1/2$.
We expand every mixed product of $c_1(L_s,h_{0,s})$ and
$c_1(L_s,h_\delta)$ as a difference of positive mixed measures
on $M_s$ and $U_s$. Its total variation is then at most
$(M_s+U_s)^d$. We put
\[
 R_0=\max_{s\in W}(M_s+U_s)^d
     =\max_{s\in W}(L_s+2f_s^*N_s)^d.
\]
By flatness, the intersection numbers in this maximum are
constant on each of the finitely many connected components of
$(S^c)_E$, so $R_0$ is finite. The mass formula and the
corrected energy identity imply
\[
 |E_{h_\delta}(s)|\le(d+1)c_*\bigl(R_0+|L_s^d|\bigr).
\]
The original finite minimum has the same bound. The bound
depends on the fixed metrics and the original twist $n$, and
is independent of $\delta$, of $q_Q$, and of the chosen compact subsets.
It holds at every Berkovich point, including non-rigid points
and points over finite extensions charged by $\mu$. The
$L^1$ error is therefore at most
$2(d+1)c_*\bigl(R_0+|L_s^d|\bigr)\delta$.

The assumed integral bound, followed by $\delta\to0$, bounds
the integral of every finite minimum below by $-B_0$.
Since $E$ has a countable dense subset, $W$ is second countable
\cite[Theorem 1.1]{HLP14}. The countable-subcover argument and
monotone convergence in Lemma~\ref{lem:energy-envelope-gluing}
therefore pass this bound to the full infimum.
\end{proof}

\subsubsection*{Embeddings into nonarchimedean local fields}

Let $K_0=\BQ$ or $\mathbb F_p(t)$, and let $E$ be a finite
extension of $\BQ_p$ or $\mathbb F_p((t))$, respectively,
with the standard embedding $K_0\hookrightarrow E$.
For a finitely generated extension $F/K_0$, let $K$ be the
algebraic closure of $K_0$ in $F$, and fix
$\iota:K\hookrightarrow E$. Then $K$ is a global field and
$F/K$ is regular. In positive characteristic, regularity follows
from the fact that $K$ has $p$-rank one: since $K$ is algebraically
closed in $F$, a $p$-basis element of $K$ has no $p$-th root in $F$.

We define a measure class on these embeddings by choosing a
smooth geometrically integral $K$-model $S$ with $K(S)=F$.
Haar measure in analytic coordinates defines a measure class
on $S_\iota(E)$, since coordinate changes preserve null sets.
The proper $K$-subvarieties are countable in number and each
has measure zero. Their complement consists of the $K$-generic
points, which identify with embeddings $F\hookrightarrow E$
extending $\iota$. Comparing two models on a common birational
open subset shows that the measure class is independent of $S$.
We use this class on each nonempty component corresponding to
an embedding $K\hookrightarrow E$.

For the following theorem, let $K_0=\BQ$ or $\mathbb F_p(t)$,
let $F/K_0$ be finitely
generated, and let $E$ be a fixed finite extension of $\BQ_p$
or $\mathbb F_p((t))$, respectively. Let $f:X\to C$ be a flat
projective surjective morphism of smooth geometrically integral
projective varieties over $F$, with $\dim X=d$ and $\dim C=1$.
Fix an admissible collection $(\bar A_\sigma)$ and polarized
model data inducing it, as in
Definition~\ref{def:admissible-polarizations}.
The Haar measure classes just defined give the meaning of
almost every embedding. In the metric classes, the auxiliary
base bundle may depend on the metric. For $L^d=0$, the
energy is independent of $\bar A$.

\begin{thm}
\label{thm:padic-intrinsic-energy}
Let $L$ be a geometrically nef and $f$-ample rational
line bundle.
For almost every $K_0$-embedding $\sigma:F\hookrightarrow E$
and every integrable reference metric $h_{0,\sigma}$,
\[
 \inf_{h\in\BH^+(L_\sigma)}
 \Psi_\sigma(h_{0,\sigma},h)>-\infty.
\]
The same infimum is obtained over $\CH^+(L_\sigma)$.
For the fixed polarized spread, the exceptional set is independent
of the metric in $\BH^+(L_\sigma)$ and of the integrable
reference metric.
\end{thm}

The proof uses the following nonarchimedean version of
Lemma~\ref{lem:bounded-height-volume}. Let $K$ be a global
field, $v$ a nonarchimedean place, and $E=K_v$. Let $S$ be
a smooth geometrically integral quasi-projective $K$-variety
with a normal projective compactification $S^c$, and fix an
ample height $h_B$ on $S^c$. We write $W=(S^c)_E^{\an}$.
For $a\in S^c(\bar K)$ with separable residue field $F_a=K(a)$,
let $a_w$ be the Berkovich point induced by a place $w\mid v$
of $F_a$. We use the normalized orbit measure
\[
 \nu_a=\nu_{a,v}
 =\frac1{[F_a:K]}\sum_{w\mid v}[F_{a,w}:K_v]\,\delta_{a_w}.
\]
The measure $\nu_a$ is the pushforward of the average over all
embeddings $F_a\hookrightarrow\overline{K_v}$ extending
$K\hookrightarrow K_v$. We say that $\mu|_{S(E)}$ dominates
the Haar measure class if every Borel subset of $S(E)$ of
$\mu$-measure zero has Haar measure zero.

\begin{lem}
\label{lem:nonarch-bounded-height-measure}
There are $H>0$ and a probability measure $\mu$ on $W$ such that
$\mu|_{S(E)}$ dominates the Haar measure class and $\mu$ gives zero
mass to every proper closed $K$-subvariety of $S^c$.
For every dense $K$-open subset $U\subset S$, the same $\mu$ is a
weak limit of convex combinations of normalized Galois orbit
measures of points $a\in U(\bar K)$ with separable residue field
over $K$ and $h_B(a)\le H$.
\end{lem}
\begin{proof}
\smallskip\noindent\emph{A measure on projective space.}
We use the construction on symmetric products from
Lemma~\ref{lem:bounded-height-volume}, with $n=\dim S$.
When $n=0$, the point mass and the identity projection prove
the result. Suppose $n>0$. We choose $\varpi\in K$ that is a
uniformizer at $v$, put $q_v=\#(\CO_E/\varpi)$, and define
\[
 \varphi:\BP^1\longrightarrow\BP^1,\qquad
 \varphi(z)=\frac{z^{q_v}-z}{\varpi}.
\]
On each residue ball of $\CO_E$, the map $\varphi$ is an
analytic bijection onto $\CO_E$, expanding distances by
$|\varpi|^{-1}$. Thus $R_N=(\varphi^N)^*[0]$ consists of
$q_v^N$ simple points, one in each residue class modulo
$\varpi^N$. The associated probability measures converge to
Haar probability $m$ on $\CO_E$. This also holds in
characteristic $p$, where $\varphi'(z)=-\varpi^{-1}$.
For $E=\BQ_p$, see \cite[Theorem 1 and Lemma 3]{PS19}.

We pass to $\BP^n$ through the quotient
$\sigma:(\BP^1)^n\to\operatorname{Sym}^n(\BP^1)\simeq\BP^n$.
The product $\varphi^{\times n}$ descends to an endomorphism
$\Phi$ of algebraic degree $q_v$, giving
\[
 \begin{tikzcd}[column sep=large,row sep=large]
 (\BP^1)^n \arrow[r,"\varphi^{\times n}"] \arrow[d,"\sigma"']
   & (\BP^1)^n \arrow[d,"\sigma"]\\
 \BP^n \arrow[r,"\Phi"'] & \BP^n,
 \end{tikzcd}
 \qquad
 \Phi^*\CO_{\BP^n}(1)\simeq\CO_{\BP^n}(q_v).
\]
For the fixed point $b=\sigma(0,\ldots,0)$, we put
\[
 D_N=(\Phi^N)^*[b]=\sigma_*(R_N^{\times n}),\qquad
 \deg D_N=q_v^{nN}.
\]
The identity of cycles follows from
$\sigma^*[b]=n![(0,\ldots,0)]$ by finite flat pullback and
pushforward, in every characteristic. Consequently the measures
associated with $q_v^{-nN}D_N$ converge to
$\lambda_0=\sigma_*m^{\otimes n}$. To bound the heights of
their support points, we again use the Call--Silverman canonical
height \cite{CS93}, as in Lemma~\ref{lem:bounded-height-volume}:
\[
 \widehat h_\Phi(y)=\lim_{j\to\infty}q_v^{-j}h(\Phi^j(y))
\]
It differs from the Weil height $h$ by a bounded function and
vanishes on preperiodic points. Over a global function field,
the same limit converges because $h\circ\Phi=q_vh+O(1)$.
Every point in the support of $D_N$ is preperiodic, so the
heights of all these points are uniformly bounded.

We next spread the limiting measure over $\BP^n(E)$ by
projective coordinate changes. Off the diagonals, the Jacobian
of $\sigma$ is the nonzero Vandermonde product
$\prod_{i<j}(z_i-z_j)$. Since the diagonals have product Haar
measure zero, $\lambda_0$ is absolutely continuous. We choose
disjoint small balls around $n$ distinct elements of $\CO_E$.
On their product, $\sigma$ is an analytic isomorphism onto a
nonempty compact open set $\Omega\subset\BP^n(E)$, and
$\lambda_0$ dominates the Haar measure class on $\Omega$.
Using density of $\operatorname{PGL}_{n+1}(K)$ in
$\operatorname{PGL}_{n+1}(E)$ and compactness, we choose
$\gamma_1,\ldots,\gamma_\ell\in\operatorname{PGL}_{n+1}(K)$
whose translates of $\Omega$ cover $\BP^n(E)$. The measure
\[
 \lambda=\frac1\ell\sum_{j=1}^{\ell}(\gamma_j)_*\lambda_0
\]
is then in the Haar measure class on $\BP^n(E)$. It is the
weak limit of the probability measures of
$(\ell q_v^{nN})^{-1}\sum_j(\gamma_j)_*D_N$. Every point
in the support of these cycles has uniformly bounded height,
since $h\circ\gamma_j=h+O(1)$ for each of the finitely many $j$.

\smallskip\noindent\emph{Pullback to the parameter space.}
As in Lemma~\ref{lem:bounded-height-volume}, a linear center
avoiding $S^c$ and the tangent space at a geometric smooth point
gives a finite generically \'etale $K$-projection
$u:S^c\to\BP^n$ of degree $\delta$.
We choose a proper closed subset
$Z\subset\BP^n$ containing the image of the non-\'etale
locus and $u(S^c\setminus S)$, and put $V=\BP^n\setminus Z$
and $U=u^{-1}(V)\subset S$.

\Needspace{10\baselineskip}
We have the Cartesian square
\[
 \begin{tikzcd}[column sep=large,row sep=large]
 U \arrow[r,hook] \arrow[d,"u|_U"']
   & S^c \arrow[d,"u"]\\
 V \arrow[r,hook] & \BP^n
 \end{tikzcd}
\]
whose left vertical map is finite \'etale of degree $\delta$.

Since $\lambda(Z(E))=0$, we may delete components in $Z$ from
the target cycles and renormalize without changing their weak
limit or the height bound on their points. We pull the remaining
cycles back through $U\to V$ and divide by $\delta$.
Height comparison for the finite map $u$ bounds the height of
every point in these pullbacks. The points of the root cycles
have separable residue fields over $K$. Products, pushforward
by $\sigma$ and the $\gamma_j$, and finite \'etale pullback
preserve separability. Hence the resulting measures are convex
combinations of the normalized orbit measures $\nu_a$.

To prove convergence on $(S^c)_E^{\an}$, let $\psi$ be a
continuous function on this compact Berkovich space. We take
its trace on $V(E)$ with residue-field degrees and extend it
by zero on $Z(E)$. The trace is bounded by
$\delta\|\psi\|_\infty$ and is continuous off the
$\lambda$-null set $Z(E)$. Its integrals therefore converge
under the weak convergence on $\BP^n(E)$. The limiting
probability measure $\mu$ is the normalized trace pullback
of $\lambda$ on $V$.

We check domination on $S(E)$ locally. For $x\in U(E)$,
the inverse function theorem gives compact open neighborhoods
$V_x\subset U(E)$ and $W_x\subset V(E)$ on which
$u:V_x\to W_x$ is an analytic isomorphism. For a Borel
subset $W\subset V_x$, we have
\[
 \mu(W)=\frac{\lambda(u(W))}{\delta}.
\]
Thus $\mu|_{S(E)}$ dominates the Haar measure class, because
$S\setminus U$ is a proper $K$-subvariety and has Haar
measure zero. We regard $\mu$ as a measure on the full
Berkovich space; it may also charge points over finite
extensions of $E$. By finiteness of $u$ and the corresponding
property of $\lambda$, every proper $K$-subvariety has
$\mu$-measure zero. We can therefore delete components outside
any prescribed dense open subset and obtain the same limit,
with the same height bound on every remaining point.

\end{proof}

\begin{proof}[Proof of Theorem~\ref{thm:padic-intrinsic-energy}]
\smallskip\noindent\emph{The spread and the local field.}
We use the fixed polarized spread over a smooth geometrically
integral quasi-projective $K$-model $S$ of $F$. After shrinking
$S$, we may assume that $f$ is flat, that $X\to S$ and
$C\to S$ are smooth with geometrically integral fibers, and
that $L$ is $f$-ample. We fix an embedding $K\hookrightarrow E$
and denote the induced place by $v$.

We may assume $K_v=E$ after a finite extension of $K$ inside
$E$. For a separable local extension, we use polynomial
approximation and Krasner's lemma. In positive characteristic,
we first realize the maximal separable subextension $E_s/K_v$
in this way and choose a global uniformizer $u$ of $E_s$.
If $[E:E_s]=p^a$, then $E=E_s(u^{1/p^a})$, since $E_s$
has $p$-rank one. Adjoining $u^{1/p^a}$ globally therefore
realizes the remaining purely inseparable extension. For the
chosen inclusion of the enlarged field $K'$ into $E$, we have
\[
 (S\times_K K')\times_{K'}E\simeq S\times_K E
\]
Thus the analytic parameter space and its Haar measure class
are preserved. By regularity of $F/K$, the compositum is a
field, and this change of constants preserves the embeddings
under consideration.

Over either global field, Lemma~\ref{lem:polarized-arithmetic-spread}
gives a projective flat tower carrying a relatively ample model
of $L$, with normal projective parameter base $S^c$.

The relatively ample model defines the reference metrics, which
will stay fixed away from $v$. We use the common tower constructed
after Definition~\ref{def:admissible-polarizations}, so the pullback
of $\bar A_0$ induces the prescribed polarization with $A_s^d=1$.
Let
$X^c\xrightarrow{f^c}C^c\xrightarrow{\pi^c}S^c$ be the
generic projective tower. For $s\in S(E)$, we have
\[
 \begin{tikzcd}[column sep=large,row sep=large]
 X_s \arrow[r] \arrow[d,"f_s"']
   & X_E^c \arrow[r] \arrow[d,"f_E"]
   & X^c \arrow[d,"f^c"]\\
 C_s \arrow[r] \arrow[d]
   & C_E^c \arrow[r] \arrow[d,"\pi_E"]
   & C^c \arrow[d,"\pi^c"]\\
 \Spec E \arrow[r,"s"'] & S_E^c \arrow[r] & S^c.
 \end{tikzcd}
\]
All four squares are Cartesian.

Choose the measure $\mu$ of
Lemma~\ref{lem:nonarch-bounded-height-measure} on $(S^c)_E^{\an}$.

By \eqref{eq:nonarch-local-infimum} and
Lemma~\ref{lem:nonarch-energy-selection}, it suffices to bound
the average energy on this full Berkovich space for strongly
relatively nef $h\in\BH^+(L_E)$. Such metrics can replace
the component at $v$ by Corollary~\ref{cor:adelic-replacement}.

\smallskip\noindent\emph{The average energy bound.}
We use the Haar-normalized absolute value on $E=K_v$, with
$|\varpi|_E=q_v^{-1}$ for a uniformizer $\varpi$ and residue
cardinality $q_v$. Absolute point heights have global denominator
$[K:K_0]$. Over a function field, the unnormalized adelic
intersection degrees are geometric degrees over $\mathbb F_p$
multiplied by $\log p$; to obtain a point height, we divide by
the degree of its field of definition over $K_0$.

For a strongly relatively nef $h\in\BH^+(L_E)$ and a rational
$\epsilon>0$, we put
\[
 \bar L_\epsilon(h)=\bar L_0^v(h)+\epsilon\bar A_0,
 \qquad m_\epsilon=(L_s+\epsilon A_{0,s})^d.
\]
Since $L$ is nef on the generic parameter fiber,
$L+\epsilon A_0$ is ample there. It is therefore ample on
every fiber over a dense open subset of $S$. We define
\[
 D_A(\bar L_\epsilon(h))=-D_g(\bar L_\epsilon(h))
 +(d+1)m_\epsilon\langle\bar L_\epsilon(h),\bar A^{[d]}\rangle_g.
\]
This is a rational adelic line bundle. On the ample open set,
Theorem~\ref{thm:ahi}, with the normalization above, gives
\[
 h_{D_A(\bar L_\epsilon(h))}(a)\ge0.
\]

Since all metrics away from $v$ are fixed, the variation of
the metric on $D_A(L_\epsilon)$ is given by the continuous local function
\[
 E_{\epsilon,h}(s)
 =\Psi_s(h_{0,s}+\epsilon h_{A_0,s},
                    h|_{X_s}+\epsilon h_{A_0,s}),
\]
where $h_{A_0,s}$ is the metric of $\bar A_{0,s}$. The correction
uses $(d+1)m_\epsilon c_1(\bar A_s)^d$, with the normalized
polarization held fixed.
The height comparison with the fixed reference gives
\[
 \frac1{[K:K_0]}\int E_{\epsilon,h}\,d\nu_a
 =h_{D_A(\bar L_\epsilon(h))}(a)
  -h_{D_A(\bar L_\epsilon(h_0))}(a)
 \ge -B_*.
\]
We can choose $B_*$ independently of $h$ and of rational
$0<\epsilon\le1$ for every point in the approximating cycles.
Indeed, multilinearity expresses $D_A(\bar L_\epsilon(h_0))$
as a polynomial in $\epsilon$ with finitely many fixed mixed
Deligne bundles as coefficients. The heights for each coefficient
bundle are bounded on points satisfying the fixed ample height
bound. For each $\epsilon$, we delete components outside its
ample open set; Lemma~\ref{lem:nonarch-bounded-height-measure}
preserves both the limit $\mu$ and the height bound.

For fixed $h$ and $\epsilon$, we now pass to the weak limit
of the convex combinations of orbit measures. We then let
$\epsilon\to0$. By multilinearity of the Deligne pairings,
$E_{\epsilon,h}\to E_h$ uniformly on the compact analytic
parameter space. Consequently $\int E_h\,d\mu\ge-B_0$,
where $B_0=[K:K_0]B_*$ depends on the fixed model, reference,
polarization, and height bound, and is independent of the strongly
relatively nef metric $h\in\BH^+(L_E)$.

By Lemma~\ref{lem:nonarch-energy-selection}, the infimum $e$
belongs to $L^1(\mu)$. Domination on $S(E)$ and
\eqref{eq:nonarch-local-infimum} imply finiteness outside one
Haar-null set for the full fiber class $\BH^+(L_s)$.
Identify generic $E$-points with embeddings of $F$ and apply
the cocycle identity to change the reference metric.
\end{proof}

\begin{remark}
The Haar measure class in
Theorem~\ref{thm:padic-intrinsic-energy} uses local compactness
of $E$. For $\BC_p$, the additive group has no nonzero locally
finite Radon Haar measure, since $\BC_p$ is not locally compact.
The compact Berkovich analytification of a projective variety
does carry Chambert--Loir measures, but their support can consist
of nonclassical points. For example, the standard model metric
on $\CO_{\BP^1}(1)$ gives the Dirac mass at the Gauss point;
see \cite[Proposition 2.1.1]{CL10}. To formulate a result for
embeddings $F\hookrightarrow\BC_p$, one must therefore specify
a measure and prove boundedness for that measure. The finite
extensions $E/\BQ_p$ account for embeddings with image in
such fields; completion also gives elements transcendental over
$\BQ_p$, and hence further embeddings in general.

\end{remark}

\clearpage
\phantomsection

\end{document}